\documentclass[a4paper,reqno, 11pt]{amsart}  
\usepackage[DIV=12, oneside]{typearea}
\usepackage[english]{babel}
\usepackage[centertags]{amsmath}
\usepackage{amstext,amssymb,amsopn,amsthm}
\usepackage{mathrsfs}
\usepackage{dsfont}
\usepackage{bbm}
\usepackage{thmtools}
\usepackage{graphicx}
\usepackage[titletoc,title]{appendix}
\usepackage{etexcmds}
\usepackage{esint}
\usepackage{multicol}

\usepackage[backgroundcolor=white, bordercolor=blue,
linecolor=blue]{todonotes}

\usepackage[colorlinks=true, linkcolor=black, citecolor=black]{hyperref}
\usepackage{enumitem}
\setlist[enumerate]{itemsep=0mm}

\theoremstyle{plain}
\declaretheorem[title=Theorem, parent=section]{theorem}
\declaretheorem[title=Lemma,sibling=theorem]{lemma}
\declaretheorem[title=Proposition,sibling=theorem]{proposition}
\declaretheorem[title=Corollary,sibling=theorem]{corollary}

\theoremstyle{definition}
\declaretheorem[title=Definition,sibling=theorem]{definition}
\declaretheorem[title=Remark,sibling=theorem]{remark}
\declaretheorem[title=Remark, numbered=no]{remark*}

\declaretheorem[title=Assumption, numbered=no]{assumption*}

\numberwithin{equation}{section}

\newcommand{\N}{\mathds{N}}
\newcommand{\R}{\mathds{R}}

\newcommand{\1}{\mathbbm{1}}

\newcommand{\eps}{\varepsilon}

\DeclareMathOperator{\dist}{dist}

\DeclareMathOperator{\supp}{supp}

\DeclareMathOperator*{\osc}{osc}

\renewcommand{\d}{\textnormal{d}}

\newcommand{\U}{\tilde{u}}

\newcommand{\F}{\tilde{f}}

\def\dashint{\fint}

\begin{document}

\title{Semiconvexity estimates for nonlinear integro-differential equations}

\author{Xavier Ros-Oton}
\author{Clara Torres-Latorre}
\author{Marvin Weidner}

\address{ICREA, Pg. Llu\'is Companys 23, 08010 Barcelona, Spain \& Universitat de Barcelona, Departament de Matem\`atiques i Inform\`atica, Gran Via de les Corts Catalanes 585, 08007 Barcelona, Spain \& Centre de Recerca Matem\`atica, Barcelona, Spain}
\email{xros@icrea.cat}

\address{Departament de Matemàtiques i Informàtica, Universitat de Barcelona, Gran Via de les Corts Catalanes 585, 08007 Barcelona, Spain}
\email{claratorreslatorre@ub.edu}

\address{Departament de Matemàtiques i Informàtica, Universitat de Barcelona, Gran Via de les Corts Catalanes 585, 08007 Barcelona, Spain}
\email{mweidner@ub.edu}





\begin{abstract}
In this paper we establish for the first time local semiconvexity estimates for fully nonlinear equations and for obstacle problems driven by integro-differential operators with general kernels. 
Our proof is based on the Bernstein technique, which we develop for a natural class of nonlocal operators and consider to be of independent interest. 
In particular, we solve an open problem from Cabr\'e-Dipierro-Valdinoci \cite{CDV22}.
As an application of our result, we establish optimal regularity estimates and smoothness of the free boundary near regular points for the nonlocal obstacle problem on domains.
Finally, we also extend the Bernstein technique to parabolic equations and nonsymmetric operators.
\end{abstract}

\maketitle

\section{Introduction} 

The aim of this work is to establish semiconvexity estimates for solutions to nonlinear equations driven by integro-differential operators of the form
\begin{equation}\label{eq:L_intro-Bernstein}
Lu(x) = \text{p.v.}\int_{\R^n} \big(u(x) - u(y)\big) K(x-y) \d y,
\end{equation}
where $K : \R^n \to [0,\infty]$ is comparable to the kernel of the fractional Laplacian. To be precise, we consider the following natural class of symmetric jumping kernels $K$ (see \cite{CaSi11a}, \cite{CaSi11b}) satisfying the classical uniform ellipticity condition
\begin{align}
\tag{$K_{\asymp}$}
\lambda |y|^{-n-2s} \le K(y) \le \Lambda |y|^{-n-2s}
\end{align}
for some constants $0 < \lambda \le \Lambda$ and $s \in (0,1)$, and the smoothness conditions
\begin{align}
\tag{$C^1$}
|\nabla K(y)| &\le \Lambda |y|^{-1} K(y),\\
\tag{$C^2$}
|D^2 K(y)| &\le \Lambda |y|^{-2} K(y).
\end{align}
We denote the family of all such operators by $\mathcal{L}_s(\lambda,\Lambda;2)$; see also \autoref{def:classL}.

\subsection{Fully nonlinear equations}

The regularity theory for  fully nonlinear nonlocal equations was developed by Caffarelli and Silvestre in their celebrated series of papers \cite{CaSi09}, \cite{CaSi11a}, and \cite{CaSi11b}.
They established a nonlocal counterpart of the Krylov-Safonov theorem, stating that solutions are $C^{1+\eps}$, and an Evans-Krylov theorem, which yields $C^{2s+\eps}$-regularity of solutions to concave equations.
Let us also refer to \cite{ChDa12}, \cite{GuSc12}, \cite{Kri13}, \cite{Ser15a}, \cite{Ser15b}, \cite{ChDa16}, \cite{ScSi16}, \cite{JiXi16}, \cite{ChKr17}, \cite{DoZh19} for extensions of the aforementioned results, e.g., to operators with coefficients and to parabolic problems; see also \cite{FeRo24} and \cite{BCI08,BCI11,BaIm08}.

It remains an intriguing open problem after \cite{CaSi11b} to establish higher regularity of solutions to nonlocal Bellman-type equations. 
So far, it is still unknown whether solutions are more regular than $C^{1 + \eps}\cup C^{2s + \eps}$, even if the underlying class of operators possesses only smooth kernels.

Here we prove for the first time semiconvexity estimates for solutions to these equations:

\begin{theorem}
\label{thm:intro-FN-sc}
Let $s\in(0,1)$, and let $u$ be any viscosity solution to a fully nonlinear equation 
\begin{equation}\label{eq:intro-fully-nonlinear}
\inf_{\gamma \in \Gamma} \{ L_\gamma u  \} = 0 ~~ \text{ in } B_1, 
\end{equation}
where $\{L_\gamma\}_{\gamma\in\Gamma} \subset \mathcal{L}_s(\lambda,\Lambda;2)$. 
Then, $u$ satisfies in the convexity sense\footnote{Here, and throughout the entire article we say that $\partial_{ee}^2 u \ge -C$ holds true ``in the convexity sense'', if $x \mapsto u(x) + \frac{C}{2n} |x|^2$ is convex.}
\begin{align*}
\partial_{ee}^2 u \ge -C\Vert u \Vert_{L^{\infty}(\R^n)} ~~ \text{ in } B_{1/2}
\end{align*}
for all $e\in \mathbb S^{n-1}$, where $C$ depends only on $n$, $s$, $\lambda$, $\Lambda$.
\end{theorem}

Our \autoref{thm:intro-FN-sc} establishes one-sided $C^{1,1}$-regularity estimates and can be seen as the first contribution after \cite{CaSi11b} to the higher regularity for solutions to \eqref{eq:intro-fully-nonlinear}. Note that the class of operators $\mathcal{L}_s(\lambda,\Lambda;2)$ is also considered in \cite{CaSi11b}.
A more general version of the result, where we allow for an $x$-dependent right hand side, will be proved in \autoref{thm:fully-nonlinear_one-sided}.

Semiconvexity estimates play a crucial role in the study of nonlinear elliptic PDE.
For instance, semiconvexity estimates for solutions to second order fully nonlinear PDE imply two-sided $C^{1,1}$ regularity estimates, if the operator under consideration is uniformly elliptic (see \cite{CaCa95}).
In light of this, it is an interesting question to ask whether \autoref{thm:intro-FN-sc} also implies $C^{1,1}$-regularity estimates for \eqref{eq:intro-fully-nonlinear}\footnote{Even in the local case, optimal regularity for Bellman equations remains an open problem; see e.g. \cite[Chapter 4.5]{FeRo22}. The only partial result in this direction is due to Caffarelli, De Silva, and Savin \cite{CDS18}, who established optimal $C^{2,1}$ regularity in case $n = 2$ and for equations of the type $\min\{L_1 u,\,L_2u\}=0$.}, at least when $s \ge 1/2$.

Let us emphasize that, for nonlocal fully nonlinear equations, the only known semiconvexity estimate was that of Cabr\'e, Dipierro, and Valdinoci \cite{CDV22}, who proved (a priori) semiconvexity estimates for fully nonlinear equations built from affine transformations of the fractional Laplacian.

\subsection{Obstacle problems}

Semiconvexity estimates play a crucial role in the study of obstacle problems
\[
\label{eq:intro-OP}
\min\{L u , u - \phi\} = 0 ~~ \text{ in } \Omega \subset \R^n.
\]
While the classical obstacle problem (corresponding to $L = -\Delta$) is very well understood \cite{Caf77,PSU12,FeRo22}, the case of integro-differential operators $L$ is significantly more complicated, and several questions remain open. 
The regularity theory for nonlocal obstacle problems was initiated in the seminal works \cite{ACS08,Sil07,CSS08} for $L = (-\Delta)^s$, and further developed in \cite{GaPe09,PePo15,KPS15,JhNe17,GPP17,FoSp18,FeRo18,CSV20,CFR24,FeJa21,Kuk21,SaYu21}.
In that case, one can identify the obstacle problem for $(-\Delta)^s$ with a (weighted) thin obstacle problem in $\R^{n+1}_+$ thanks to the celebrated Caffarelli-Silvestre extension \cite{CaSi07}, and this gives access to many local techniques.
In particular, it allows to prove local semiconvexity estimates in the extended variables; see \cite{AtCa06,ACS08} for the case $s=\frac12$, and \cite{CDV22,FeJa21} for all $s\in(0,1)$.

The regularity theory for obstacle problems \eqref{eq:intro-OP} with general integro-differential operators $L$ requires quite different methods compared to the case $L=(-\Delta)^s$, and has been developed in \cite{CRS17,CDS17,AbRo20,FRS23}.
The best known results so far establish the optimal $C^{1+s}$-regularity of solutions \cite{FRS23}, as well as the smoothness of free boundaries near regular points \cite{CRS17,AbRo20}, whenever solutions are \emph{semiconvex}.
The semiconvexity property holds true for any global solution\footnote{Also, if one assumes that the contact set $\{u=\phi\}$ is compactly contained in $\Omega$, then one can use a cutoff argument to transform the problem into a global one. However, this does not give local regularity estimates, and all constants would depend on the distance from $\{u=\phi\}$ to $\partial\Omega$.} (i.e. when $\Omega=\R^n$), which follows by a simple translation argument based on the maximum principle \cite{Sil07,CRS17}.

An important open problem was to establish \emph{local} regularity estimates for solutions of \eqref{eq:intro-OP}, not relying on any a priori assumptions on the boundary data.
We solve this by proving \emph{local} semiconvexity estimates, which is exactly the content of our next result.

\begin{theorem}
\label{thm:intro-OP-sc}
Let $s\in(0,1)$, $L \in \mathcal{L}_s(\lambda,\Lambda;2)$, and  $u$ be any solution to the nonlocal obstacle problem 
\begin{equation}\label{eq:intro-OP2}
\min\{L u , u - \phi\} = 0 ~~ \text{ in } B_1,
\end{equation}
where  $\phi \in C^{2+2s+\eps}(\R^n)$ for some $\eps > 0$. 
Then, $u$ satisfies in the convexity sense
\begin{align*}
\partial_{ee}^2 u \ge -C(\Vert u \Vert_{L^{\infty}(\R^n)} + \Vert L \phi \Vert_{C^{1,1}(B_1)}) ~~ \text{ in } B_{1/2}
\end{align*}
for all $e\in \mathbb S^{n-1}$, where $C$ depends only on $n$, $s$, $\lambda$, $\Lambda$.
\end{theorem}

As said above, semiconvexity estimates are a key tool in the regularity theory for nonlocal obstacle problems, mainly because they imply the convexity of blow-ups.
Without this, one cannot establish the optimal regularity of solutions and the regularity of free boundaries.
Moreover, as we will see next, this also allows us to study obstacle problems under minimal regularity assumptions on the obstacle $\phi$ (which is new even in the global case $\Omega=\R^n$).

\subsubsection*{Regularity of solutions and free boundaries}

Let us explain now the main applications of our new semiconvexity estimate for nonlocal obstacle problems.
We consider general nonlocal operators belonging to the class $\mathcal{L}_s(\lambda,\Lambda;1)$, i.e., having symmetric jumping kernels $K$ satisfying\footnote{Note that, in contrast to global solutions, optimal regularity for local solutions to nonlocal obstacle problems cannot hold true in general without some regularity assumption on $K$, or  prescribing regularity on the complement data (see \cite[Proposition 6.1]{RoSe16b}, \cite[Section 5.1]{Ser15b}).} the ellipticity condition \eqref{eq:coercive} and \eqref{eq:L1}.
As in \cite{CRS17,AbRo20,FRS23}, we also need to assume that $K$ is homogeneous, i.e., 
\begin{align}
\label{eq:hom}
K(y) =  \frac{K(y/|y|)}{|y|^{n+2s}} \qquad \textrm{for all} \quad y \in \R^n\setminus\{0\}.
\end{align}
Our first result in this direction is to obtain for the first time \emph{local} $C^{1+s}$-estimates.

\begin{corollary}[Optimal regularity]
\label{thm:OP_opt-reg}
Let $s\in(0,1)$ and $L \in \mathcal{L}_s(\lambda,\Lambda;1)$ satisfying \eqref{eq:hom}. 
Let $u$ be any solution to the obstacle problem 
\begin{align*}
\min\{Lu , u-\phi\} = 0 ~~ \text{ in } B_1,
\end{align*}
with $\phi \in C^{\beta}(B_1)$ for some $\beta > 0$. 
Then, 
\begin{align*}
\qquad \qquad \Vert u \Vert_{C^{\beta}(B_{1/2})} \le C (\Vert \phi \Vert_{C^{\beta}(B_1)} + \Vert u \Vert_{L^{\infty}(\R^n)} )\qquad \textrm{for}\quad \beta<1+s,
\end{align*}
while
\begin{align*}
\quad\ \qquad \qquad \Vert u \Vert_{C^{1+s}(B_{1/2})} \le C (\Vert \phi \Vert_{C^{\beta}(B_1)} + \Vert u \Vert_{L^{\infty}(\R^n)} ) \qquad \textrm{for}\quad \beta>1+2s.
\end{align*}
The constant $C$ depends only on $n$, $s$, $\beta$, $\lambda$, and $\Lambda$.
\end{corollary}

The key point of this result is that the estimates are completely local, i.e., they do not depend on any assumption on the boundary data nor on the obstacle on $\partial B_1$.
Still, even in case of global solutions, observe that \autoref{thm:OP_opt-reg} only requires\footnote{In case $L = (-\Delta)^s$, it was proved in \cite[Theorem 6.1]{CDS17} that optimal $C^{1+s}$ estimates already hold when $\phi \in C^{\beta}$ with $\beta > 1+s$. 
Their proof is based on a truncated version of Almgren's monotonicity formula for the extension. 
Due to the lack of monotonicity formulas in the more general setting of nonlocal operators $L\in \mathcal{L}_s(\lambda,\Lambda;1)$, such proof cannot work in our context.
Still, we expect that the same result for $\beta>1+s$ can be established for general kernels using our new result on the convexity of blow-ups, \autoref{thm:blowups-convex}.}  obstacles $\phi \in C^{\beta}$ for some $\beta > 1+2s$ in order to obtain optimal $C^{1+s}$-regularity estimates for solutions, which even improves the results in \cite{CRS17,FRS23} for the global nonlocal obstacle problem in case $s<\frac12$. 
Furthermore, the optimal $C^{\beta}$-regularity in case of non-smooth obstacles for $\beta < 1+s$ is also new even in case of global solutions; it was only known for $\beta\leq \max\{1+\epsilon,2s+\epsilon\}$ \cite[Theorem 5.1]{CDS17}.

Our second result in this direction establishes the regularity of the free boundary near regular points for local solutions to \eqref{thm:OP_opt-reg}.

\begin{corollary}[Free boundary regularity]
\label{thm:OP_bdry-reg}
Let $s\in(0,1)$ and $L \in \mathcal{L}_s(\lambda,\Lambda;1)$ satisfying \eqref{eq:hom}. 
Let $u$ be any solution to the obstacle problem \eqref{thm:OP_opt-reg}
with  $\phi \in C^{\beta}(B_1)$ for some $\beta > 1+2s$, and let $\alpha \in (0,s) \cap (0,1-s)$. 
Then, near any free boundary point $x_0 \in \partial \{ u > \phi \}\cap B_{1/2}$, there exist $c_0 \ge 0$ and $e \in \mathbb S^{n-1}$ such that for any $x \in B_{1/4}(x_0)$:
\begin{align*}
\left| u(x) - \phi(x) - c_0 \left( (x-x_0) \cdot e \right)_+^{1+s} \right| \le C (\Vert \phi \Vert_{C^{\beta}(B_1)} + \Vert u \Vert_{L^{\infty}(\R^n)} ) |x-x_0|^{1+s+\alpha},
\end{align*}
where $C$ depends only on $n$, $s$, $\lambda$, $\Lambda$, and $\alpha$.\\
Moreover, if $c_0 > 0$, then the free boundary is a $C^{1,\alpha}$-graph in a ball $B_{\rho_0}(x_0)$ with $C \rho_0^{\alpha} \ge c_0$ and $C$ depending only on $n$, $s$, $\lambda$, $\Lambda$, and $\alpha$.
\end{corollary}

If we assume in addition that $\phi\in C^\infty(B_1)$ and $K|_{\mathbb S^{n-1}}\in C^\infty(\mathbb S^{n-1})$, then the free boundary is actually $C^{\infty}$ near regular points (i.e., near those points $x_0$ for which $c_0 > 0$); see \cite{AbRo20}.

\subsubsection*{Convexity of blow-ups}

As explained before, the semiconvexity of solutions to nonlocal obstacle problems is crucial in the study of the regularity of the free boundary and of optimal regularity estimates. 
The approach in \cite{ACS08,Sil07,CSS08,CRS17,FRS23} is to show first that solutions are semiconvex and to prove that blow-ups, as limits of correctly rescaled semiconvex solutions, are then convex.
The analysis of the free boundary heavily relies on the classification of blow-ups, and their convexity is of central importance in this procedure; see \cite{CRS17,CSS08}.

Here, we are able to establish the convexity of the blow-ups directly, as follows.

\begin{theorem}[Convexity of blow-ups]
\label{thm:blowups-convex}
Let $s\in(0,1)$,  $L \in \mathcal{L}_s(\lambda,\Lambda;1)$, and  $\alpha \in (0,s) \cap (0,1-s)$. 
Let $u_0 \in C^{0,1}_{loc}(\R^n) \cap C^{2s+\eps}_{loc}(\R^n)$ be such that $u_0 \ge 0$ in $\R^n$,
\begin{align}
\label{eq:blowups-convex-ass1}
\Vert \nabla u_0 \Vert_{L^{\infty}(B_R)} \le c R^{s+\alpha}\qquad \textrm{for all}\quad R\geq1,
\end{align}
\begin{align}
\label{eq:blowups-convex-ass2}
\quad L \left( \frac{u_0(x+h) - u_0(x)}{|h|} \right) \ge 0 \qquad \textrm{for all}\quad x \in \{ u_0 > 0 \},\quad  h\in\R^n.
\end{align}
Then, $u_0$ is convex in $\R^n$.

In particular, if in addition \eqref{eq:hom} holds, $0 \in \partial \{ u _0 > 0 \}$, and $u_0 \in C^1(\{ u_0 > 0 \})$ solves in the viscosity sense $L(\nabla u_0) = 0$ in $\{ u_0 > 0\}$, then $u_0(x)=\kappa (x\cdot \nu)_+^{1+s}$, with $\nu\in \mathbb S^{n-1}$, $\kappa\geq0$.
\end{theorem}

A key advantage of establishing the convexity of blow-ups without relying on semiconvexity of solutions is that we can relax the regularity assumptions of the obstacle $\phi$ in order to obtain optimal regularity estimates, as in \autoref{thm:OP_opt-reg} above. 
Moreover, it also opens the door to proving regularity estimates for $x$-dependent operators; we plan to study this in future works.

\subsection{Bernstein technique for nonlocal operators}

A central contribution of this article is the development of the Bernstein technique to obtain semiconvexity estimates for general integro-differential operators, solving an important open problem left in \cite{CDV22}. 
This technique plays an important role in the proofs of our main results concerning fully nonlinear equations and obstacle problems (see \autoref{thm:intro-FN-sc} and \autoref{thm:intro-OP-sc} above).
Still, we emphasize that, as in the local case, our nonlocal Bernstein technique works in a rather general framework, so it could be useful in completely different contexts (see e.g. \cite{ChSh95} and \cite{ASW12} for applications of the classical Bernstein technique to Hamilton-Jacobi equations and double obstacle problems).

The main insight behind the Bernstein technique is that, if derivatives of the solution $u$ are also subsolutions to an equation, then the maximum principle can be used in order to obtain regularity estimates for these solutions.
This observation can be traced back to Serge Bernstein (see \cite{Ber06,Ber10}), who noticed that
\begin{align*}
u \text{ harmonic in } B_1 ~~ \Rightarrow ~~ |\nabla u|^2 \text{ subharmonic in } B_1.
\end{align*}
By the maximum principle, this yields an estimate of $\|\nabla u\|_{L^\infty(B_1)}$ by $\|\nabla u\|_{L^\infty(\partial B_1)}$. 
This idea can be generalized by considering, instead of $|\nabla u|^2$,  the following auxiliary function 
\begin{align}
\label{eq:auxpsi}
\psi = \eta^2(\partial_e u)^2 + \sigma u^2,
\end{align}
which allows to prove interior regularity estimates for $\partial_e u$ by lower order terms. Here, $\eta \in C_c^{\infty}(\R^n)$ is a cut-off function, and $\sigma > 0$ is a suitably chosen constant. 
Such choice of $\psi$ has already been applied in \cite{CaCa95} in the context of (local) fully nonlinear PDEs. Different auxiliary functions appear for instance in the context of the mean curvature equation (see \cite{Ser69}, \cite{Wan98}). Moreover, we refer to \cite{LaUr68}, \cite{Eva79}, \cite{Bar91}, \cite{GiTr01}, and  \cite{AtCa06}, \cite{Fer16}, \cite{FeJa21}, \cite{CDV22}, \cite{ASW12} for further references on the Bernstein technique applied to elliptic equations of second order and to free boundary problems, respectively.

So far, the only work on the Bernstein technique for integro-differential operators is \cite{CDV22}, where the authors consider the auxiliary function \eqref{eq:auxpsi} and are able to treat the family of nonlocal operators that arise as affine transformations of the fractional Laplacian (i.e., kernels of the form $K(y) = |Ay|^{-n-2s}$ for some symmetric, uniformly elliptic matrix $A$). 
Moreover, they obtain Lipschitz estimates for operators belonging to the class $\mathcal{L}_s(\lambda,\Lambda;2)$ --- but not semiconvexity estimates like the ones in Theorems \ref{thm:intro-FN-sc} and \ref{thm:intro-OP-sc} above.
A central contribution of our work is to extend and improve the results of \cite{CDV22} by establishing the Bernstein technique for general nonlocal operators belonging to the natural class $\mathcal{L}_s(\lambda,\Lambda;1)$, and to solve all the problems that were left open after \cite{CDV22}. Moreover, as explained in more detail in Section \ref{sec:extensions}, we are able to extend our method to treat 
\begin{itemize}
\item parabolic equations,
\item nonsymmetric operators with drifts, and
\item nonlocal operators not necessarily comparable to $(-\Delta)^s$.
\end{itemize}

\subsubsection*{Key estimates}

The main ingredient in the proof of our semiconvexity estimates via the Bernstein technique is the following:

\begin{theorem}
\label{thm:keyest}
Let $s\in(s_0,1)$, with $s_0>0$, and $L \in \mathcal{L}_s(\lambda,\Lambda;1)$. 
Let $\eta \in C^{1,1}(\R^n)$ be such that $\eta \ge 0$.  Then, there exists $\sigma_0 = \sigma_0(n, s, \Lambda/\lambda, \Vert \eta \Vert_{C^{1,1}(\R^n)}) > 0$ such that for every $\sigma \ge \sigma_0$ and every $u \in C_{loc}^{1+2s+\eps}(\R^n) \cap L^{\infty}(\R^n)$
\begin{align}
\label{eq:keyest}
\qquad\qquad L\big(\eta^2 (\partial_e u)^2 + \sigma u^2\big) \le 2 \eta^2 L(\partial_e u) \partial_e u + 2\sigma L(u) u \qquad \textrm{in}\quad \R^n.
\end{align}
\end{theorem}

We emphasize that \autoref{thm:keyest} has so far only been established in \cite{CDV22} for affine transformations of the fractional Laplacian, by using the Caffarelli-Silvestre extension. Finding a proof of \eqref{eq:keyest} without using the extension remained open after \cite{CDV22}, even in case $L = (-\Delta)^s$ (see Open problem 1.6 in \cite{CDV22}). 
As we explain below, our proof merely relies on elementary arguments using the precise shape of $L$ (see \eqref{eq:L_intro-Bernstein}), thereby solving Open problem 1.6 in \cite{CDV22}. Moreover, our approach allows us to obtain \autoref{thm:keyest} for the natural class of operators $\mathcal{L}_s(\lambda,\Lambda;1)$ without any additional remainder terms as in \cite{CDV22}. 
This solves Open problem 1.7 in \cite{CDV22}.

It is important to emphasize that, for linear operators of second order, \eqref{eq:keyest} is a simple consequence of the product rule.
However, for nonlocal operators such an estimate is very far from trivial, and was left as an open problem in \cite{CDV22}.

\begin{remark}
The key estimate \autoref{thm:keyest} is robust with respect to the limit $s \nearrow 1$, i.e., the constants $\sigma_0$ and $C$ depend only on $n, s_0, \lambda, \Lambda$. 
Since any operator $\mathcal{L}$ of the form
\begin{align}
\label{eq:local-op}
\mathcal{L}u(x) = \sum_{i,j = 1}^n a_{ij} \partial_{ij}u(x), \qquad \text{ with} \quad \lambda |\xi|^2 \le \sum_{i,j = 1}^n a_{ij} \xi_i \xi_j \le \Lambda |\xi|^2 ~~ \forall \xi \in \R^n
\end{align}
can be approximated by a sequence of operators $L_s \in \mathcal{L}_s((1-s)\lambda,(1-s)\Lambda;1)$ as $s \nearrow 1$, our results are a true generalization of the corresponding ones for operators \eqref{eq:local-op}.
\end{remark}

Observe that the key estimate \autoref{thm:keyest} is not suitable for proving one-sided regularity estimates, such as \autoref{thm:intro-FN-sc} and \autoref{thm:intro-OP-sc}. Therefore, in order to establish semiconvexity estimates, we rely on one-sided key estimates of the following form:

\begin{theorem}
\label{thm:keyest_pos-part}
Let $s\in(s_0,1)$, with $s_0>0$, and $L \in \mathcal{L}_s(\lambda,\Lambda;1)$. 
Let $\eta \in C^{1,1}(\R^n)$ be such that $\eta \ge 0$.  Then, there exists $\sigma_0 = \sigma_0(n, s, \Lambda/\lambda, \Vert \eta \Vert_{C^{1,1}(\R^n)}) > 0$ such that for every $\sigma \ge \sigma_0$ and every $v \in C_{loc}^{1+2s+\eps}(\R^n) \cap L^{\infty}(\R^n)$
\begin{align}
\label{eq:keyest_pos-part}
\qquad\qquad L\big(\eta^2 (\partial_e v)_+^2 + \sigma v^2\big) \le 2 \eta^2 L(\partial_e v) (\partial_e v)_+ + 2\sigma L(v) v \qquad \textrm{in}\quad \R^n.
\end{align}
\end{theorem}

We will prove the one-sided estimate \eqref{eq:keyest_pos-part} by a slight modification of the proof of \autoref{thm:keyest}. 
Before our work, \eqref{eq:keyest_pos-part} was only known for affine transformations of the fractional Laplacian (see \cite{CDV22}).  Therefore, after \cite{CDV22}, it even remained unclear, whether semiconvexity estimates such as \autoref{thm:intro-FN-sc} for nonlocal fully nonlinear equations driven by operators from the class $\mathcal{L}_s(\lambda,\Lambda;1)$ do hold true (see Open problem 1.8 in \cite{CDV22}). \autoref{thm:keyest_pos-part} is the main ingredient in our proofs of the desired semiconvexity estimates \autoref{thm:intro-FN-sc}, \autoref{thm:intro-OP-sc}, the first of which solves Open problem 1.8 in \cite{CDV22}.

\subsubsection*{Difference quotients}

In  \eqref{eq:keyest}, it is apparent that one needs to consider sufficiently smooth solutions $u$ in order for $L(\partial_e u)$ to be well-defined ($u \in C_{loc}^{1+2s+\eps}(\R^n)$ is sufficient, for any $\eps > 0$). 
In case of fully nonlinear equations, this is not a problem, as one can approximate any viscosity solution $u$ by smooth solutions $u^\epsilon$, thanks to the results in \cite{Fer24}.
However, in the obstacle problem, solutions are never more regular than $C^{1+s}$, and therefore, \autoref{thm:keyest} and \autoref{thm:keyest_pos-part} cannot be used directly to derive semiconvexity estimates.
In order to be able to apply the Bernstein technique to the obstacle problem (and also to other equations where such approximation argument is not possible), we prove an estimate reminiscent of \eqref{eq:keyest} for difference quotients. 
To this end, we introduce for $h \in \R^n$
\begin{align*}
u_h(x) := \int_0^1 u(x + th) \d t \qquad \textrm{and} \qquad  D_{h}u(x) = \frac{u(x+h) - u(x)}{|h|},
\end{align*}
and establish the following.

\begin{proposition}
\label{lemma:keyest_diff-quot}
Let $s\in(s_0,1)$, with $s_0>0$, and $L \in \mathcal{L}_s(\lambda,\Lambda;1)$. 
Let $\eta \in C^{1,1}(\R^n)$ be such that $\eta \ge 0$, and  $|h| \le 1/8$. 
Then, there exists $\sigma_0 = \sigma_0(n, s, \Lambda/\lambda, \Vert \eta \Vert_{C^{1,1}(\R^n)}) > 0$ such that for every $\sigma \ge \sigma_0$ and every $u,v \in C^{2s+\eps}_{loc}(\R^n) \cap L^{\infty}(\R^n)$
\begin{align}
\label{eq:keyest_diff-quot}
L(\eta^2 (D_{h} u)^2 + \sigma u_{h}^2) &\le 2 \eta^2 L(D_{h}u) D_{h}u + 2\sigma L(u_{h}) u_{h},\\
\label{eq:keyest_diff-quot_pos-part}
L(\eta^2 (D_{h} v)_+^2 + \sigma v_{h}^2) &\le 2 \eta^2 L(D_{h}v) (D_{h}v)_+ + 2\sigma L(v_{h}) v_{h}.
\end{align}
\end{proposition}

For another version of such an estimate for difference quotients, we refer the reader to \autoref{lemma:keyest_diff-quot2}. Moreover, we mention \autoref{cor:keyest_Holder-diff-quot} and \autoref{cor:keyest_Holder-diff-quot2}, where such estimates are established for H\"older difference quotients.
Let us point out that the idea to prove Bernstein estimates for difference quotients has already been presented in an earlier version of \cite{CDV22} (see Section 3 in \cite{CDV20}). 
However, \autoref{lemma:keyest_diff-quot} was only proved for the local Laplacian in \cite{CDV20}.
We refer to \cite{DoKr07} for a version of the Bernstein technique for finite difference operators.

Finally, as we pointed out before, we are able to extend our proof of the aforementioned Bernstein key estimates in several directions. In Subsection \ref{subsec:parabolic}, we explain how to modify the Bernstein technique in order to obtain a priori regularity estimates for nonlocal parabolic equations. Moreover, in Subsection \ref{subsec:nonsymm}, we prove Bernstein key estimates for nonlocal operators with nonsymmetric kernels under the presence of additional drift terms, and in Subsection \ref{subsec:Levy} we consider nonlocal operators that are not necessarily comparable to the fractional Laplacian.

\subsection{Acknowledgments}

{The authors were supported by the European Research Council (ERC) under the Grant Agreement No 801867, and by the AEI project PID2021-125021NA-I00 (Spain). In addition, the first author was supported by the SGR project 2021 SGR 00087 (Catalunya), the AEI grant RED2022-134784-T (Spain), and the Spanish State Research Agency, through the Mar\'ia de Maeztu Program for Centers and Units of Excellence in R\&D (CEX2020-001084-M)}.

\subsection{Outline}
This article is structured as follows: In Section \ref{sec:prelim} we present and prove several auxiliary results which will be crucial in the derivation of our main results. Section \ref{sec:keyest} is devoted to the study of Bernstein key estimates for nonlocal operators and contains the proof of our main results, \autoref{thm:keyest} and \autoref{thm:keyest_pos-part}.  Moreover, in Section \ref{sec:fully-nonlinear}, we prove semiconvexity estimates for solutions to nonlocal fully nonlinear equations (see \autoref{thm:intro-FN-sc}). The Bernstein key estimates for difference quotients, and in particular \autoref{lemma:keyest_diff-quot}, are established in Section \ref{sec:diffquot}. In Section \ref{sec:obstacle}, we prove our main results on the nonlocal obstacle problem, \autoref{thm:OP_opt-reg}, and \autoref{thm:OP_bdry-reg}. By application of the Bernstein technique, we prove our main auxiliary result, the convexity of blow-ups (see \autoref{thm:blowups-convex}) and establish semiconvexity estimates for solutions (see \autoref{thm:intro-OP-sc}). Finally, in Section \ref{sec:extensions}, we discuss several extensions of our technique to parabolic equations, nonlocal equations with drifts, and nonlocal operators that are not necessarily comparable to the fractional Laplacian.

\section{Preliminaries}
\label{sec:prelim}

Let us introduce the classes of operators we will be working with throughout this article.

\begin{definition}[Regularity classes]
\label{def:classL}
Let $L$ be an integro-differential operator of the form \eqref{eq:L_intro-Bernstein} where $K :\R^n \to [0,\infty]$ is symmetric, i.e., $K(y) = K(-y)$.
\begin{itemize}
\item[(i)] We say that $L \in \mathcal{L}_s(\lambda,\Lambda)$ for some $s \in (0,1)$ and $0 < \lambda \le \Lambda$ if $L$  satisfies the following ellipticity condition:
\begin{align}
\label{eq:coercive}\tag{$K_{\asymp}$}
\lambda |y|^{-n-2s} \le K(y) \le \Lambda |y|^{-n-2s}.
\end{align}
\item[(ii)] We say that $L \in \mathcal{L}_s(\lambda,\Lambda;1)$ if $L \in \mathcal{L}_s(\lambda,\Lambda)$ and satisfies in addition
\begin{align}
\label{eq:L1}\tag{$C^1$}
|\nabla K(y)| \le \Lambda |y|^{-1} K(y).
\end{align}
\item[(iii)] We say that $L \in \mathcal{L}_s(\lambda,\Lambda;2)$ if $L \in \mathcal{L}_s(\lambda,\Lambda;1)$ and satisfies in addition
\begin{align}
\label{eq:L2}\tag{$C^2$}
|D^2 K(y)| \le \Lambda |y|^{-2} K(y).
\end{align}
\end{itemize}
\end{definition}

\begin{remark}
(i) Notice that if $L \in \mathcal{L}_s(\lambda,\Lambda)$, then $Lu(x_0)$ is well-defined for any $x_0 \in \R^n$ and $u \in C^{2s+\eps}(B_{\delta}(x_0)) \cap L^{\infty}(\R^n)$, for some $\eps, \delta > 0$.

(ii) The assumptions \eqref{eq:L1} and \eqref{eq:L2} are common in the study of higher regularity for nonlocal equations and appeared first in \cite{CaSi09,CaSi11a,CaSi11b}.
\end{remark}

Let us associate $L$ with a bilinear form $B$, defined by
\begin{align*}
B(u,v)(x) = \int_{\R^n} \big(u(x) - u(y)\big)\big(v(x) - v(y)\big) K(x-y) \d y.
\end{align*}

Sometimes, we write $L = L_{K}$ or $B = B_K$ in order to emphasize the corresponding kernel through the notation.

We observe the following nonlocal product rule:
\begin{lemma}
\label{lemma:productrule}
Let $s\in(0,1)$ and $L \in \mathcal{L}_s(\lambda,\Lambda)$.
Then, for any $u,v \in C^{2s+\eps}(B_{\delta}(x_0)) \cap L^{\infty}(\R^n)$ we have
\begin{align*}
L(uv) = uLv + vLu - B(u,v) \quad \textrm{in} \quad B_\delta(x_0).
\end{align*}
In particular,
\begin{align*}
L(u^2) = 2uLu - B(u,u) \quad \textrm{in} \quad B_\delta(x_0).
\end{align*}
\end{lemma}

\begin{proof}
We compute
\begin{align*}
u(x)v(x) - u(y)v(y) = u(x)(v(x) - v(y)) + v(x)(u(x) - u(y)) - (u(x) - u(y))(v(x) - v(y)),
\end{align*}
and the result follows.
\end{proof}

\section{Key estimates for the Bernstein technique}
\label{sec:keyest}

In this section we establish the key estimates for smooth functions. 
For this, we will need the following result, which allows us to separately consider the singularity at the origin and the behavior at infinity. 

\begin{lemma}[Kernel decomposition]
\label{lemma:kernel-decomp}
Let $s\in(0,1)$, $L \in \mathcal{L}_s(\lambda,\Lambda;1)$ and let $\eps \in (0,1)$. Then, there exist $K_1, K_2 : \R^n \to [0,\infty]$ such that $K = K_1 + K_2$ and the following properties hold true:
\begin{itemize}
\begin{multicols}{2}
\item[(i)] $\supp(K_1) \subset B_{\eps}$,
\item[(ii)] $K_1 \equiv K$ in $B_{\eps/2}$,
\item[(iii)] $K_1 \le K$, $K_2 \le K$,
\item[(iv)] $\supp(K_2) \subset \R^n \setminus B_{\eps/2}$,
\item[(v)] $|\nabla K_2| \le c_1 \eps^{-1} K$, 
\item[(vi)] $c_2\mu_{K_2}(\R^n) \le \mu_K(\R^n \setminus B_{\eps/2}) \le c_3 \mu_{K_2}(\R^n)$, 
\end{multicols}
\end{itemize}
where we denote $\mu_{K} = K(y) \d y$ and $\mu_{K_2} = K_2(y)\d y$, and the constants $c_1,c_2,c_3 > 0$ depend only on $n, s, \lambda, \Lambda$, but not on $\eps$.
\end{lemma}

\begin{proof}
Let $\psi \in C^{\infty}([0,\infty))$ be a cutoff function satisfying $0 \le \psi \le 1$, $\psi \equiv 0$ in $B_{1/2}$ and $\psi \equiv 1$ in $\R^n \setminus B_1$. Moreover, assume that $|\psi'| \le 4$. We define
\begin{align*}
K_1(y) = \left[1-\psi\left(\frac{|y|}{\eps}\right)\right]K(y), ~~ K_2(y) = \psi\left(\frac{|y|}{\eps}\right)K(y).
\end{align*}
Then, clearly, $K = K_1 + K_2$ and properties $(i)$, $(ii)$, $(iii)$, and $(iv)$ follow immediately by construction. Moreover, note that
\begin{align*}
|\nabla K_2(y)| \le \eps^{-1}|\psi'(|y|/\eps)|K(y) + \psi(|y|/\eps)|\nabla K(y)| \le 4 \eps^{-1}K(y) + 2\Lambda \eps^{-1}K(y),
\end{align*}
where we used $(iv)$ and \eqref{eq:L1}. This proves $(v)$. Note that the first estimate in $(vi)$ is a direct consequence of $(iii)$ and $(iv)$. To show the second inequality in $(vi)$, we compute using \eqref{eq:coercive}
\begin{align*}
\mu_K(\R^n \setminus B_{\eps/2}) \le c \eps^{-2s} \le c \int_{B_{2\eps} \setminus B_{\eps}} K(y) \d y = \mu_{K_2}(B_{2\eps} \setminus B_{\eps}) \le c \mu_{K_2}(\R^n),
\end{align*}
which concludes the proof.
\end{proof}

We also need the following simple observation concerning cutoff functions.

\begin{lemma}
\label{lemma:cutoff-est}
Let $s\in (0,1)$ and let $K$ be symmetric, with
\begin{align*}
K(y) \le \Lambda |y|^{-n-2s}, ~~ \supp(K) \subset B_{\eps}
\end{align*}
for some $\Lambda > 0$ and $\eps \in (0,1)$. Let $\eta \in C^{1,1}(B_1)$. Then, for any $x \in B_{1-\eps}$
\begin{align*}
L(\eta^2)(x) &\le c_1 \Vert D^2 \eta^2 \Vert_{L^{\infty}(B_{\eps}(x))} \eps^{2-2s},\\
B(\eta,\eta)(x) &\le c_2 \Vert \nabla \eta \Vert_{L^{\infty}(B_{\eps}(x))}^2 \eps^{2-2s},
\end{align*}
where $c_1, c_2 > 0$ are constants depending only on $n, s, \Lambda$.
\end{lemma}

\begin{proof}
For the first estimate we compute
\begin{align*}
L(\eta^2)(x) &= \int_{B_{\eps}(x)}(\eta^2(x) - \eta^2(y) + \nabla \eta^2(x)(x-y)) K(x-y) \d y\\
&\le \Lambda\Vert D^2 \eta^2 \Vert_{L^{\infty}(B_{\eps}(x))} \int_{B_{\eps}(x)} |x-y|^{2-n-2s} \d y\\
&\le c\Lambda\Vert D^2 \eta^2 \Vert_{L^{\infty}(B_{\eps}(x))} \eps^{2-2s}.
\end{align*}
Note that we used the even symmetry of $K$ in the first identity. For the second estimate, we observe
\begin{align*}
B(\eta,\eta)(x) &= \int_{B_{\eps}(x)} |\eta(x) - \eta(y)|^2 K(x-y) \d y\\
&\le \Lambda\Vert \nabla \eta \Vert_{L^{\infty}(B_{\eps}(x))}^2 \int_{B_{\eps}(x)} |x-y|^{2-n-2s} \d y\\
&\le c\Lambda\Vert \nabla \eta \Vert_{L^{\infty}(B_{\eps}(x))}^2 \eps^{2-2s}.
\end{align*}
\end{proof}

With this at hand, we can now start the proof of our key estimates in Theorems \ref{thm:keyest} and \ref{thm:keyest_pos-part}.

\subsection{First order estimates}

The goal of this section is to establish the key estimate for the Bernstein technique (see \autoref{thm:keyest}), which will be used in order to prove first derivative estimates.

Before we prove \autoref{thm:keyest} let us list two equivalent formulations of  \eqref{eq:keyest} that will turn out to be more convenient to prove:

\begin{remark}
The following two estimates are equivalent to \eqref{eq:keyest}:
\begin{align}
\label{eq:equivkey0}
L(\eta^2) (\partial_e u)^2 - B(\eta^2, (\partial_e u)^2) &\le \eta^2 B(\partial_e u , \partial_e u) + \sigma B(u,u),\\
\label{eq:equivkey1}
\int_{\R^n} (\eta^2(x) - \eta^2(y)) (\partial_e u(y))^2 K(x-y) \d y &\le \eta^2(x) B(\partial_e u , \partial_e u)(x) + \sigma B(u,u)(x).
\end{align}
This can be seen as follows: With the help of the nonlocal product rule \autoref{lemma:productrule} we compute:
\begin{align*}
L(\eta^2 (\partial_e u)^2 + \sigma u^2) &= L(\eta^2) (\partial_e u)^2 + 2\eta^2 L(\partial_e u)\partial_e u - \eta^2 B(\partial_e u , \partial_e u) - B(\eta^2, (\partial_e u)^2)\\
&\quad + 2\sigma L(u) u - \sigma B(u,u).
\end{align*}
Therefore, \eqref{eq:keyest} is equivalent to \eqref{eq:equivkey0}. Moreover, the left hand side in \eqref{eq:equivkey0} can be rewritten as follows:
\begin{align*}
L(\eta^2) (\partial_e u)^2 - B(\eta^2, (\partial_e u)^2) = \int_{\R^n} (\eta^2(x) - \eta^2(y)) (\partial_e u(y))^2 K(x-y) \d y.
\end{align*}
This is a simple consequence of the following identity:
\begin{align*}
(\eta^2(x) - \eta^2(y)) (\partial_e u(x))^2 - (\eta^2(x) - \eta^2(y))((\partial_e u(x))^2 - (\partial_e u(y))^2) = (\eta^2(x) - \eta^2(y)) (\partial_e u(y))^2.
\end{align*}
Thus, \eqref{eq:keyest} and \eqref{eq:equivkey0} are both equivalent to \eqref{eq:equivkey1}.
\end{remark}

Moreover, the following interpolation inequality will turn out to be crucial in the proof of \autoref{thm:keyest}.

\begin{lemma}
\label{lemma:interpol}
Let $s\in(0,1)$ and $\delta > 0$. Assume that $K : \R^n \to [0,\infty]$ satisfies for some $0 < \lambda \le \Lambda$:
\begin{align}
\label{eq:interpol-coercive}
\lambda |y|^{-n-2s} &\le K(y) \le \Lambda |y|^{-n-2s}~~ \forall y \in B_{\delta},\\
\label{eq:interpol-C1}
|\nabla K(y)| &\le \Lambda |y|^{-1} K(y)~~ \forall y \in B_{\delta}.
\end{align}
Then, for every $x \in \R^n$ and $u \in C^{0,1}(B_{\delta}(x))$ it holds
\begin{align*}
\big(\partial_e u(x)\big)^2 \le \delta^{2s} B(\partial_e u, \partial_e u)(x) + c\delta^{2s-2} B(u,u)(x),
\end{align*}
where $c = c(n,s,\lambda,\Lambda) > 0$ does not depend on $\delta$.
\end{lemma}

\begin{proof}
First, given any $\delta > 0$, we construct an auxiliary kernel $K_{\delta} : \R^n \to [0,\infty]$ satisfying the following properties:
\begin{itemize}
\begin{multicols}{2}
\item[(1)] $K_{\delta}^2(y) \le c_1 K(y) |y|^2$ for $y \in B_{\delta}$,
\item[(2)] $|\nabla K_{\delta}(y)|^2 \le c_2 K(y)$ for $y \in B_{\delta}$,
\item[(3)] $\supp(K_{\delta}) \subset B_{\delta}(0)$,
\item[(4)] $c_3 \delta^{\frac{n}{2} - s + 1} \le \mu_{K_{\delta}}(B_{\delta}) \le c_4 \delta^{\frac{n}{2} - s + 1}$ ,
\end{multicols}
\end{itemize}
where we define $\mu_{K_{\delta}} = K_{\delta}(y) \d y$, and $c_1, c_2, c_3, c_4 > 0$ are constants depending only on $n, s, \lambda, \Lambda$.\\
To do so, we proceed in the same way as in the proof of \autoref{lemma:kernel-decomp}. Indeed, let $\psi \in C^{\infty}([0,\infty))$ by a cutoff function satisfying $\psi \equiv 1$ in $B_{1/2}$, $\psi \equiv 0$ in $\R^n \setminus B_1$, $0 \le \psi \le 1$ and $|\psi'| \le 4$. Then, we define
\begin{align*}
K_{\delta}(y) = \psi \left( \frac{|y|}{\delta} \right)K(y)|y|^{\frac{n}{2}+s+1}
\end{align*} 
and observe that the properties (1), (3) and (4) follow immediately from the construction and \eqref{eq:interpol-coercive}. To prove (2), we compute for $y \in B_{\delta}$, using \eqref{eq:interpol-coercive} and \eqref{eq:interpol-C1}
\begin{align*}
|\nabla K_{\delta}(y)|^2 &\lesssim \left(\frac{|y|}{\delta}\right)^2 \left|\psi'\left(\frac{|y|}{\delta}\right)\right|^2 K^2(y) |y|^{n+2s} + \psi^2\left(\frac{|y|}{\delta}\right)\left(|\nabla K(y)|^2 |y|^{n+2s+2} + K^2(y) |y|^{n+2s} \right)\\
&\lesssim K(y).
\end{align*}

Having constructed $K_{\delta}$, let us turn to proving the desired interpolation estimate.\\
We compute, using (4), and the notation $\mu_{K_{\delta}}(x,\d y) = K_{\delta}(x-y) \d y$:
\begin{align*}
\partial_e u(x) &= \dashint_{B_{\delta}(x)} (\partial_e u(x) - \partial_e u(y)) \mu_{K_{\delta}}(x,\d y) + \dashint_{B_{\delta}(x)} \partial_e u(y) \mu_{K_{\delta}}(x,\d y) \\
&=\frac{|B_\delta(x)|}{\mu_{K_\delta}(B_\delta(x))} \left[\dashint_{B_{\delta}(x)} (\partial_e u(x) - \partial_e u(y)) K_{\delta}(x-y) \d y + \dashint_{B_{\delta}(x)} \partial_e u(y) K_{\delta}(x-y) \d y\right] \\
&\lesssim \delta^{\frac{n}{2} + s - 1} \dashint_{B_{\delta}(x)} (\partial_e u(x) - \partial_e u(y)) K_{\delta}(x-y) \d y + \delta^{\frac{n}{2} + s - 1} \dashint_{B_{\delta}(x)} \partial_e u(y) K_{\delta}(x-y) \d y.
\end{align*}
By Jensen's inequality, integration by parts, and (1), (2), and (3):
\begin{align*}
&|\partial_e u(x)|^2 \lesssim \delta^{n + 2s - 2} \left[ \dashint_{B_{\delta}(x)} \hspace{-0.1cm}(\partial_e u(x) - \partial_e u(y))^2 K_{\delta}^2(x-y) \d y + \left(\dashint_{B_{\delta}(x)} \hspace{-0.1cm} \partial_e u(y) K_{\delta}(x-y) \d y \right)^2 \right]\\
&\lesssim \delta^{n + 2s - 2} \left[ \dashint_{B_{\delta}(x)} \hspace{-0.1cm} (\partial_e u(x) - \partial_e u(y))^2 K_{\delta}^2(x-y) \d y + \left(\dashint_{B_{\delta}(x)} \hspace{-0.1cm} (u(y) - u(x)) \partial_e K_{\delta}(x-y) \d y \right)^2 \right]\\
&\lesssim \delta^{2s - 2}  \left[\int_{B_{\delta}(x)} \hspace{-0.1cm} (\partial_e u(x) - \partial_e u(y))^2 K(x-y)|x-y|^{2} \d y + \int_{B_{\delta}(x)} \hspace{-0.1cm} (u(y) - u(x))^2 \left(\partial_e K_{\delta}(x-y)\right)^2 \d y \right]\\
&\lesssim \delta^{2s} \int_{B_{\delta}(x)} \hspace{-0.1cm} (\partial_e u(x) - \partial_e u(y))^2 K(x-y) \d y + \delta^{2s - 2} \int_{B_{\delta}(x)} \hspace{-0.1cm} (u(y) - u(x))^2 K(x-y) \d y.
\end{align*}
Therefore, we obtain
\begin{align*}
\big(\partial_e u(x)\big)^2 \le c_1\delta^{2s} B(\partial_e u, \partial_e u)(x) + c_2\delta^{2s-2} B(u,u)(x).
\end{align*}
To conclude, we may repeat the computation with $c_1^{-1/2s}\delta$ instead of $\delta$ if $c_1 > 1$.
\end{proof}

\begin{remark}
\label{remark:robust-1}
Note that if we replace the constants $\lambda, \Lambda$ in \eqref{eq:interpol-coercive} and \eqref{eq:interpol-C1} by $\lambda (1-s), \Lambda (1-s)$, where $s \in (s_0, 1)$, then the constant $c > 0$ in \autoref{lemma:interpol} would only depend on $n, s_0, \lambda, \Lambda$. 
To see this, one redefines $K_{\delta}(y) = \sqrt{1-s}\,\psi(|y|/\delta)K(y)|y|^{\frac{n}{2}+s+1}$. 
\end{remark}

Now, we are in a position to prove the key estimate \autoref{thm:keyest}.

\begin{proof}[Proof of \autoref{thm:keyest}]
Throughout the proof, we will denote by $c > 0$ any constant that only depends on $n,s,\lambda,\Lambda$ and whose value might change from line to line. Note that we can assume without loss of generality that $\eta(x) > 0$, since otherwise the desired estimate \eqref{eq:equivkey1} is trivially satisfied.

For every $x$, let us decompose the kernel $K$ into two parts $K = K_1 + K_2$, taking care separately of the behavior at zero and at infinity. In order to do so, we choose a parameter 
\[\eps := \gamma \eta(x) > 0,\]
where $\gamma = \gamma(\Vert \eta \Vert_{C^{1,1}(\R^n)}) \in (0,1)$ will be determined later (see \eqref{eq:gamma_def_A}), and choose $K_1$ and $K_2$ as in \autoref{lemma:kernel-decomp} with respect to $\eps$. Consequently $K_1$ and $K_2$ satisfy the properties $(i)-(vi)$.

\textbf{Step 1:} Let us first explain how to treat the integrals involving $K_2$. To be precise, we will show that for some constant $\sigma_2 = \sigma_2(n,s,\lambda,\Lambda) > 0$:
\begin{align}
\label{eq:genK_2-estimate_A}
\begin{split}
\int_{\R^n} &(\eta^2(x) - \eta^2(y)) (\partial_e u(y))^2 K_2(x-y) \d y  \\
&\le \eta^2(x) \int_{\R^n} (\partial_e u(x) - \partial_e u(y))^2 K_2(x-y) \d y + \sigma_2 \frac{\eta^2(x)}{\eps^2} B_{K}(u,u)(x).
\end{split}
\end{align}
Note that the claim \eqref{eq:genK_2-estimate_A} follows if we manage to prove: 
\begin{align*}
\eta^2(x) &\int_{\R^n} (\partial_e u(y))^2 K_2(x-y) \d y\\
&\le \eta^2(x)\int_{\R^n} (\partial_e u(x) - \partial_e u(y))^2 K_2(x-y) \d y + \sigma_2 \frac{\eta^2(x)}{\eps^2} B_K(u,u)(x),
\end{align*}
which is equivalent to 
\begin{align}
\label{eq:A3help_A}
\begin{split}
2 \partial_e u(x) \eta^2(x) &\int_{\R^n} \partial_e u(y) K_2(x-y) \d y\\
&\le (\partial_e u(x))^2 \eta^2(x) \int_{\R^n} K_2(x-y) \d y + \sigma_2\frac{\eta^2(x)}{\eps^2} B_K(u,u)(x).
\end{split}
\end{align}
Note that the first term on the right hand side in \eqref{eq:A3help_A} is finite, since $K_2$ is integrable due to $(iii)$ and $(iv)$ in \autoref{lemma:kernel-decomp}.
To prove \eqref{eq:A3help_A}, we introduce the measure $\mu_{K_2}(x,\d y) = K_2(x-y)\d y$ and use the Young's inequality as follows:
\begin{align*}
2 \partial_e u(x) \eta^2(x) \int_{\R^n} \partial_e u(y) K_2(x-y) \d y &\le \mu_{K_2}(x,\R^n) (\partial_e u(x))^2 \eta^2(x)\\
&\quad + \eta^2(x) \mu_{K_2}(x,\R^n)^{-1} \left( \int_{\R^n}  \partial_e u(y) \mu_{K_2}(x, \d y) \right)^2\\
&=: J_1 + J_2.
\end{align*}
For $J_1$ we obtain
\begin{align*}
J_1 = \mu_{K_2}(x,\R^n) (\partial_e u(x))^2 \eta^2(x)  = (\partial_e u(x))^2 \eta^2(x) \int_{\R^n} K_2(x-y) \d y,
\end{align*}
which coincides with the first term on the right hand side of \eqref{eq:A3help_A}. In order to estimate $J_2$, let us recall that by $(iv)$, $(v)$, and $(vi)$ in \autoref{lemma:kernel-decomp}, we have
\begin{align*}
\supp(K_2(x-\cdot)) \subset \R^n \setminus B_{\eps/2}(x), ~~ |\nabla K_2(x-\cdot)| \le c\eps^{-1} K(x-\cdot), ~~ \frac{\mu_{K}(x,\R^n \setminus B_{\eps/2}(x))}{\mu_{K_2}(x,\R^n)} \le c.
\end{align*}
Thus, using integration by parts and Jensen's inequality:
\begin{align*}
J_2 &= \eta^2(x) \mu_{K_2}(x,\R^n)^{-1} \left(\int_{\R^n} (u(y) - u(x)) \partial_e K_2(x-y) \d y\right)^2\\
&\le c \frac{\eta^2(x)}{\eps^2} \mu_{K_2}(x,\R^n)^{-1} \left(\int_{\R^n \setminus B_{\eps/2}(x)} |u(y) - u(x)| K(x-y) \d y\right)^2\\
&=  c \frac{\eta^2(x)}{\eps^2} \frac{\mu_{K}(x,\R^n \setminus B_{\eps/2}(x))^2}{\mu_{K_2}(x,\R^n)} \left(\dashint_{\R^n \setminus B_{\eps/2}(x)} |u(y) - u(x)| \mu_{K}(x,\d y) \right)^2\\
&\le c \frac{\eta^2(x)}{\eps^2} \int_{\R^n} (u(x) - u(y))^2 K(x-y) \d y.
\end{align*}
Altogether, by combining the estimates for $J_1$ and $J_2$, we have shown \eqref{eq:A3help_A}, and therefore \eqref{eq:genK_2-estimate_A}, as desired.\\

\textbf{Step 2:} We claim that
\begin{align}
\label{eq:K_1-estimate_A}
L_{K_1}(\eta^2)(x) (\partial_e u(x))^2  - B_{K_1}(\eta^2, (\partial_e u)^2)(x) \le \eta^2(x) B_{K_1}(\partial_e u , \partial_e u)(x) + \sigma_1 B_K(u,u)(x),
\end{align}
where $\sigma_1 = \sigma_1(n,s,\lambda,\Lambda,\Vert \eta \Vert_{C^{1,1}(\R^n)})> 0$ is a constant. Note that by combining \eqref{eq:K_1-estimate_A} with \eqref{eq:genK_2-estimate_A} and using that $\eps = \gamma \eta(x)$, we obtain that for every $x \in \R^n$:
\begin{align}
\label{eq:result}
L_K(\eta^2)(x) (\partial_e u(x))^2 - B_K(\eta^2, (\partial_e u)^2)(x) \le \eta^2(x) B_K(\partial_e u , \partial_e u)(x) + \sigma B_K(u,u)(x),
\end{align}
where $\sigma = \sigma_1 + \sigma_2 \gamma^{-2} > 0$. This proves the desired result.

Let us now prove \eqref{eq:K_1-estimate_A}. Recall that $\supp(K_1(x - \cdot)) \subset B_{\eps}(x)$ by $(i)$ in \autoref{lemma:kernel-decomp}. By Young's again, we compute for $A = \frac{1}{8}(\Vert \nabla \eta \Vert_{L^{\infty}(\R^n)} + 2)^{-2} > 0$:
\begin{align*}
- B_{K_1}(\eta^2, (\partial_e u)^2)(x) &= -\int_{B_{\eps}(x)} (\eta^2(x) - \eta^2(y))[(\partial_e u(x))^2-(\partial_e u(y))^2] K_1(x-y) \d y \\
&\le A\int_{B_{\eps}(x)} (\eta(x) + \eta(y))^2 (\partial_e u(x) - \partial_e u(y))^2 K_1(x-y) \d y\\
&\quad + \frac{1}{4A} \int_{B_{\eps}(x)} (\eta(x) - \eta(y))^2 (\partial_e u(x) + \partial_e u(y))^2 K_1(x-y) \d y\\
&=: I_1 + I_2.
\end{align*}
For $I_1$, we use that for $y \in B_{\eps}(x)$
\begin{align*}
|\eta(x) - \eta(y)| \le \Vert \nabla \eta \Vert_{L^{\infty}(B_{\eps}(x))} |x-y| \le \Vert \nabla \eta \Vert_{L^{\infty}(\R^n)} \eps = \gamma \Vert \nabla \eta \Vert_{L^{\infty}(\R^n)}  \eta(x),
\end{align*}
and therefore
\begin{align*}
    (\eta(x)+\eta(y))^2 \leq (2+\gamma \Vert \nabla \eta \Vert_{L^{\infty}(B_{\eps}(x))})^2\eta(x)^2 \leq (2+\Vert \nabla \eta \Vert_{L^{\infty}(B_{\eps}(x))})^2\eta(x)^2.
\end{align*}
We obtain
\begin{align*}
I_1 &\le A (\Vert \nabla \eta \Vert_{L^{\infty}(\R^n)} + 2)^2 \eta^2(x) \int_{B_{\eps}(x)} (\partial_e u(x) - \partial_e u(y))^2 K_1(x-y) \d y \le \frac{\eta^2(x)}{8} B_{K_1}(\partial_e u, \partial_e u)(x).
\end{align*}
For $I_2$, we make use of the following algebraic inequality
\begin{align*}
(a+b)^2 \le (a+b)^2 + (3a - b)^2 = 8 a^2 + 2(a - b)^2
\end{align*}
and apply it to $a = \partial_e u(x)$ and $b = \partial_e u(y)$. This yields, for some $c_1 > 0$,
\begin{align*}
I_2 &\le 2 A^{-1} (\partial_e u(x))^2 \int_{B_{\eps}(x)} (\eta(x) - \eta(y))^2 K_1(x-y) \d y\\
&\quad + \frac{1}{2A} \int_{B_{\eps}(x)}  (\eta(x) - \eta(y))^2 (\partial_e u(x) - \partial_e u(y))^2 K_1(x-y) \d y\\
&\le \left[2c_1 A^{-1} \gamma^{2-2s} \Vert \nabla \eta \Vert_{L^{\infty}(\R^n)}^2\right] \eta^{2-2s}(x) (\partial_e u(x))^2\\
&\quad + \left[(2A)^{-1} \Vert \nabla \eta \Vert_{L^{\infty}(\R^n)}^2  \gamma^2\right] \eta^2(x) B_{K_1}(\partial_e u , \partial_e u)(x)\\
&= I_{2,1} + I_{2,2},
\end{align*}
where we applied \autoref{lemma:cutoff-est}.
By choosing 
\begin{align}
\label{eq:gamma_def_A}
\gamma = (4A^{-1} \Vert \nabla \eta \Vert_{L^{\infty}(\R^n)}^2)^{-\frac{1}{2}} \wedge 1,
\end{align}
we estimate
\begin{align*}
I_{2,2} \le \frac{\eta^2(x)}{8} B_{K_1}(\partial_e u , \partial_e u)(x).
\end{align*}
We apply \autoref{lemma:interpol} to $K_1$ with $\delta = [8c_1 A^{-1}\gamma^{2-2s} \Vert \nabla \eta \Vert_{L^{\infty}(\R^n)}^2]^{-\frac{1}{2s}} \eta(x) \wedge \frac{\gamma}{2} \eta(x)$ in order to estimate $I_{2,1}$. Note that $K_1$ satisfies the assumptions of \autoref{lemma:interpol} due to \autoref{lemma:kernel-decomp}. This yields
\begin{align*}
I_{2,1} \le \frac{\eta^2(x)}{4} B_{K_1}(\partial_e u , \partial_e u)(x) + \sigma_{1,1} B_{K_1}(u,u)(x),
\end{align*}
where, for $c_2 > 0$ being the constant in \autoref{lemma:interpol},
\begin{align*}
\sigma_{1,1} = \frac{c_2}{4}\left[8c_1A^{-1}\gamma^{2-2s} \Vert \nabla \eta \Vert_{L^{\infty}(\R^n)}^2\right]^{\frac{1}{s}} \vee 2^{3-2s}c_1c_2A^{-1}\Vert \nabla \eta \Vert_{L^{\infty}(\R^n)}^2 > 0.
\end{align*}
By combination of the estimates for $I_1$, $I_{2,1}$ and $I_{2,2}$, we obtain
\begin{align}
\label{eq:-BK_upper_bound}
- B_{K_1}(\eta^2, (\partial_e u)^2)(x) &\le \frac{1}{2} \eta^2(x)B_{K_1}(\partial_e u, \partial_e u)(x) + \sigma_{1,1} B(u,u)(x).
\end{align}
Finally, we observe that for some $c_3 > 0$, by \autoref{lemma:cutoff-est}
\begin{align*}
L_{K_1}(\eta^2)(x) \le c_3 \Vert D^2 \eta^2 \Vert_{L^{\infty}(\R^n)} \eps^{2-2s} = \left[c_3 \Vert D^2 \eta^2 \Vert_{L^{\infty}(\R^n)} \gamma^{2-2s}\right] \eta^{2-2s}(x),
\end{align*}
and apply \autoref{lemma:interpol} to $K_1$ with $\delta = [4c_3 \Vert D^2 \eta^2 \Vert_{L^{\infty}(\R^n)} \gamma^{2-2s}]^{-\frac{1}{2s}} \eta(x) \wedge \frac{\gamma}{2}\eta(x)$ and obtain
\begin{align*}
L_{K_1}(\eta^2)(x)(\partial_e u(x))^2 \le \frac{1}{4} \eta^2(x)B_{K_1}(\partial_e u, \partial_e u)(x) + \sigma_{1,2} B(u,u)(x),
\end{align*}
where, for $c_2 > 0$ being the constant in \autoref{lemma:interpol}
\begin{align*}
\sigma_{1,2} = \frac{c_2}{4}\left[4c_3 \Vert D^2 \eta^2 \Vert_{L^{\infty}(\R^n)} \gamma^{2-2s}\right]^{\frac{1}{s}} \vee 2^{2-2s}c_2c_3\Vert D^2 \eta^2 \Vert_{L^{\infty}(\R^n)} > 0.
\end{align*}
Consequently, we obtain \eqref{eq:K_1-estimate_A} with $\sigma_1 = \sigma_{1,1} + \sigma_{1,2}$, as desired.
\end{proof}

\begin{remark}
\label{remark:robust-2}
Note that if we choose $L \in \mathcal{L}_s((1-s)\lambda,(1-s)\Lambda;1)$ and $s \ge s_0$ for some $s_0 \in (0,1)$, then the constant $\sigma = \sigma_1+\sigma_2\gamma^{-2} > 0$ in \eqref{eq:keyest} will depend only on $n,s_0,\lambda,\Lambda$.
This is because the constants in Step 1 will not depend on $s$, and the constant in Step 2 will be of the form 
\begin{align*}
\sigma_1 = \sigma_{1,1}+\sigma_{1,2},
\end{align*}
with
\begin{align*}
\sigma_{1,1} = \frac{c_2}{4}\left[8c_1A^{-1}\gamma^{2-2s} \Vert \nabla \eta \Vert_{L^{\infty}(\R^n)}^2\right]^{\frac{1}{s}} \vee 2^{3-2s}c_1c_2A^{-1}\Vert \nabla \eta \Vert_{L^{\infty}(\R^n)}^2, \\
\sigma_{1,2} = \frac{c_2}{4}\left[4c_3 \Vert D^2 \eta^2 \Vert_{L^{\infty}(\R^n)} \gamma^{2-2s}\right]^{\frac{1}{s}} \vee 2^{2-2s}c_2c_3\Vert D^2 \eta^2 \Vert_{L^{\infty}(\R^n)},
\end{align*}
and these expressions remain bounded as $s \nearrow 1$. The constants $\gamma, \sigma_2, c_1, c_2, c_3, A > 0$ depend only on $n, s_0, \lambda, \Lambda$.
\end{remark}

\subsection{One-sided second order estimates}

In this section, we prove another key estimate, reminiscent of \autoref{thm:keyest}, which will allow us to prove one-sided second derivative bounds for solutions to certain PDEs driven by $L$.

\begin{remark}
The following two estimates are equivalent to \eqref{eq:keyest_pos-part}:
\begin{align}
\label{eq:equivkey_pos-part0}
\begin{split}
L(\eta^2) (\partial_e v)_+^2 &- B(\eta^2, (\partial_e v)_+^2) \\
&\le \eta^2 B((\partial_e v)_+ , (\partial_e v)_+) + 2\eta^2 (\partial_e v)_+ \left[ L(\partial_e v) - L((\partial_e v)_+)  \right] + \sigma B(v,v)\\
&= \eta^2 B((\partial_e v)_+ , (\partial_e v)_+) - 2\eta^2 (\partial_e v)_+  L((\partial_e v)_-) + \sigma B(v,v),
\end{split}
\end{align}
\begin{align}
\label{eq:equivkey_pos-part1}
\begin{split}
\int_{\R^n}& (\eta^2(x) - \eta^2(y)) ((\partial_e v(y))_+)^2 K(x-y) \d y\\
&\le \eta^2 B((\partial_e v)_+ , (\partial_e v)_+)(x) - 2\eta^2 (\partial_e v)_+  L((\partial_e v)_-)(x) + \sigma B(v,v)(x).
\end{split}
\end{align}
The proof goes by the same arguments as for \eqref{eq:equivkey0} and \eqref{eq:equivkey1}.
\end{remark}

The strategy of our proof is similar to the one of the first key estimate \autoref{thm:keyest}. We will again split the kernel into two parts $K_1$ and $K_2$, taking care of the singularity at zero, and the decay at infinity, respectively. In order to treat $K_1$ we need to prove an interpolation inequality similar to \autoref{lemma:interpol}:

\begin{lemma}
\label{lemma:interpol_pos-part}
Let $\delta > 0$. Assume that $K$ satisfies \eqref{eq:interpol-coercive} and \eqref{eq:interpol-C1}.
Then, for every $x \in \R^n$ and $v \in C^{0,1}(B_{\delta}(x))$ it holds
\begin{align*}
(\partial_e v(x))_+^2 \le \delta^{2s} B((\partial_e v)_+, (\partial_e v)_+)(x) - \delta^{2s} L((\partial_e v)_-)(x)(\partial_e v)_+(x) + c\delta^{2s-2} B(v,v)(x).
\end{align*}
where $c = c(n,s,\lambda) > 0$ does not depend on $\delta$.
\end{lemma}

\begin{proof}
Assume that $\partial_e v(x) > 0$, otherwise the estimate is trivial.
Let $K_{\delta}$ be as in the proof of \autoref{lemma:interpol}. By following the proof of \autoref{lemma:interpol}, we obtain
\begin{align*}
\partial_e v(x) &\le c\delta^{\frac{n}{2} + s - 1} \left[ \dashint_{B_{\delta}(x)} \hspace{-0.1cm}(\partial_e v(x) - \partial_e v(y)) K_{\delta}(x-y) \d y + \dashint_{B_{\delta}(x)} \hspace{-0.1cm} \partial_e v(y) K_{\delta}(x-y) \d y \right]\\
&=c \delta^{\frac{n}{2} + s - 1} \left[ \dashint_{B_{\delta}(x)} \hspace{-0.1cm} ((\partial_e v(x))_+ - (\partial_e v(y))_+) K_{\delta}(x-y) \d y + \dashint_{B_{\delta}(x)} \hspace{-0.1cm} (\partial_e v(y))_- K_{\delta}(x-y) \d y \right]\\
&\quad + c\delta^{\frac{n}{2} + s - 1} \dashint_{B_{\delta}(x)} \hspace{-0.1cm} \partial_e v(y) K_{\delta}(x-y) \d y.
\end{align*}
By multiplication of the aforementioned estimate with $\partial_e v(x)$ on both sides, and applying Young's inequality, we obtain
\begin{align*}
(\partial_e v(x))^2 &\le c\left(\delta^{\frac{n}{2} + s - 1} \dashint_{B_{\delta}(x)} ((\partial_e v(x))_+ - (\partial_e v(y))_+) K_{\delta}(x-y) \d y \right)^2\\
&\quad + c\delta^{\frac{n}{2} + s - 1} \partial_e v(x) \dashint_{B_{\delta}(x)} (\partial_e v(y))_- K_{\delta}(x-y) \d y\\
&\quad + c\left( \delta^{\frac{n}{2} + s - 1} \dashint_{B_{\delta}(x)} \partial_e v(y) K_{\delta}(x-y) \d y\right)^2 + \frac{1}{2}(\partial_e v(x))^2.
\end{align*}
From here, with the first and third term, we proceed as in the proof of \autoref{lemma:interpol}. The fourth term can be absorbed to the left hand side. This yields
\begin{align*}
(\partial_e v(x))^2 &\lesssim \delta^{2s} \int_{B_{\delta}(x)} ((\partial_e v(x))_+ - (\partial_e v(y))_+)^2 K(x-y) \d y\\
&\quad + \delta^{\frac{n}{2} + s - 1} \partial_e v(x) \dashint_{B_{\delta}(x)} (\partial_e v(y))_- K_{\delta}(x-y) \d y \\
&\quad + \delta^{2s - 2} \int_{B_{\delta}(x)} (v(y) - v(x))^2 K(x-y) \d y.
\end{align*}
In order to treat the second term on the right hand side, we recall the definition of $K_{\delta}$ and deduce the following estimate:
\begin{align}
\label{eq:help-nonsymm-one-sided}
\begin{split}
\delta^{\frac{n}{2} + s - 1} \partial_e v(x) \dashint_{B_{\delta}(x)} & (\partial_e v(y))_- K_{\delta}(x-y) \d y \\
&\le \delta^{-\frac{n}{2} + s - 1} \partial_e v(x) \int_{B_{\delta}(x)} (\partial_e v(y))_- K(x-y)|x-y|^{\frac{n}{2} + s + 1} \d y \\
&\le \delta^{2s} \partial_e v(x)\int_{B_{\delta}(x)} (\partial_e v(y))_- K(x-y) \d y \\
&= - \delta^{2s} \partial_e v(x) \int_{B_{\delta}(x)} ((\partial_e v(x))_- - (\partial_e v(y))_-) K(x-y) \d y\\
&\le - \delta^{2s} L((\partial_e v)_-)(x)\partial_e v(x).
\end{split}
\end{align}
Altogether, taking $\delta$ smaller if necessary as in \autoref{lemma:interpol} we obtain the desired result.
\end{proof}

We are now ready to give the proof of \autoref{thm:keyest_pos-part}.

\begin{proof}[Proof of \autoref{thm:keyest_pos-part}]
Let $\eps = \gamma \eta(x)$ as in the proof of \autoref{thm:keyest}. Moreover, define $K_1$ and $K_2$ as in \autoref{lemma:kernel-decomp} with respect to $\eps$.\\
\textbf{Step 1:} Let us first explain how to treat the integrals involving $K_2$. We show that for some uniform constant $\sigma_2 > 0$:
\begin{align}
\label{eq:genK_2-estimate_pos-part}
\begin{split}
\int_{\R^n} (\eta^2(x) - \eta^2(y)) & (\partial_e v(y))_+^2  K_2(x-y) \d y \le \eta^2(x) B_{K_2}((\partial_e v)_+,(\partial_e v)_+)(x)\\
& - 2\eta^2(x) (\partial_e v)_+(x) L_{K_2}((\partial_e  v)_-)(x) + \sigma_2 \frac{\eta^2(x)}{\eps^2} B_{K}(v,v)(x).
\end{split}
\end{align}
The proof follows by the same idea as in \autoref{thm:keyest}. Note that because our key estimates contain positive parts, there appears the additional term $- 2\eta^2(x) (\partial_e v)_+(x) L_{K_2}((\partial_e  v)_-)(x)$ in \eqref{eq:genK_2-estimate_pos-part}. Since this term is nonnegative, we can compensate the possible smallness of the other terms on the right hand side due to the consideration of positive parts.\\

The estimate \eqref{eq:genK_2-estimate_pos-part} follows if we manage to prove: 
\begin{align*}
\eta^2(x) \int_{\R^n} (\partial_e v(y))_+^2  K_2(x-y) \d y \le & \ \eta^2(x)\int_{\R^n} ((\partial_e v(x))_+ - (\partial_e v(y))_+)^2 K_2(x-y) \d y\\
& - 2\eta^2(x) (\partial_e v(x))_+ L_{K_2}((\partial_e  v)_-)(x) + c \frac{\eta^2(x)}{\eps^2} B_K(v,v)(x),
\end{align*}
which is equivalent to 
\begin{align}
\label{eq:A3help_pos-part}
\begin{split}
2 (\partial_e v(x))_+ \eta^2(x) &\int_{\R^n} (\partial_e v(y))_+ K_2(x-y) \d y \le (\partial_e v(x))_+^2 \eta^2(x) \int_{\R^n} K_2(x-y) \d y\\
& + 2(\partial_e v(x))_+ \eta^2(x) \int_{\R^n} (\partial_e v(y))_- K_2(x-y) \d y + c\frac{\eta^2(x)}{\eps^2} B_K(v,v)(x).
\end{split}
\end{align}
Here, we used that $(\partial_e v(x))_+(\partial_e v(x))_- = 0$, and hence
\begin{align*}
	- 2 \eta^2(x) (\partial_e v(x))_+ L_{K_2}((\partial_e  v)_-)(x)
	 &= 2\eta^2(x) (\partial_e v(x))_+\int_{\R^n} ((\partial_e v(y))_--(\partial_e v(x))_-) K_2(x-y) \d y\\
	&= 2\eta^2(x) (\partial_e v(x))_+\int_{\R^n} (\partial_e v(y))_- K_2(x-y) \d y.
\end{align*}

To prove \eqref{eq:A3help_pos-part}, we introduce the measure $\mu_{K_2}(x,\d y) = K_2(x-y)\d y$ and estimate
\begin{align*}
2 &(\partial_e v(x))_+ \eta^2(x) \int_{\R^n} (\partial_e v(y))_+ K_2(x-y) \d y \\
&= 2 (\partial_e v(x))_+ \eta^2(x) \int_{\R^n} (\partial_e v(y))_- K_2(x-y) \d y + 2 (\partial_e v(x))_+ \eta^2(x) \int_{\R^n} \partial_e v(y) K_2(x-y) \d y \\
&\le 2 (\partial_e v(x))_+ \eta^2(x) \int_{\R^n} (\partial_e v(y))_- K_2(x-y) \d y  + \mu_{K_2}(x,\R^n) ((\partial_e v(x))_+)^2 \eta^2(x)\\
&\quad + \eta^2(x) \mu_{K_2}(x,\R^n)^{-1} \left( \int_{\R^n}  \partial_e v(y) \mu_{K_2}(x, \d y) \right)^2\\
&=: J_0 + J_1 + J_2.
\end{align*}
Note that $J_0$ already coincides with second term on the right hand side of \eqref{eq:A3help_pos-part}. In order to estimate $J_1$ and $J_2$ we proceed precisely by the same arguments as in the estimation of $J_1$ and $J_2$ in Step 1 of the proof of \autoref{thm:keyest}. This proves \eqref{eq:A3help_pos-part}, and therefore \eqref{eq:genK_2-estimate_pos-part}, as desired.

\textbf{Step 2:} We claim that
\begin{align}
\label{eq:K_1-estimate_pos-part}
\begin{split}
L_{K_1}(\eta^2)(x) (\partial_e v(x))_+^2 -& B_{K_1}(\eta^2, (\partial_e v)_+^2)(x) \le \eta^2(x) B_{K_1}((\partial_e v)_+ , (\partial_e v)_+)(x)\\
& - 2 \eta^2(x) (\partial_e v(x))_+ L_{K_1}((\partial_e v)_-)(x) + \sigma_1 B_K(v,v)(x),
\end{split}
\end{align}
where $\sigma_1 > 0$ is a constant.\\
To establish \eqref{eq:K_1-estimate_pos-part}, the same proof as for \eqref{eq:K_1-estimate_A} in Step 2 for \autoref{thm:keyest} goes through. One only needs to replace $\partial_e u$ by $(\partial_e v)_+$ and apply \autoref{lemma:interpol_pos-part} instead of \autoref{lemma:interpol} to deduce \eqref{eq:-BK_upper_bound}. Hence, we only need to prove
\begin{align*}
	L_{K_1}(\eta^2)(x)(\partial_e v(x))_+^2 &\le \frac{1}{2} \eta^2(x)B_{K_1}((\partial_e v)_+,(\partial_e v)_+)(x)\\
	&\quad - 2 \eta^2(x) (\partial_e v(x))_+ L_{K}((\partial_e v)_-)(x) + \sigma_{1,2} B(v,v)(x).
\end{align*}
If $\partial_e v(x) \le 0$ the estimate is trivial. Otherwise, following the computations in Step 2 using \autoref{lemma:interpol_pos-part} instead of \autoref{lemma:interpol} gives
\begin{align*}
	L_{K_1}(\eta^2)(x)(\partial_e v(x))_+^2 &\le \frac{1}{4} \eta^2(x)B_{K_1}((\partial_e v)_+,(\partial_e v)_+)(x)\\
	&\quad - \frac{1}{4} \eta^2(x) (\partial_e v(x))_+ L_{K}((\partial_e v)_-)(x) + \sigma_{1,2} B(v,v)(x),
\end{align*}
so, it suffices to note that $L_K((\partial_e v)_-)(x) < 0$ whenever $\partial_e v(x) > 0$.

Note that by combining \eqref{eq:K_1-estimate_pos-part} with \eqref{eq:genK_2-estimate_pos-part} and using that $\eps = \gamma \eta(x)$, we obtain that for every $x \in \R^n$:
\begin{align}
\label{eq:result_pos-part}
\begin{split}
L_K(\eta^2)(x) (\partial_e v(x))^2 -& B_K(\eta^2, (\partial_e v)^2)(x) \le \eta^2(x) B_K(\partial_e v , \partial_e v)(x)\\
& - 2 \eta^2(x) (\partial_e v(x))_+ L_{K}((\partial_e v)_-)(x) + \sigma B_K(v,v)(x),
\end{split}
\end{align}
where $\sigma = \sigma(n,s,\Lambda/\lambda,\Vert \eta \Vert_{C^{1,1}(\R^n)}) > 0$. This proves the desired result.
\end{proof}

Let us now prove a key estimate involving second derivatives. It is a straightforward corollary of \autoref{thm:keyest_pos-part}.

\begin{corollary}
\label{cor:keyest_snd-order}
Let $s\in(0,1)$, $L \in \mathcal{L}_s(\lambda,\Lambda;1)$. Let $\bar{\eta} \in C_c^{\infty}(B_{1/4})$ be such that $\bar{\eta} \equiv 1$ in $B_{1/8}$ and $0 \le \bar{\eta} \le 1$.  Then, there exists $\sigma_0 = \sigma_0(n, s, \lambda, \Lambda, \Vert \bar{\eta} \Vert_{C^{1,1}(\R^n)}) > 0$ such that for every $\sigma \ge \sigma_0$ and $\U \in C^{2+2s+\eps}_c(B_1)$ the following estimate holds true:
\begin{align}
\label{eq:keyest_snd-order}
    L(\bar{\eta}^2 (-\partial^2_{ee} \U)_+^2 + \sigma (\partial_e \U)^2) &\le 2\bar{\eta}^2 L(-\partial^2_{ee} \U)(-\partial^2_{ee} \U)_+ + 2 \sigma L(\partial_e \U) \partial_e \U \quad \textrm{in}\quad \R^n.
\end{align}
\end{corollary}

\begin{proof}
The result follows by application of \autoref{thm:keyest_pos-part} with $-\partial_e\U$.
\end{proof}

\section{Application to fully nonlinear equations}
\label{sec:fully-nonlinear}

The goal of this section is to establish \autoref{thm:intro-FN-sc}. The main tool in our proof is the Bernstein technique for integro-differential operators, which we develop in this article. 

First, let us state a more  general version of \autoref{thm:intro-FN-sc}. To do so, let us introduce the class of fully nonlinear operators $\mathcal{J}_s(\lambda,\Lambda)$:

\begin{definition}
We define $\mathcal{J}_s(\lambda,\Lambda)$ to be the set of all operators $\mathcal{I}$ of the form
\begin{align*}
\mathcal{I} u = \inf_{\gamma \in \Gamma} \left\{ L_{\gamma} u  - c_{\gamma} \right\},
\end{align*}
where $\Gamma$ is an index set and for any $\gamma \in \Gamma$ it holds $L_{\gamma} \in \mathcal{L}_s(\lambda,\Lambda;2)$, and $c_{\gamma} \in C^{1,1}$ satisfying
\begin{align*}
\sup_{\gamma \in \Gamma} \Vert c_{\gamma} \Vert_{C^{1,1}} \le \Lambda.
\end{align*}
\end{definition}

We prove the following.

\begin{theorem}
\label{thm:fully-nonlinear_one-sided}
Let $s\in(0,1)$ and $\mathcal{I} \in \mathcal{J}_s(\lambda,\Lambda)$. 
Let $u \in C(B_1) \cap L^{\infty}(\R^n)$ be a viscosity solution to
\begin{align}
\label{eq:fully-nonlinear}
\mathcal{I} u = 0 ~~ \text{ in}\quad B_1.
\end{align}
Then, in the convexity sense,
\begin{align}
\label{eq:fully-nonlinear_one-sided}
\inf_{B_{1/8}} \partial_{ee}^2 u \ge -C \left( \sup_{\gamma \in \Gamma} \Vert c_{\gamma} \Vert_{C^{1,1}(B_1)} + \Vert u \Vert_{L^{\infty}(\R^n)}\right),
\end{align}
where $C$ depends only on $n$, $s$, $\lambda$, $\Lambda$.
\end{theorem}

\begin{proof}
We start by proving that for any smooth solution $u\in C^\infty(B_1)\cap L^\infty(\R^n)$ of an equation of the type \eqref{eq:fully-nonlinear} we have 
\begin{align}
\label{eq:fully-nonlinear-approx_one-sided}
\sup_{B_{1/8}} -\partial_{ee}^2 u \le C \left( \sup_{\gamma \in \Gamma} \Vert \nabla c_{\gamma}\Vert_{L^\infty(B_1)}+\sup_{\gamma \in \Gamma} \Vert D^2 c_{\gamma}\Vert_{L^\infty(B_1)} + \Vert \partial_e u \Vert_{L^{\infty}(B_{1})} + \Vert u \Vert_{L^{\infty}(\R^n)} \right),
\end{align}
where $C = C(n,s,\lambda,\Lambda) > 0$.

For this, a key observation is that, by the regularity of $u$, for every $x \in B_1$ there exist $L_{\gamma(x)} \in \mathcal{L}_s(\lambda,\Lambda;2)$ and $c_{\gamma(x)} \in C^{1,1}$ with $\Vert c_{\gamma(x)} \Vert_{C^{1,1}} \le \Lambda$ such that 
\begin{align*}
L_{\gamma(x)} u(x)  = c_{\gamma(x)}(x).
\end{align*}
To see this, we take minimizing sequences $(L_n)$ and $(c_n)$ such that $-L_n u(x) + c_n(x) \searrow 0$. Note that by Arzel\`a-Ascoli, the sequence $(c_n)$ converges up to a subsequence to $c_{\gamma(x)}$ with the aforementioned properties and moreover, the sequence $(L_n)$ weakly convergences up to a subsequence to $L_{\gamma(x)}$ satisfying $c_{\gamma(x)}(x) = \lim_{n \to \infty} c_n(x) = \lim_{n \to \infty} L_n u(x) = L_{\gamma(x)}u(x)$.

\textbf{Step 1:} We claim that the following holds true in the classical sense:
\begin{align}
\label{eq:fully-nonlinear-int_interchange-derivative}
L_{\gamma(\cdot)} (\partial_e u) = \partial_e c_{\gamma(\cdot)}, ~~ L_{\gamma(\cdot)} (-\partial^2_{ee} u) \le -\partial^2_{ee} c_{\gamma(\cdot)} ~~ \text{ in } B_{3/4}.
\end{align}
Here, whenever we write $\partial_e c_{\gamma(\cdot)}$, we mean the partial derivative of the function $(x,y) \mapsto c_{\gamma(x)}(y)$ with respect to $y$, not with respect to $x$. The same convention applies to the expression $\partial_{ee}^2 c_{\gamma(\cdot)}$.
To prove \eqref{eq:fully-nonlinear-int_interchange-derivative}, we observe that since we have $L_{\gamma(\cdot)}u \ge c_{\gamma(\cdot)}$ in $B_1$, we find that for every small enough $h=te$, with $t\in\R$, and any $x \in B_{3/4}$
\begin{align*}
L_{\gamma(x)}(D_{h} u)(x) &= \frac{L_{\gamma(x)}u(x + h) - L_{\gamma(x)}u(x)}{|h|} \ge D_{h}c_{\gamma(x)}(x),\\
D_{h}c_{\gamma(x)}(x-h) &\ge \frac{L_{\gamma(x)}u(x - h) - L_{\gamma(x)}u(x)}{-|h|} = L_{\gamma(x)}(-D_{-h} u)(x),
\end{align*}
and 
\begin{align*}
&L_{\gamma(x)}(D_{h} D_{-h}u)(x) = \frac{2L_{\gamma(x)}u(x) - L_{\gamma(x)}u(x+h) - L_{\gamma(x)}u(x-h)}{|h|^2} \le D_{h} D_{-h}c_{\gamma(x)}(x).
\end{align*}
The aforementioned inequalities hold true in the classical sense since $u \in C^{\infty}(B_1)$. Therefore, \eqref{eq:fully-nonlinear-int_interchange-derivative} follows by taking the limit $|h| \to 0$.

\textbf{Step 2:} Next, take two cutoff functions $\eta \in C_c^{\infty}(B_1)$ and $\bar{\eta} \in C_c^{\infty}(B_{1/4})$ with $\eta, \bar{\eta} \ge 0$ and $\eta \equiv 1$ in $B_{1/2}$, $\bar{\eta} \equiv 1$ in $B_{1/8}$. By application of \autoref{cor:keyest_snd-order} to $\U = \eta u$, we derive that the following estimate holds true in a pointwise sense for any $x \in B_{1/4}$:
\begin{align}
\label{eq:fully-nonlinear-Bernstein-mollified}
\begin{split}
L_{\gamma(x)}(\bar{\eta}^2 & (-\partial^2_{ee} \U)_+^2 + \sigma (\partial_e \U)^2 )(x)\\
&\le 2\bar{\eta}^2 L_{\gamma(x)}(-\partial^2_{ee} \U)(-\partial^2_{ee} \U)_+(x) + 2 \sigma L_{\gamma(x)}(\partial_e \U) \partial_e \U(x).
\end{split}
\end{align}
Note that 
    \begin{align*}
        L_{\gamma(\cdot)}(\partial_e \U) = L_{\gamma(\cdot)}(\partial_e u) - L_{\gamma(\cdot)}(\partial_e[(1-\eta)u]), \quad L_{\gamma(\cdot)}(\partial_{ee}^2 \U) = L_{\gamma(\cdot)}(\partial^2_{ee} u) - L_{\gamma(\cdot)}(\partial^2_{ee}[(1-\eta)u]),
    \end{align*}
    and that we can estimate for any $x \in B_{1/4}$ using \eqref{eq:L1}
    \begin{align*}
        |L_{\gamma(x)}(\partial_e[(1-\eta)u])(x)| &= \left| \int_{\R^n \setminus B_{1/4}(x)} \partial_e[(1-\eta)u](y) K_{\gamma(x)}(x-y) \d y \right|\\
        &= \left| \int_{\R^n \setminus B_{1/4}(x)}  (1-\eta)u(y) \partial_e K_{\gamma(x)}(x-y)\d y \right|\\
        &\le C \Vert u \Vert_{L^{\infty}(\R^n)}.
    \end{align*}
    Note that also here, $\partial_e K_{\gamma(x)}(x-y)$ has to be interpreted as the partial derivative with respect to the argument of $K_{\gamma(x)}$.
    An analogous computation, using \eqref{eq:L2} instead of \eqref{eq:L1} yields
    \begin{align*}
        |L_{\gamma(x)}(\partial^2_{ee}[(1-\eta)u])(x)| \le C \Vert u \Vert_{L^{\infty}(\R^n)}.
    \end{align*}
    By combination of \eqref{eq:fully-nonlinear-Bernstein-mollified} with the previous observations and \eqref{eq:fully-nonlinear-int_interchange-derivative}, we derive 
    \begin{align}
    \label{eq:fully-nonlinear_F-estimate}
    \begin{split}
    L_{\gamma(x)}&(\bar{\eta}^2 (-\partial^2_{ee} \U)_+^2 + \sigma (\partial_e \U)^2)(x) \\
    &\le C \bar{\eta}^2 [L_{\gamma(x)}(-\partial^2_{ee} u)(x) + \Vert u \Vert_{L^{\infty}(\R^n)}](-\partial^2_{ee} u)_+(x) + C [L_{\gamma(x)}(\partial_e u) + \Vert u \Vert_{L^{\infty}(\R^n)}]|\partial_e u|(x)\\
    &\le C \bar{\eta}^2 [-\partial^2_{ee} c_\gamma|_{\gamma=\gamma(x)} + \Vert u \Vert_{L^{\infty}(\R^n)}](-\partial^2_{ee} u)_+(x) + C [\partial_e c_\gamma|_{\gamma=\gamma(x)} + \Vert u \Vert_{L^{\infty}(\R^n)}]|\partial_e u|(x)\\
    &=: F(x)
    \end{split}
    \end{align}

\textbf{Step 3:} Let us set 
\[\phi := \bar{\eta}^2 (-\partial^2_{ee} \U)_+^2 + \sigma (\partial_e \U)^2.\]
By the maximum principle  (see for instance \cite{CaSi09,CDV22}),
\begin{align*}
\sup_{B_{1/4}} \phi \le C \sup_{B_{1/4}} F + \sup_{\R^n \setminus B_{1/4}} \phi.
\end{align*}
Note that we can estimate by Young's inequality
\begin{align*}
C\sup_{B_{1/4}} F \le \frac{1}{2} \sup_{B_{1/4}} \phi + \tilde C \sup_{B_{1/4}} \left[ \bar{\eta}^2 (\partial^2_{ee} c_{\gamma(\cdot)})^2 + (\partial_e c_{\gamma(\cdot)})^2 + \Vert u \Vert^2_{L^{\infty}(\R^n)} \right].
\end{align*}
This yields 
\begin{align}
\label{eq:fully-nonlinear-max-principle}
\sup_{B_{1/4}} \phi \le C \left( \sup_{\gamma \in \Gamma} \Vert \nabla c_{\gamma}\Vert_{L^\infty(B_1)}^2+\sup_{\gamma \in \Gamma} \Vert D^2 c_{\gamma}\Vert_{L^\infty(B_1)}^2 + \Vert u \Vert_{L^{\infty}(\R^n)}^2 + \sup_{\R^n \setminus B_{1/4}} \phi \right).
\end{align}

By definition of $\bar{\eta}$ and $\eta$, we have
\begin{align*}
\sup_{B_{1/8}} (-\partial^2_{ee} u)_+^2 \le \sup_{B_{1/8}} \left[(-\partial^2_{ee} u)_+^2 + \sigma (\partial_e u)^2 \right] = \sup_{B_{1/8}} \phi \le \sup_{B_{1/4}} \phi.
\end{align*}
Moreover,
\begin{align*}
\sup_{\R^n \setminus B_{1/4}} \phi = \sup_{\R^n \setminus B_{1/4}} \left[ \sigma (\partial_e \U)^2 \right] \le \sigma \sup_{B_{1}}\left[ (\partial_e u)^2 + (\partial_e \eta)^2u^2 \right].
\end{align*}
By combination of the previous observations with \eqref{eq:fully-nonlinear-max-principle}, we find
\begin{align*}
\sup_{B_{1/8}} (-\partial^2_{ee} u)_+^2 &\le C \left( \sup_{\gamma \in \Gamma} \Vert c_{\gamma} \Vert_{C^{1,1}(B_{1/4})}^2 + \sup_{B_{1}}(\partial_e u)^2 + \sup_{\R^n} u^2 \right).
\end{align*}
This yields the desired estimate \eqref{eq:fully-nonlinear-approx_one-sided} in case $u \in C^{\infty}(B_1) \cap L^{\infty}(\R^n)$.

\textbf{Step 4:} Finally, let us prove that the same estimate holds for any viscosity solution $u$ of \eqref{eq:fully-nonlinear}.
Thanks to the results in \cite{Fer24}, for any viscosity solution $u$ of \eqref{eq:fully-nonlinear} there exists a sequence $u^{(k)} \in C^{\infty}(B_1) \cap L^{\infty}(\R^n)$ of solutions to the same class of equations converging to $u$ locally uniformly in $B_1$ and satisfying a global $L^{\infty}$-estimate. 
Moreover, due to \cite{CaSi11a}, the solutions $u^{(k)}$ satisfy a Lipschitz regularity estimate. Thanks to this, and since we already proved \eqref{eq:fully-nonlinear-approx_one-sided} for $C^\infty$-solutions, we obtain
\begin{align*}
\sup_{B_{1/8}} -\partial_{ee}^2 u^{(k)} &\le C_1 \left(\sup_{\gamma \in \Gamma^{(k)}} \Vert c_{\gamma} \Vert_{C^{1,1}(B_1)} + \Vert \partial_e u^{(k)} \Vert_{L^{\infty}(B_{1})} + \Vert u^{(k)} \Vert_{L^{\infty}(\R^n)} \right)\\
&\le C_2 \left(\sup_{\gamma \in \Gamma} \Vert c_{\gamma} \Vert_{C^{1,1}(B_1)} + \Vert u \Vert_{L^{\infty}(\R^n)}\right),
\end{align*}
where $C_1, C_2 > 0$ do not depend on $k$. 
Therefore, we can take the limit $k \to \infty$, and deduce the desired result.
\end{proof}

\begin{remark}
	The semiconvexity estimate \eqref{eq:fully-nonlinear_one-sided} is robust with respect to the limit $s \nearrow 1$, i.e., if we choose $L \in \mathcal{L}_s((1-s)\lambda,(1-s)\Lambda;1)$ and $s \ge s_0$ for some $s_0 \in (0,1)$, then the constant $C$ will depend only on $n, s_0, \lambda, \Lambda$.
\end{remark}

\begin{remark}
Instead of \eqref{eq:L2}, our proof of the semiconvexity estimate remains true under 
\begin{align}
\label{eq:weakL2}
\int_{\R^n \setminus B_{1/8}} |D^2K(y)| \d y \le \Lambda.
\end{align}
\end{remark}

\section{Key estimates in terms of difference quotients}
\label{sec:diffquot}

Recall that the key estimates which were established in the previous section all require some a priori smoothness assumption on $u$, such as $u \in C^{2+2s+\eps}_{loc}$ in order to make sense of expressions like $L(\partial_e u)$ and $L(\partial^2_{ee}u)$. This way, the Bernstein technique cannot be used in order to establish smoothness of solutions, but only to prove estimates for solutions that are already known to be smooth.\\
Moreover, in many applications, as for example the obstacle problem, solutions are known not to possess the required $C^{2+2s+\eps}$ regularity. In order to apply our technique to such situations, in this section we establish Bernstein key estimates for difference quotients. Let us mention that this idea was also announced in a version of a preprint of \cite{CDV22} (see Section 3 in \cite{CDV20}).

For $h \in \R^n$, and $u \in C(\R^n)$ let us define
\begin{align*}
u_{h}(x) := \int_0^1 u(x + th) \d t.
\end{align*}

Moreover, we introduce difference quotients
\begin{align*}
D_{h}u(x) := \frac{u(x+h) - u(x)}{|h|} = \partial_{e} (u_h)(x),
\end{align*}
where $e = h/|h|$.
Clearly, by the Lebesgue differentiation theorem, the following identity holds:
\begin{align*}
D_{h}u(x) = |h|^{-1} \lim_{\eps \to 0} \left( \dashint_{|h|}^{|h|+\eps} u(x + te) \d t  - \dashint_0^{\eps} u(x + te) \d t \right) =\partial_e(u_{h})(x).
\end{align*}
Moreover, note that
\begin{align*}
    D_{-h} D_{h}u(x) = \frac{2u(x) - u(x+h) - u(x-h)}{|h|^2}.
\end{align*}

\begin{remark}
    Note that we have the following integration by parts identity for difference quotients, whenever $f,g$ are such that the integrals below are well-defined:
    \begin{align*}
        \int_{\R^n} D_{h} f(x) g(x) \d x = \int_{\R^n} f(x) D_{-h} g(x) \d x.
    \end{align*}
\end{remark}

\subsection{Key estimate for non-smooth functions}

We claim the following analogues of \autoref{thm:keyest} and \autoref{thm:keyest_pos-part} for difference quotients:

\begin{proof}[Proof of \autoref{lemma:keyest_diff-quot}]
The proofs of \eqref{eq:keyest_diff-quot} and \eqref{eq:keyest_diff-quot_pos-part} go in the exact same way as the proofs of \autoref{thm:keyest} and \autoref{thm:keyest_pos-part}, respectively, upon replacing $u$ and $v$ by $u_{h}$ and $v_{h}$, respectively, and $\partial_e u$ and $\partial_e v$ by $D_{h}u$ and $D_{h}v$, respectively. Moreover, in order to prove \eqref{eq:keyest_diff-quot}, we need the following interpolation estimate, reminiscent of \autoref{lemma:interpol}:
\begin{align}
\label{eq:discrete-interpolation}
(D_{h} u(x))^2 \le \delta^{2s} B(D_{h} u, D_{h} u)(x) + c\delta^{2s-2} B(u_{h},u_{h})(x).
\end{align}
For the proof of \eqref{eq:keyest_diff-quot_pos-part} we have the following interpolation estimate:
\begin{align}
\label{eq:discrete-interpolation-2}
\begin{split}
(D_{h} v(x))_+^2 &\le \delta^{2s} B((D_{h} v)_+, (D_{h} v)_+)(x) \\
&\quad - \delta^{2s} L((D_h v)_-)(x)(D_{h}v)_+(x) + c\delta^{2s-2} B(v_{h},v_{h})(x).
\end{split}
\end{align}
The proof of the discrete interpolation estimate \eqref{eq:discrete-interpolation} also goes as before, however, the following computation has to be employed: if we denote $e=h/|h|$ then
\begin{align}
\label{eq:discrete-ibp}
\begin{split}
\int_{\R^n} D_{h} u(y) K_{\delta}(x-y) \d y &= \int_{\R^n} \partial_e (u_h(y) - u_{h}(x)) K_{\delta}(x-y) \d y\\
&= - \int_{\R^n} (u_{h}(y) - u_{h}(x)) \partial_e K_{\delta}(x-y) \d y,
\end{split}
\end{align}
using $D_hu = \partial_e (u_h)$ and integrating by parts. The same modification yields \eqref{eq:discrete-interpolation-2}.
Having at hand the discrete interpolation estimates \eqref{eq:discrete-interpolation} and \eqref{eq:discrete-interpolation-2}, the terms involving $K_1$ can be treated exactly as before. For $K_2$, we also proceed as in Step 1 of the proof of \autoref{thm:keyest}, employing a similar computation as \eqref{eq:discrete-ibp} with respect to $K_2$ instead of $K_{\delta}$.
\end{proof}

As in the non-discrete case, it is possible to use \eqref{eq:keyest_diff-quot_pos-part} in order to get a key estimate that is suitable for the application to second order estimates. 

\begin{corollary}
\label{cor:keyest_diff-quot_snd-order}
Let $s\in(s_0,1)$, with $s_0>0$, and $L \in \mathcal{L}_s(\lambda,\Lambda;1)$. Let $\bar{\eta} \in C_c^{1,1}(B_{1/4})$ be such that $\bar{\eta} \equiv 1$ in $B_{1/8}$ and $0 \le \bar{\eta} \le 1$. Let $|h| \le 1/8$. Then, there exists $\sigma_0 = \sigma_0(n, s_0, \Lambda / \lambda, \Vert \bar{\eta} \Vert_{C^{1,1}(\R^n)}) > 0$ such that for every $\sigma \ge \sigma_0$ and $u \in C^{2s+\eps}(B_1)$ the following estimate holds true
\begin{align}
\label{eq:keyest_diff-quot_snd-order}
\begin{split}
L(\bar{\eta}^2 & (D_{-h} D_{h} u)_+^2 + \sigma (D_{h} u_{-h})^2)\\
&\le 2\bar{\eta}^2 L(D_{-h} D_{h} u)(D_{-h} D_{h} u)_+ + 2 \sigma L(D_{h} u_{-h}) D_{h} u_{-h}.
\end{split}
\end{align}
\end{corollary}

\begin{proof}
    The proof follows by application of \eqref{eq:keyest_diff-quot_pos-part} with $-h$ and $v = D_{h} u$ and using
    \begin{align*}
        D_{h} (u_{-h}) = \frac{1}{h} \int_0^1 (u(\cdot - th + h) - u(\cdot - th)) \d t = (D_{h} u)_{-h}.
    \end{align*} 
\end{proof}

Finally, we observe that the key estimates \eqref{eq:keyest_diff-quot} and \eqref{eq:keyest_diff-quot_pos-part} can also be obtained for H\"older difference quotients defined as
\begin{align*}
D^{\alpha}_{h}u(x) = \frac{u(x+h) - u(x)}{|h|^{\alpha}}, ~~ \alpha \in (0,1).
\end{align*}

\begin{corollary}
\label{cor:keyest_Holder-diff-quot}
Let $s\in(s_0,1)$, with $s_0>0$, and $L \in \mathcal{L}_s(\lambda,\Lambda;1)$. Let $\eta \in C_c^{1,1}(B_1)$ be such that $\eta \equiv 1$ in $B_{1/2}$ and $0 \le \eta \le 1$.  Let $\alpha \in (0,1)$ and $|h| \le 1/8$. Then, there exists $\sigma_0 = \sigma_0(n, s_0, \Lambda / \lambda, \Vert \eta \Vert_{C^{1,1}(\R^n)}) > 0$ such that for every $\sigma \ge \sigma_0$ and $u,v \in C^{2s+\eps}_{loc}(\R^n) \cap L^{\infty}(\R^n)$:
\begin{align}
\label{eq:keyest_Holder-diff-quot}
L(\eta^2 (D^{\alpha}_{h} u)^2 + \sigma u_{h}^2) &\le 2 \eta^2 L(D^{\alpha}_{h}u) D^{\alpha}_{h}u + 2\sigma L(u_{h}) u_{h},\\
\label{eq:keyest_Holder-diff-quot_pos-part}
L(\eta^2 (D^{\alpha}_{h} v)_+^2 + \sigma v_{h}^2) &\le 2 \eta^2 L(D^{\alpha}_{h}v) (D^{\alpha}_{h}v)_+ + 2\sigma L(v_{h}) v_{h}.
\end{align}
\end{corollary}

\begin{proof}
We only explain how to prove \eqref{eq:keyest_Holder-diff-quot}, since the proof of \eqref{eq:keyest_Holder-diff-quot_pos-part} can be proved in the same way. We multiply \eqref{eq:keyest_diff-quot} on both sides by $|h|^{2-2\alpha}$. Then, for any $\sigma \ge \sigma_0$
\begin{align*}
L(\eta^2 (D^{\alpha}_{h} u)^2 + |h|^{2-2\alpha} \sigma u_{h}^2) &\le 2 \eta^2 L(D^{\alpha}_{h}u) D^{\alpha}_{h}u + 2|h|^{2-2\alpha} \sigma L(u_{h}) u_{h}.
\end{align*}
In particular, we obtain \eqref{eq:keyest_Holder-diff-quot} for any $\sigma \ge \sigma_0 \ge \sigma_0 |h|^{2-2\alpha}$.
\end{proof}

\subsection{Improved key estimate for Lipschitz continuous functions}
Moreover, we can prove the following key estimate which produces slightly different averages on the lower order terms. The price we have to pay is that we need to assume certain smoothness of $u$:

\begin{lemma}
\label{lemma:keyest_diff-quot2}
Let $s\in(s_0,1)$, with $s_0>0$, and $L \in \mathcal{L}_s(\lambda,\Lambda;1)$. Let $\eta \in C_c^{1,1}(B_1)$ be such that $\eta \ge 0$. Let $|h| \le 1/8$. Then, there exists $\sigma_0 = \sigma_0(n, s_0, \Lambda / \lambda, \Vert \eta \Vert_{C^{1,1}(\R^n)}) > 0$ such that for every $\sigma \ge \sigma_0$ and $u,v \in C^{2s+\eps}(B_{1/2}) \cap C^{0,1}(B_{1/2}) \cap L^{\infty}(\R^n)$, the following estimate holds true in $B_{1/4}$:
\begin{align}
\label{eq:keyest_diff-quot2}
L(\eta^2 (D_{h} u)^2 + \sigma [u^2]_{h}) &\le 2 \eta^2 L(D_{h}u) D_{h}u + 2\sigma [L(u) u]_{h},\\
\label{eq:keyest_diff-quot2_pos-part}
L(\eta^2 (D_{h} v)_+^2 + \sigma [v^2]_{h}) &\le 2 \eta^2 L(D_{h}v) (D_{h}v)_+ + 2\sigma [L(v) v]_{h}.
\end{align}
\end{lemma}

\begin{proof}
The proof goes exactly as the proof of \autoref{lemma:keyest_diff-quot}. The only difference is how we integrate by parts with difference quotients. This is contained in the following identity, which replaces \eqref{eq:discrete-ibp}. We claim that, denoting $e=h/|h|$,
\begin{align}
\label{eq:discrete-ibp2-claim}
\int_{\R^n} D_{h} v(y) J(x-y) \d y = \left( \int_{\R^n} \left(v(\cdot) - v(y) \right) \partial_e J(\cdot - y) \d y \right)_{h}(x), ~~ x \in B_{1/4},
\end{align}
holds true for $J = K_2$ (as in \autoref{thm:keyest}) and $J = K_{\delta}$ (as in \autoref{lemma:interpol}), where $v$ is Lipschitz continuous in $B_{1/2}$ and bounded. 
Let us first verify \eqref{eq:discrete-ibp2-claim} for globally Lipschitz continuous $v$:
\begin{align}
\label{eq:discrete-ibp2}
\begin{split}
\int_{\R^n} D_{h} v(y) J(x-y) \d y &= \int_{\R^n} \left[ \int_0^1 \partial_e v(y + th) \d t \right] J(x-y) \d y \\
&=  \int_{\R^n} \left[ \int_0^1 \partial_e \left(v(y + th) - v(x+th) \right) \d t \right] J(x-y) \d y \\
&=  \int_0^1 \int_{\R^n} \partial_e \left(v(y + th) - v(x+th) \right) J(x-y) \d y  \d t\\
&= - \int_0^1 \int_{\R^n} \left(v(y + th) - v(x+th) \right) \partial_e J(x-y) \d y  \d t\\
&= - \int_0^1 \int_{\R^n} \left(v(y) - v(x+th) \right) \partial_e J((x+th) - y) \d y  \d t\\
&= \left( \int_{\R^n} \left(v(\cdot) - v(y) \right) \partial_e J(\cdot - y) \d y \right)_{h}(x).
\end{split}
\end{align}
Note that we used Lipschitz continuity of $v$ in the first step. Moreover, we used $L([w]_{h}) = [L w]_{h}$, which is a direct consequence of Fubini's theorem.

We remark that in case $J = K_{\delta}$, it suffices to have Lipschitz continuity of $v$ in $B_{1/2}$ in order to carry out the computation in \eqref{eq:discrete-ibp2}, because $\supp(K_{\delta}) \subset B_{\delta}(0)$ and $\delta > 0$ can be chosen to satisfy $\delta < 1/8$ in the application. This yields the following discrete interpolation estimate:
\begin{align}
\label{eq:discrete-interpolation2}
(D_{h} u(x))^2 \le \delta^{2s} B(D_{h} u, D_{h} u)(x) + c\delta^{2s-2} [B(u,u)]_{h}(x),
\end{align}
which allows us to treat all the terms containing $K_1$.

In order to prove \eqref{eq:discrete-ibp2-claim} with $J = K_2$, we approximate any  function $v \in L^\infty(\R^n)$ that is Lipschitz continuous in $B_{1/2}$ by a sequence of globally Lipschitz continuous functions $(v_k)$, satisfying $\Vert v_k \Vert_{L^{\infty}(\R^n)} \le 2 \Vert v \Vert_{L^{\infty}(\R^n)}$, $v_k \equiv v$ in $B_{1/4}$, and converging to $v$ almost everywhere, as $k \to \infty$.  Then, by dominated convergence, for almost every $x \in \R^n$:
\begin{align*}
\int_{\R^n} D_{h} v(y) K_2(x-y) \d y &= \lim_{k \to \infty} \int_{\R^n} D_{h} v_k(y) K_2(x-y) \d y\\
&= \lim_{k \to \infty} \left( \int_{\R^n} \left(v_k(y) - v_k(\cdot) \right) \partial_e K_2(\cdot - y) \d y \right)_{h}(x)\\
&= \left( \int_{\R^n} \left(v(y) - v(\cdot) \right) \partial_e K_2(\cdot - y) \d y \right)_{h}(x),
\end{align*} 
where we used that $K_2$ and $|\nabla K_2|$ are integrable by properties $(iv)$, $(v)$, and $(vi)$ in \autoref{lemma:kernel-decomp}, and that $D_{h} v_k \to D_{h} v$ a.e. and $v_k(y) - v_k(x + th) \to v(y) - v(x+th)$ for a.e. $y \in \R^n$ and every $x \in B_{1/4}$. This proves \eqref{eq:discrete-ibp2-claim} for $J = K_2$ and allows us to treat all terms involving $K_2$ in the proof of the key estimate. This concludes the proof.
\end{proof}

Finally, we explain how to obtain second order estimates from \eqref{eq:keyest_diff-quot2_pos-part}:

\begin{corollary}
\label{cor:keyest_diff-quot2_snd-order}
Let $s\in(s_0,1)$, with $s_0>0$, and $L \in \mathcal{L}_s(\lambda,\Lambda;1)$. Let $\bar{\eta} \in C_c^{1,1}(B_{1/4})$ be such that $\bar{\eta} \equiv 1$ in $B_{1/8}$ and $0 \le \bar{\eta} \le 1$. Let $|h| \le 1/8$. Then, there exists $\sigma_0 = \sigma_0(n, s_0, \Lambda / \lambda, \Vert \bar{\eta} \Vert_{C^{1,1}(\R^n)}) > 0$ such that for every $\sigma \ge \sigma_0$ and $u \in C^{2s+\eps}(B_{1/2}) \cap C^{0,1}(B_{1/2}) \cap L^{\infty}(B_{1})$, the following estimate holds true in $B_{1/4}$:
\begin{align}
\label{eq:keyest_diff-quot2_snd-order}
\begin{split}
L(\bar{\eta}^2 &(D_{-h} D_{h} u)_+^2 + \sigma [(D_{h} u)^2]_{-h})\\
&\le 2\bar{\eta}^2 L(D_{-h} D_{h} u)(D_{-h} D_{h} u)_+ + 2 \sigma [ L(D_{h} u) D_{h} u]_{-h}.
\end{split}
\end{align}
\end{corollary}

\begin{proof}
    The proof follows by application of \eqref{eq:keyest_diff-quot2_pos-part} with $-h$ and $v = D_{h} u$.
\end{proof}

As in the previous section, we can also obtain an improved key inequality for H\"older difference quotients:

\begin{corollary}
\label{cor:keyest_Holder-diff-quot2}
Let $s\in(s_0,1)$, with $s_0>0$, and $L \in \mathcal{L}_s(\lambda,\Lambda;1)$. Let $\eta \in C_c^{1,1}(B_1)$ be such that $\eta \equiv 1$ in $B_{1/2}$ and $0 \le \eta \le 1$.  Let $\alpha \in (0,1)$ and $|h| \le 1/8$. Then, there exists $\sigma_0 = \sigma_0(n, s_0, \Lambda / \lambda, \Vert \eta \Vert_{C^{1,1}(\R^n)}) > 0$ such that for every $\sigma \ge \sigma_0$ and $u,v \in C^{2s+\eps}(B_{1/2}) \cap C^{0,1}(B_{1/2}) \cap L^{\infty}(\R^n)$, the following estimate holds true in $B_{1/4}$:
\begin{align}
\label{eq:keyest_Holder-diff-quot2}
L(\eta^2 (D^{\alpha}_{h} u)^2 + \sigma [u^2]_{h}) &\le 2 \eta^2 L(D^{\alpha}_{h}u) D^{\alpha}_{h}u + 2\sigma [L(u) u]_{h},\\
\label{eq:keyest_Holder-diff-quot2_pos-part}
L(\eta^2 (D^{\alpha}_{h} v)_+^2 + \sigma [v^2]_{h}) &\le 2 \eta^2 L(D^{\alpha}_{h}v) (D^{\alpha}_{h}v)_+ + 2\sigma [L(v) v]_{h}.
\end{align}
\end{corollary}

\begin{proof}
The proof goes exactly as the proof of \autoref{cor:keyest_Holder-diff-quot}, but starting with \eqref{eq:keyest_diff-quot2} and \eqref{eq:keyest_diff-quot2_pos-part}.
\end{proof}

\section{Application to the obstacle problem}
\label{sec:obstacle}

In this section we apply the Bernstein technique to the obstacle problem \eqref{eq:intro-OP2} with  $L \in \mathcal{L}_s(\lambda,\Lambda;1)$ (see \autoref{def:classL}) and with $\phi$ being a sufficiently smooth function.

First, we establish semiconvexity estimates for solutions to \eqref{eq:intro-OP2} in Subsection \ref{subsec:Bernstein_OP}. In Subsection \ref{subsec:blowups-convex} we prove that blow-ups are convex (see \autoref{thm:blowups-convex}). This result is crucial in the proof of our main results \autoref{thm:OP_bdry-reg} and \autoref{thm:OP_opt-reg}, which will be carried out in Subsection \ref{subsec:OP-main}.

Let us start our discussion by several remarks.

\begin{remark}
It was proved in \cite[Theorem 5.1]{CDS17} that any solution $u$ to the obstacle problem \eqref{eq:intro-OP2} belongs to $C^{\beta}$ whenever $\phi \in C^{\beta}$ for some $\beta <\max\{1+\epsilon,  2s+\epsilon\}$. 
In particular, if $\beta > 2s$, we may assume that $u$ is a solution to \eqref{eq:intro-OP2} in the classical sense.
\end{remark}

\begin{remark}\label{remark-change}
Let us also point out that once $\phi \in C^{\beta}$ for some $\beta > 2s$, we can rewrite \eqref{eq:intro-OP2} as an obstacle problem with a zero obstacle and an inhomogeneity, as follows
\begin{align}
\label{eq:OP}
\min\{Lv -f , v\} = 0 ~~ \text{ in } B_1.
\end{align}
This can be achieved by extending $\phi$ in a smooth way to $\R^n$ and defining $f= L \phi \in C^{\beta - 2s}$ and $v = u-\phi$ (as long as $\beta-2s$ is not an integer). 
Sometimes it will be more convenient to work with the formulation \eqref{eq:OP} instead of \eqref{eq:intro-OP2}.
\end{remark}

\subsection{Semiconvexity estimates}
\label{subsec:Bernstein_OP}

Using the Bernstein technique for difference quotients (see Section \ref{sec:diffquot}), we prove the semiconvexity estimate in  \autoref{thm:intro-OP-sc}. 
Note that we will not apply these results in the proof of optimal regularity for solutions and regularity of the free boundary. However, we consider both, the result in itself and its proof using Bernstein technique, of independent interest.

In order to apply the Bernstein technique, we need to be able to evaluate $L u$ and $L D_{h} u$ in a pointwise sense. We rely on Theorem 5.1 in \cite{CDS17}, which states that solutions to the obstacle problem are classical once the data is smooth enough.

\begin{proof}[Proof of \autoref{thm:intro-OP-sc}]
Thanks to \autoref{remark-change}, we may assume that $u$ is a solution of \eqref{eq:OP}, with $f=L\phi$.\\
Let us denote $\U = \eta u$, where $\eta \in C_c^{\infty}(B_{3/4})$, $\eta \equiv 1$ in $B_{1/2}$ and $0 \le \eta \le 1$. Let $0 < |h| < 1/16$.

\textbf{Step 1:} We estimate the quantity $[L(D_{h} \U)D_{h}\U]_{-h}$ in $B_{1/4 - |h|}$. Let us first prove the following claim: For all $x \in B_{1/4}$,
\begin{align}
\label{eq:application_OP_first-order1}
L(D_{h} u)(x)D_{h}u(x) \le \Vert \nabla f \Vert_{L^{\infty}(B_1)}  |D_{h}u(x)|.
\end{align}
To prove \eqref{eq:application_OP_first-order1}, we first assume that $u(x) > 0$. We distinguish between two cases:\\
Case 1: $u(x+h) > 0$. Then, $L(D_{h} u)(x) = D_{h}f(x)$, and therefore
\begin{align*}
L(D_{h} u)(x) D_{h}u(x) = D_{h}f(x)D_{h}u(x) \le \Vert \nabla f \Vert_{L^{\infty}(B_1)}  |D_{h}u(x)|.
\end{align*}
Case 2: $u(x+h) = 0$. Then, 
\begin{align*}
D_{h}u(x) = - \frac{u(x)}{h} < 0, ~~ L(D_{h} u)(x) = \frac{Lu(x+h) - f(x)}{h} \ge  D_{h} f(x).
\end{align*}
Thus,
\begin{align*}
L(D_{h} u)(x)D_{h}u(x) \le D_{h} f(x)D_{h}u(x)\le \Vert \nabla f \Vert_{L^{\infty}(B_1)}  |D_{h}u(x)|.
\end{align*}
On the other hand, let us now assume that $u(x) = 0$. We distinguish between two cases:\\
Case 1: $u(x+h) = 0$. Then, $D_{h}u(x) = 0$, and therefore
\begin{align*}
L(D_{h} u)(x)D_{h}u(x) = 0.
\end{align*}
Case 2: $u(x+h) > 0$. Then,
\begin{align*}
D_{h}u(x) = \frac{u(x+h)}{h} > 0, ~~ L(D_{h}u)(x) = \frac{f(x+h) - Lu(x)}{h} \le D_{h}f(x).
\end{align*}
Thus,
\begin{align*}
L(D_{h} u)(x)D_{h}u(x) \le D_{h} f(x)D_{h}u(x)\le \Vert \nabla f \Vert_{L^{\infty}(B_1)}  |D_{h}u(x)|.
\end{align*}
All in all, we get for any $x \in B_1$:
\begin{align*}
L(D_{h} u)(x)D_{h}u(x) \le D_{h} f(x)D_{h}u(x)\le \Vert \nabla f \Vert_{L^{\infty}(B_1)}  |D_{h}u(x)|,
\end{align*}
which yields \eqref{eq:application_OP_first-order1}, as desired.\\
Now let us turn to estimating $[L(D_{h} \U)D_{h}\U]_{-h}$ in $B_{1/4 - |h|}$. For this, we note that for $x \in B_{1/4}$, using \eqref{eq:L1}:
\begin{align*}
    |L(D_{h}[(1-\eta)u](x)| &= \left| \int_{\R^n \setminus B_{\frac{1}{4} - |h|}(x)} D_{h}[(1-\eta)u](y) K(x-y) \d y \right| \\
    &= \left| \int_{\R^n \setminus B_{\frac{1}{4} - |h|}(x)} (1-\eta)u(y) D_{-h}K(x-y) \d y \right| \\
    &\le C \Vert u \Vert_{L^{\infty}(\R^n)}.
\end{align*}
Therefore, we obtain
\begin{align*}
L(D_{h} \U)D_{h}\U(x) &= L(D_{h} u)D_{h}u - L(D_{h}[1-\eta]u)D_{h}u\\
&\le C \left(\Vert \nabla f \Vert_{L^{\infty}(B_1)} + \Vert u \Vert_{L^{\infty}(\R^n)} \right) |D_{h}u(x)|
\end{align*}
Using interior Lipschitz estimates for solutions of \eqref{eq:OP} (see \cite[Theorem 5.1]{CDS17}), we obtain the following estimate in $B_{1/4-|h|}$:
\begin{align}
\label{eq:application_OP_snd-order2}
\begin{split}
[L(D_{h} \U) D_{h}\U]_{-h} &\le C \left( \Vert \nabla f \Vert_{L^{\infty}(B_1)} + \Vert u \Vert_{L^{\infty}(\R^n)} \right)(\Vert f \Vert_{C^{0,1}(B_1)} + \Vert u \Vert_{L^{\infty}(\R^n)})\\
&\le C(\Vert f \Vert_{C^{0,1}(B_1)}^2 + \Vert u \Vert_{L^{\infty}(\R^n)}^2).
\end{split}
\end{align}

\textbf{Step 2:} We give an estimate for $L (D_{-h}D_{h}\U)(x) (D_{-h}D_{h}\U(x))_+$ in $B_{1/4}$:
First of all, let us assume that $u(x) > 0$. In this case
\begin{align*}
L (D_{-h}D_{h}u)(x) &= \frac{2 Lu(x) - L u(x+h) - Lu(x-h)}{h^2}\\
&= \frac{2 f(x) - L u(x+h) - Lu(x-h)}{h^2} \le D_{-h}D_{h}f(x).
\end{align*}
Moreover, since in $\{ u = 0 \}$, it holds $D_{-h}D_{h}u \le 0$, we obtain for $x \in B_1$:
\begin{align}
\label{eq:application_OP_snd-order1}
L (D_{-h}D_{h}u)(x) (D_{-h}D_{h}u(x))_+ \le \Vert D^2 f \Vert_{L^{\infty}(B_1)} (D_{-h}D_{h}u(x))_+.
\end{align}
Next, note that, by a similar argument as in Step 1, but using \eqref{eq:L2} instead of \eqref{eq:L1}, we obtain for $x \in B_{1/4 -2|h|}$:
\begin{align}
\label{eq:application_OP_snd-order5}
    |L (D_{-h}D_{h}[(1-\eta)u])(x)| \le C \Vert u \Vert_{L^{\infty}(\R^n)}.
\end{align}
Therefore, using \eqref{eq:application_OP_snd-order1} and \eqref{eq:application_OP_snd-order5}, we have the following estimate:
\begin{align}
\label{eq:application_OP_snd-order3}
\begin{split}
L&(D_{-h}D_{h}\U) (D_{-h}D_{h}\U)_+(x)\\
&= L (D_{-h}D_{h}u) (D_{-h}D_{h}u(x))_+ - L (D_{-h}D_{h}[(1-\eta)u]) (D_{-h}D_{h}u(x))_+\\
&\le C \left(\Vert D^2 f \Vert_{L^{\infty}(B_1)} + \Vert u \Vert_{L^{\infty}(\R^n)} \right)(D_{-h}D_{h}u(x))_+
\end{split}
\end{align}

\textbf{Step 3:} Finally, it remains to apply \autoref{cor:keyest_diff-quot2_snd-order} to $\U$ and to combine it with the estimates \eqref{eq:application_OP_snd-order2} and \eqref{eq:application_OP_snd-order3} we obtained in the previous steps. 
Note that we can apply \autoref{cor:keyest_diff-quot2_snd-order}, since we know that $u$ is Lipschitz continuous in $B_{1/2}$.\\
Let $\bar{\eta} \in C_c^{\infty}(B_{1/4})$ be such that $\bar{\eta} \equiv 1$ in $B_{1/8}$ and $0 \le \bar{\eta} \le 1$. This yields the following estimate for every $x \in B_{1/4}$, if $|h|$ is small enough:
\begin{align*}
\begin{split}
L(\bar{\eta}^2& (D_{-h} D_{h} \U)_+^2 + \sigma [(D_{h} \U)^2]_{-h})(x)\\
&\le 2\bar{\eta}^2(x) L(D_{-h} D_{h} \U)(x)(D_{-h} D_{h} \U)_+(x) + 2 \sigma [L(D_{h} \U) D_{h} \U]_{-h}(x)\\
&\le C \bar{\eta}^2(x) \left(\Vert D^2 f \Vert_{L^{\infty}(B_1)} + \Vert u \Vert_{L^{\infty}(\R^n)} \right)(D_{-h}D_{h}u(x))_+ +C(\Vert f \Vert_{C^{0,1}(B_1)}^2 + \Vert u \Vert_{L^{\infty}(\R^n)}^2).
\end{split}
\end{align*}
Consequently, by the maximum principle (see Corollary 5.2 in \cite{Ros16}), we have
\begin{align}
\label{eq:obstacle-proof-reference}
\begin{split}
\sup_{B_{1/4}} \left( \bar{\eta}^2 (D_{-h} D_{h} u)_+^2 \right) &\le \sup_{B_{1/4}} \left( \bar{\eta}^2 (D_{-h} D_{h} \U)_+^2 + \sigma [(D_{h} \U)^2]_{-h}\right)\\
&\le C \sup_{B_{1/4}} \left(\bar{\eta}^2 \left(\Vert D^2 f \Vert_{L^{\infty}(B_1)} + \Vert u \Vert_{L^{\infty}(\R^n)} \right)(D_{-h}D_{h}u)_+ \right)\\
&\quad + C(\Vert f \Vert_{C^{0,1}(B_1)}^2 + \Vert u \Vert_{L^{\infty}(\R^n)}^2)\\
&\quad + C \sup_{\R^n \setminus B_{1/4}} \left( \bar{\eta}^2 (D_{-h} D_{h} \U)_+^2 + \sigma [(D_{h} \U)^2]_{-h} \right)\\
&\le \frac{1}{2}\sup_{B_{1/4}} \left( \bar{\eta}^2 (D_{-h} D_{h} u)_+^2 \right) + C (\Vert f \Vert_{C^{1,1}(B_1)}^2 + \Vert u \Vert_{L^{\infty}(\R^n)}^2)\\
&\quad + C \sup_{\R^n \setminus B_{1/4}}  [(D_{h} \U)^2]_{-h},
\end{split}
\end{align}
where we used Young's inequality. By absorption of the first term on the right hand side, as well as 
\begin{align*}
\sup_{\R^n \setminus B_{1/4}} [(D_{h} \U)^2]_{-h} \le C \Vert \nabla u \Vert_{L^{\infty}(B_{7/8})}^2 + C \Vert u \Vert_{L^{\infty}(\R^n)}^2 \le C(\Vert f \Vert_{C^{0,1}(B_1)}^2 + \Vert u \Vert_{L^{\infty}(\R^n)}^2),
\end{align*}
which is again due to the interior Lipschitz estimates in \cite[Theorem 5.1]{CDS17}, we end up with
\begin{align}
\label{eq:application_OP_snd-order4}
\begin{split}
\sup_{B_{1/8}} \left((D_{-h} D_{h} u)_+^2 \right) \le \sup_{B_{1/4}} \left( \bar{\eta}^2 (D_{-h} D_{h} u)_+^2 \right) \le C(\Vert f \Vert_{C^{1,1}(B_1)}^2 + \Vert u \Vert_{L^{\infty}(\R^n)}^2).
\end{split}
\end{align}
Since this estimate is uniform in $h$, we can pass to the limit $h \to 0$ and deduce 
\begin{align*}
\sup_{B_{1/8}} (-\partial^2_{ee} u)_+ \le C (\Vert f \Vert_{C^{1,1}(B_1)} + \Vert u \Vert_{L^{\infty}(\R^n)}).
\end{align*}
This concludes the proof.
\end{proof}

\begin{remark}
The semiconvexity estimate we just proved is robust with respect to the limit $s \nearrow 1$, i.e., if we choose $L \in \mathcal{L}_s((1-s)\lambda,(1-s)\Lambda;1)$ and $s \ge s_0$ for some $s_0 \in (0,1)$, then the constant $C$ will depend only on $n, s_0, \lambda, \Lambda$
\end{remark}

\begin{remark}
Note that our proof of the semiconvexity estimate remains true under \eqref{eq:weakL2} instead of \eqref{eq:L2}.
\end{remark}

\subsection{Convexity of blow-ups}
\label{subsec:blowups-convex}

In this section, we show how to apply the Bernstein technique in order to prove \autoref{thm:blowups-convex}, which states that blow-ups for the nonlocal obstacle problem are necessarily convex. As explained before, this result is a central ingredient in our proof of \autoref{thm:OP_bdry-reg} and \autoref{thm:OP_opt-reg}.

\begin{proof}[Proof of \autoref{thm:blowups-convex}]
Let us take two cutoff functions $\eta \in C_c^{\infty}(B_2)$ and $\bar{\eta} \in C_c^{\infty}(B_{1/2})$ with $\eta, \bar{\eta} \ge 0$ and $\eta \equiv 1$ in $B_1$, $\bar{\eta} \equiv 1$ in $B_{1/4}$. Then, by \autoref{cor:keyest_diff-quot2_snd-order} applied to $\eta u_0$ with $-h$, we obtain for $x \in B_{1/2}$:
\begin{align*}
L(\bar{\eta}^2 & (D_{-h}D_{h}(\eta u_0))_+^2 + \sigma [(D_{h}(\eta u_0))^2]_{-h})(x)\\
 &\le 2 \bar{\eta}^2 L(D_{-h}D_{h}(\eta u_0)) (D_{-h}D_{h}(\eta u_0))_+(x) + 2\sigma [L(D_{h}(\eta u_0)) D_{h}(\eta u_0)]_{-h}(x).
\end{align*}
Since for $|h|$ small enough and $x \in B_{1/2}$
\begin{align*}
D_{-h}D_{h}(u_0)(x) = \frac{2u_0(x) - u_0(x+h) - u_0(x-h)}{|h|^2} = \frac{\frac{u_0(x) - u_0(x+h)}{|h|} + \frac{u_0(x) - u_0(x-h)}{|h|}}{|h|},
\end{align*}
using \eqref{eq:blowups-convex-ass2} for $x \in B_{1/2} \cap \{ u_0 > 0 \}$, and $D_{-h}D_{h} u_0 \le 0$ in $B_{1/2} \cap \{u_0 = 0\}$, we obtain
\begin{align*}
L(D_{-h}D_{h} u_0)(D_{-h}D_{h} u_0)_+(x) \le 0, ~~ x \in B_{1/2}.
\end{align*}
Moreover, note that by \eqref{eq:blowups-convex-ass1}
\begin{align*}
\Vert D_{h}[(1-\eta) u_0] \Vert_{L^{\infty}(B_R)} \le C(\Vert u_0\Vert_{L^{\infty}(B_3)} + \Vert \nabla u_0 \Vert_{L^{\infty}(B_{2R})}) \le C(\Vert u_0\Vert_{L^{\infty}(B_3)} + R^{s+\alpha}) ~~ \forall R \ge \frac{1}{8},
\end{align*}
and moreover that due to \eqref{eq:L1} it holds $D_{h}K(x-y) \le c|x-y|^{-n-1-2s} \le C|y|^{-n-1-2s}$ for $x \in B_{1/8}$ and $y \notin B_{1/4-2|h|}$. Hence, after integrating by parts and using the aforementioned observations, we get
\begin{align*}
    \left|L(D_{-h}D_{h} [(1-\eta) u_0])\right| &= \left|\int_{\R^n \setminus B_{1/4-2|h|}} -D_{-h}D_{h} [(1-\eta) u_0](y) K(x-y) \d y\right|\\
    &= \left|\int_{\R^n \setminus B_{1/4-2|h|}} D_{h} [(1-\eta) u_0](y) D_{h}K(x-y) \d y\right|\\
    &\le C\int_{\R^n \setminus B_{1/4-2|h|}} (\Vert u_0 \Vert_{L^\infty(B_3)}+|y|^{s+\alpha})|x - y|^{-n-1-2s} \d y\\
    &\le C ( \Vert u_0 \Vert_{L^{\infty}(B_3)} + 1).
\end{align*}
Therefore, we obtain
\begin{align}
\label{eq:convex-blowup-1}
\begin{split}
    L&(D_{-h}D_{h}(\eta u_0))(D_{-h}D_{h}(\eta u_0))_+(x)\\
    &= \left( L(D_{-h}D_{h} u_0) - L(D_{-h}D_{h} [(1-\eta) u_0])\right) (D_{-h}D_{h}(\eta u_0))_+(x)\\
    &\le C ( \Vert u_0 \Vert_{L^{\infty}(B_3)} + 1) (D_{-h}D_{h}(\eta u_0)(x))_+.
    \end{split}
\end{align}
Moreover, we obtain
\begin{align*}
L(D_{h}u_0) D_{h}u_0(x) \le 0 \quad \textrm{in}\quad  B_{3/4}
\end{align*}
by using the same argument as in Step 1 of the proof of  \autoref{thm:intro-OP-sc} together with \eqref{eq:blowups-convex-ass2}.
Indeed, if $u_0(x), u_0(x+h) > 0$, \eqref{eq:blowups-convex-ass2} applied with $x$ and with $x+h$ implies that $L(D_{h}u_0) = 0$, and if $u_0(x) = u_0(x+h) = 0$, then $D_{h}u_0(x) = 0$. Moreover, if $u_0(x) > u_0(x+h) = 0$, then $L(D_{h}u_0) \ge 0$ and $D_{h}u_0 < 0$. The case $u_0(x+h) > u_0(x) = 0$ follows by changing roles of $x$ and $x+h$.\\
Therefore, we deduce
\begin{align}
\label{eq:convex-blowup-2}
\begin{split}
    L(D_{h}(\eta u_0)) D_{h}(\eta u_0)(x) &= \left(L(D_{h}u_0) - L(D_{h}[(1-\eta)u_0]) \right) D_{h}u_0(x) \\
    &\le C( \Vert u_0 \Vert_{L^{\infty}(B_3)} + 1) |D_{h}u_0(x)|,
    \end{split}
\end{align}
where we used \eqref{eq:blowups-convex-ass1} again as above to estimate $- L(D_{h}[(1-\eta)u_0]) \le C( \Vert u_0 \Vert_{L^{\infty}(B_3)} + 1)$.
Altogether, by combining \eqref{eq:convex-blowup-1} and \eqref{eq:convex-blowup-2} we obtain for any $x \in B_{1/4}$
\begin{align*}
L&(\bar{\eta}^2 (D_{-h}D_{h}(\eta u_0))_+^2 + \sigma [(D_{h}(\eta u_0))^2]_{-h})(x)\\
&\le C \bar{\eta}^2 ( \Vert u_0 \Vert_{L^{\infty}(B_3)} + 1) (D_{-h}D_{h}(\eta u_0)(x))_+ + C( \Vert u_0 \Vert_{L^{\infty}(B_3)} + 1) |D_{h}u_0(x)|.
\end{align*}
Consequently, by the maximum principle (see Corollary 5.2 in \cite{Ros16}), the definitions of $\eta, \bar{\eta}$, and an application of Young's inequality (as in \eqref{eq:obstacle-proof-reference}), we obtain
\begin{align*}
\sup_{B_{1/4}} (D_{-h}D_{h}(\eta u_0))_+^2 &\le C \left(\Vert u_0 \Vert_{L^{\infty}(B_3)}^2 + \sup_{\R^n \setminus B_{1/2}} [(D_{h}(\eta u_0))^2]_{-h} \right) \\
&\le C \left( \Vert u_0 \Vert_{L^{\infty}(B_3)}^2 + \Vert \nabla u_0 \Vert_{L^{\infty}(B_{3})}^2 \right).
\end{align*}
In particular, upon sending $h \to 0$ and using \eqref{eq:blowups-convex-ass1}, we have shown
\begin{align*}
-D^2 u_0 \le C \left( \Vert u_0 \Vert_{L^{\infty}(B_3)} + \Vert \nabla u_0 \Vert_{L^{\infty}(B_3)}\right) < \infty ~~ \text{ in } B_{1/4}.
\end{align*}
Next, we observe that for any $r \ge 1$, the function
\begin{align*}
u^{(r)}_0(x) := \frac{u_0(rx)}{r^{1+s+\alpha}}
\end{align*}
satisfies all the assumptions of the theorem, and therefore, by application of the same arguments as before, we obtain for $x \in B_{1/4}$:
\begin{align*}
-D^2u^{(r)}_0(x) \le C \left( \Vert u^{(r)}_0 \Vert_{L^{\infty}(B_3)}  + \Vert \nabla u^{(r)}_0 \Vert_{L^{\infty}(B_3)} \right) \le C \frac{\Vert u_0 \Vert_{L^{\infty}(B_{3r})}}{r^{s+\alpha+1}} + \frac{\Vert \nabla u_0 \Vert_{L^{\infty}(B_{3r})}}{r^{s+\alpha}} \le C,
\end{align*}
where we used \eqref{eq:blowups-convex-ass1} and the fact that \eqref{eq:blowups-convex-ass1} also implies a growth control on $u_0$ itself:
\begin{align*}
    |u_0(x)| \le |u_0(0)| + r \Vert\nabla u_0\Vert_{L^{\infty}(B_{r})} \le c r^{s+\alpha+1} ~~ \forall x \in B_r, ~~ \forall r \ge 1.
\end{align*}
Here, $c,C > 0$ might depend on $u_0$, but not on $r$, once $r$ is large enough.
Consequently, for any $x \in \R^n$, by choosing $r = 8R|x|$ for $R > 1$, where $\widetilde{x} = \frac{x}{8R|x|} \in B_{1/4}$:
\begin{align*}
D^2u_0(x) = r^{s+\alpha-1} D^2u^{(r)}_0(\widetilde{x}) \le C r^{s+\alpha-1} \le C |x|^{s+\alpha-1} R^{s+\alpha-1} \to 0 ~~ \text{ as } R \to \infty,
\end{align*}
where we used that $s + \alpha < 1$. This implies the desired result.

Finally, notice that when $K$ is homogeneous, $0 \in \partial \{ u_0 > 0 \}$, and $u_0 \in C^1(\{ u_0 > 0 \})$ solves in the viscosity sense $L(\nabla u_0) = 0$ in $\{ u_0 > 0 \}$, then we can then apply the classification of blow-ups in \cite[Theorem 4.1]{CRS17} with $\Omega = \{ u_0  = 0 \}$ to obtain $u_0(x)=\kappa (x\cdot \nu)_+^{1+s}$ for some $\nu\in\mathbb S^{n-1}$, $\kappa\geq0$.
\end{proof}

\subsection{Optimal regularity for solutions}
\label{subsec:OP-main}

We next turn our attention to proving \autoref{thm:OP_bdry-reg} and \autoref{thm:OP_opt-reg}, i.e.,   optimal $C^{1+s}$ regularity for solutions to the obstacle problem \eqref{eq:OP} and regularity of the free boundary near regular points. 
Moreover, we also prove optimal regularity estimates in the presence of non-smooth obstacles $\phi \in C^{\beta}$ with $\beta <1+s$.
Throughout this subsection, we assume that $L \in \mathcal{L}_s(\lambda,\Lambda;1)$ and is homogeneous, i.e., \eqref{eq:hom} holds.

The proofs in this section follow the same overall strategy as the proof of optimal regularity for the global obstacle problem in \cite{FRS23} (see also \cite{CFR24}). However, as we explained before, we do not deduce the convexity of the blow-ups from semiconvexity estimates of solutions, but apply \autoref{thm:blowups-convex}, where we proved it directly. In the sequel we will sketch the proof of our main results, following the scheme in \cite{FRS23} and emphasizing the main differences.

As a preparation, we need the following slightly improved version of Theorem 5.1 from \cite{CDS17}, whose proof uses several ideas from Lemma 3.6 in \cite{AbRo20}. For the sake of readability, we postpone it to the Appendix (see Section \ref{sec:appendix}).

\begin{lemma}
\label{lemma:CDS}
Let $s\in(0,1)$, $L \in \mathcal{L}_s(\lambda,\Lambda;1)$. Let $f \in C^{\beta-2s}(B_1)$ for some $\beta \in (2s, 1+s)$, and let $\alpha \in (0,s)$. 
Let $u$ be a viscosity solution to the obstacle problem \eqref{eq:OP}.
Then, $u \in C^{\max\{2s+\epsilon,1+\epsilon\}}(B_{1/2})$ and it satisfies the following estimate
\begin{align*}
\Vert u \Vert_{C^{\max\{2s+\epsilon,1+\epsilon\}}(B_{1/2})} \le 
C \left([f]_{C^{\beta-2s}(B_1)} + \left\Vert \frac{u}{(1+|x|)^{1 + s+\alpha}} \right\Vert_{L^{\infty}(\R^n)} \right),
\end{align*}
where $C = C(n,s,\lambda,\Lambda) > 0$ is a uniform constant.
\end{lemma}

Using this, we can now give the:

\begin{proof}[Proof of \autoref{thm:OP_opt-reg} for $\beta<1+s$]
First of all, recall that the case $\beta<\max\{2s+\epsilon,1+\epsilon\}$ was done in \cite[Theorem 5.1]{CDS17}.
In particular, we may assume $\beta>\max\{2s,1\}$ and that $u \in C^{2s + \eps}(B_1)$ for some small $\eps > 0$, i.e., $u$ is a classical solution. 
Our goal is to prove that, when $\beta \in (\max\{2s,1\} , 1+s)$, we have
\begin{align}
\label{eq:rough-obstacle-help1}
|\nabla (u-\phi)(x)| \le C ( \Vert \phi \Vert_{C^{\beta}(B_1)} + \Vert u \Vert_{L^{\infty}(\R^n)}) |x-x_0|^{\beta-1}
\end{align}
for any free boundary point $x_0 \in \partial \{ u > \phi \}$.\\
Note that \eqref{eq:rough-obstacle-help1} implies the desired result upon combining it with interior regularity estimates in the same way as in Step 2 of the proof of \cite[Proposition 2.6.4]{FeRo24}.
Our proof of \eqref{eq:rough-obstacle-help1} is based on the ideas of the proof of \cite[Corollary 2.12]{FRS23}. 
After a normalization, we can assume without loss of generality that $x_0 = 0$ and that $\Vert \phi \Vert_{C^{\beta}(B_1)} + \Vert u \Vert_{L^{\infty}(\R^n)} = 1$.
Moreover, let us replace $u$ by $w := v \eta := (u-\phi)\eta$ for some cutoff-function $\eta \in C_c^{\infty}(B_1)$ with $0 \le \eta \le 1$ and $\eta \equiv 1$ in $B_{7/8}$. 
Clearly, $w$ satisfies
\begin{align*}
\min\{L w - f , w\} = 0~~ \text{ in}\quad B_{3/4},
\end{align*}
where $f = - L \phi - L((1-\eta)v) \in C^{\beta-2s}(B_{3/4})$ satisfies $\Vert f \Vert_{C^{\beta-2s}(B_{3/4})} \le C$ for a uniform constant $C = C(n,s,\lambda,\Lambda,\beta) > 0$. Therefore, by \cite[Theorem 5.1]{CDS17} (see also \autoref{lemma:CDS}), it holds that $\nabla w$ is globally bounded by a uniform constant. Note that the regularity of $f$ follows since by $\phi \in C^{\beta}(B_1)$, we have that $L \phi \in C^{\beta-2s}(B_{3/4})$, and moreover, $L((1-\eta)v) \in C^{0,1}(B_{3/4})$, as one can see from \eqref{eq:L1} and the identity
\begin{align}
\label{eq:fL-computation}
L(v(1-\eta))(x) = -\int_{\R^n} [v(1-\eta)](y) K(x-y)\d y.
\end{align}

Let us now begin with the proof of \eqref{eq:rough-obstacle-help1}. We need to show that
\begin{align}
\label{eq:rough-obstacle-help2}
|\nabla w(x)| \le C |x|^{\mu},
\end{align}
where $\mu=\beta-1$ and $C$ depends only on $n$, $s$, $\beta$, $\lambda$, $\Lambda$.
Let us assume by contradiction that \eqref{eq:rough-obstacle-help2} does not hold. Then, we can find sequences $(w_k)$, $(L_k)$, $(f_k)$ with the properties discussed before such that $0 \in \partial \{ w_k > 0 \}$ and
\begin{align}
\label{eq:rough-obstacle-help7}
\sup_{k} \sup_{r > 0} \frac{\Vert \nabla w_k \Vert_{L^{\infty}(B_r)}}{r^{\mu}} = \infty.
\end{align}
By \cite[Lemma 2.12]{FRS23}, there exist sequences $(r_m)$, $(k_m)$ with $r_m \to 0$ and $k_m \to 0$, as $m \to \infty$, such that
\begin{align*}
\widetilde{w}_m = \frac{w_{k_m}(r_m x)}{r_m \Vert \nabla w_{k_m} \Vert_{L^{\infty}(B_{r_m})}}, ~~ \nabla\widetilde{w}_m = \frac{\nabla w_{k_m}(r_m x)}{\Vert \nabla w_{k_m} \Vert_{L^{\infty}(B_{r_m})}}
\end{align*}
satisfy
\begin{align}
\nonumber
\Vert \nabla \widetilde{w}_m \Vert_{L^{\infty}(B_1)} &= 1,\\
\label{eq:rough-obstacle-uniform-est-growth}
|\nabla \widetilde{w}_m (x)| &\le C (1 + |x|^{\mu}) ~~ \forall x \in \R^n,\\
\nonumber
\min\{L_m \widetilde{w}_m - \widetilde{f}_m , \widetilde{w}_m\} &= 0~~ \text{ in } B_{3/(4r_m)},
\end{align}
where we write $L_m := L_{k_m}$ and define
\begin{align*}
\widetilde{f}_m := r_m^{2s-1}\frac{f_{k_m}(r_m x)}{\Vert \nabla w_{k_m}\Vert_{L^{\infty}(B_{r_m})}}.
\end{align*}
Consequently, by the uniform control on $\nabla \widetilde{w}_m$, we have $\widetilde{w}_m \to \widetilde{w}_0$ locally uniformly for some limiting function $\widetilde{w}_0$ up to a subsequence. Moreover, observe that
\begin{align}
\label{eq:rough-obstacle-help4}
L_m(D_{h} \widetilde{w}_m) &\ge  D_{h} \widetilde{f}_m ~~ \text{ in } \{ \widetilde{w}_m > 0 \} \cap B_{3/(4r_m)},\\
\label{eq:rough-obstacle-help5}
L_m(D_{h} \widetilde{w}_m) &=  D_{h} \widetilde{f}_m ~~ \text{ in } \{ x : \dist(x, \{\widetilde{w}_m = 0\}) > |h| \} \cap B_{3/(4r_m)},
\end{align}
and that by the proof of Lemma 2.12 in \cite{FRS23} we also have
\begin{align*}
r_m^{-\mu} \Vert \nabla w_{k_m} \Vert_{L^{\infty}(B_{r_m})} \ge \frac{1}{2} \sup_{k} \sup_{r \ge r_m} r^{-\mu} \Vert \nabla w_k \Vert_{L^{\infty}(B_r)}.
\end{align*}
By \eqref{eq:rough-obstacle-help7}, we can extract a further subsequence $(l_m) \subset (k_m)$ such that for every $m \in \N$:
\begin{align}
\label{eq:rough-obstacle-help8}
\Vert \nabla w_{l_m} \Vert_{L^{\infty}(B_{r_{l_m}})} \ge m r_{l_m}^{\mu}.
\end{align} 
To simplify the presentation of the proof, let us slightly abuse notation and write again $r_m$ and $k_m$ instead of $r_{l_m}$ and $l_{m}$. Now, since $f_{k_m} \in C^{\beta-2s}(B_{1})$ we have for any $|h| > 0$ and $x \in B_{3/(4r_m)}$:
\begin{align}
\label{eq:rough-obstacle-help6}
| D_{h} \widetilde{f}_m (x) | \le m^{-1} r_m^{2s-1-\mu} |D_{h} f_{k_m}(r_m x)| \le c m^{-1} r_m^{\beta-1-\mu} |h|^{\beta-2s-1} \to 0, ~~ \text{ as } m \to \infty,
\end{align}
where we used \eqref{eq:rough-obstacle-help8} and $\mu = \beta - 1$.
Moreover, we have $\widetilde{w}_m \in C^{2s+\eps}(B_{3/(4r_m)})$ and
\begin{align}
\label{eq:rough-obstacle-uniform-est}
\Vert \widetilde{w}_m \Vert_{C^{2s+\eps}(B_{3/(4r_m)})} \le C
\end{align}
for a uniform constant $C = C(n,s,\lambda,\Lambda) > 0$. This is a consequence of \autoref{lemma:CDS}, the gradient control on $\widetilde{w}_m$ and the fact that $[\widetilde{f}_m]_{C^{\beta-2s}(B_{3/(4r_m)})} \to 0$, as $m \to \infty$, which follows from the computation
\begin{align*}
| D^{\beta - 2s}_{h} \widetilde{f}_m (x) | \le m^{-1} r_m^{2s-1-\mu} |D^{\beta - 2s}_{h} f_{k_m}(r_m x)| \le c m^{-1} r_m^{\beta-1-\mu} ~~ \forall |h| > 0.
\end{align*}
By \eqref{eq:rough-obstacle-uniform-est}, we have uniform convergence $\widetilde{w}_m \to \widetilde{w}_0 \in C^{2s+\eps}_{loc}(\R^n)$ and using also \eqref{eq:rough-obstacle-uniform-est-growth}, we get convergence in $L^1_{\omega_s}(\R^n)$. For a definition of the space $L^1_{\omega_s}(\R^n)$, we refer to \cite[Definition 1.10.1]{FeRo24}.
Therefore, using \eqref{eq:rough-obstacle-help4}, \eqref{eq:rough-obstacle-help6}, and the stability of viscosity solutions (see \cite[Proposition 3.1.12]{FeRo24}), we obtain
\begin{align*}
L_{\infty} (D_{h} \widetilde{w}_0) \ge 0 ~~ \text{ in } \{ \widetilde{w}_0 > 0 \},
\end{align*}
where $L_{\infty}$ denotes the weak limit of $(L_m)_m$. Here, we apply the stability of viscosity solutions in every ball contained in $\{ \widetilde{w}_0 > 0\}$. Note that every such ball is also contained in $\{ \widetilde{w}_m > 0 \}$ for large enough $m$ by \eqref{eq:rough-obstacle-uniform-est}. Thus, $
\widetilde{w}_0$ satisfies \eqref{eq:blowups-convex-ass2}. 

Moreover, using again the stability of viscosity solutions and  \eqref{eq:rough-obstacle-help5} and \eqref{eq:rough-obstacle-help6}, we can conclude
\begin{align*}
L_{\infty}(\nabla \widetilde{w}_0) = 0 ~~ \text{ in } \{ \widetilde{w}_0 > 0 \}.
\end{align*}
and also
\begin{align}
\label{eq:rough-obstacle-help3}
\Vert \nabla \widetilde{w}_0 \Vert_{L^{\infty}(B_R)} \le c R^{\mu}, ~~ \forall R \ge 1.
\end{align}
Note that \eqref{eq:rough-obstacle-help3} implies \eqref{eq:blowups-convex-ass1} for any $\alpha \in (0,s)$.
Hence, we are in a position to apply \autoref{thm:blowups-convex} to $\widetilde{w}_0$. This yields convexity of $\widetilde{w}_0$, and
\begin{align*}
\widetilde{w}_0(x) = \kappa (x \cdot e)^{1+s}_+
\end{align*}
for some $\kappa \ge 0$ and $e \in \mathbb S^{n-1}$. 
By convexity of $\widetilde{w}_0$ and $\Vert \nabla \widetilde{w}_m \Vert_{L^{\infty}(B_1)} = 1$, it holds that $\Vert \nabla \widetilde{w}_0 \Vert_{L^{\infty}(B_1)} \ge 1$. 
Therefore $\kappa \neq 0$. However, since $\mu = \beta - 1 < s$, this is a contradiction to \eqref{eq:rough-obstacle-help3}. Therefore, we have proved \eqref{eq:rough-obstacle-help2} and \eqref{eq:rough-obstacle-help1}. 
\end{proof}

Our proof of the regularity of the free boundary (see \autoref{thm:OP_bdry-reg}) and the optimal regularity (see \autoref{thm:OP_opt-reg} in case $\beta > 1+2s$) follows the strategy developed in \cite{FRS23}. The main ingredient is the following quantitative estimate, which is an analog of \cite[Theorem 2.2]{FRS23}. However, due to the classification of blow-ups in \autoref{thm:blowups-convex}, we do not need to assume semiconvexity of $u$ in order to establish \eqref{eq:blowup-close}. This is the only difference compared to \cite[Theorem 2.2]{FRS23}.

\begin{lemma}
\label{lemma:Thm2.2}
Let $s\in(0,1)$, $L \in \mathcal{L}_s(\lambda,\Lambda;1)$ satisfying \eqref{eq:hom}. Let $\alpha \in (0,s) \cap (0,1-s)$. Let $\eta > 0$ and assume that $u \in C^{0,1}(\R^n)$ satisfies
\begin{align*}
u &\ge 0 \text{ in } B_{1/\eta}, ~~\text{ with } 0 \in \partial \{u > 0\},\\
\min\{Lu - f, u\} &= 0 ~~ \text{ in } B_{1/\eta}, ~~\text{ with } |\nabla f| \le \eta,\\
\Vert \nabla u \Vert_{L^{\infty}(B_R)} &\le R^{s+\alpha}, ~~ \forall R \ge 1.
\end{align*}
Then, there exist $e \in \mathbb S^{n-1}$ and $\kappa \ge 0$ such that
\begin{align}
\label{eq:blowup-close}
\Vert u - \kappa (x \cdot e)_+^{1+s} \Vert_{C^{0,1}(B_1)} \le \eps(\eta),
\end{align}
where $\eps(\eta)$ is a modulus of continuity, depending only on $n,s,\lambda,\Lambda,\alpha$.\\
Moreover, for any $\kappa_0 > 0$ there is $\eps_0 > 0$ such that if $\kappa \ge \kappa_0$ and $\eps(\eta) \le \eps_0$, then $\partial \{ u > 0\}$ is a $C^{1,\gamma}$-graph in $B_{1/2}$ for some $\gamma > 0$ depending only on $n$, $s$, $\lambda$, $\Lambda$, $\alpha$, $\kappa_0$, and
\begin{align*}
|u - \kappa (x \cdot e)_+^{1+s}| \le C |x|^{1+s+\gamma}, ~~ |\nabla u - \nabla \kappa (x \cdot e)_+^{1+s}| \le C |x|^{s+\gamma} ~~ \text{ in } B_1,
\end{align*}
where $C > 0$ depends only on $n$, $s$, $\lambda$, $\Lambda$, $\alpha$, $\kappa_0$.
\end{lemma}

\begin{proof}
First, we prove \eqref{eq:blowup-close}. Let us assume by contradiction that there exists $\eps > 0$ and a sequence $(\eta_k)$ with $\eta_k \to 0$ and a sequence of operators $(L_k)$ and solutions $(u_k)$ satisfying the assumptions of \autoref{lemma:Thm2.2} for which
\begin{align}
\label{eq:contradiction}
\Vert u_k - \kappa (x \cdot e)_+^{1+s} \Vert_{C^{0,1}(B_1)} > \eps
\end{align}
for any $e \in \mathbb S^{n-1}$ and $\kappa \ge 0$. Similar to  \eqref{eq:rough-obstacle-uniform-est} in the proof of \autoref{thm:OP_opt-reg} in case $\beta < 1+s$, we can show that $\Vert u_k \Vert_{C^{1+s-\eps}(B_{1})}$ is uniformly bounded and that $u_k$ converges, up to a subsequence, in the $C^1$-norm to a limiting function $u_0$ satisfying the assumptions of \autoref{thm:blowups-convex}. Consequently, there exist $e \in \mathbb{S}^{n-1}$ and $\kappa \ge 0$ such that $u_0(x) = \kappa (x \cdot e)_+^{1+s}$. This contradicts \eqref{eq:contradiction}. Thus, we have shown \eqref{eq:blowup-close}.
With \eqref{eq:blowup-close} at hand we are in the exact same setup as in   \cite[Theorem 2.2 (ii), (iii)]{FRS23}. Note that the proof of Theorem 2.2 (ii), (iii) does not use the semiconvexity assumption, any more. Hence, we can apply \cite[Theorem 2.2 (ii), (iii)]{FRS23} and the proof is complete.
\end{proof}

We are now ready to establish the dichotomy between regular and degenerate free boundary points, show $C^{1,\gamma}$-regularity of the free boundary near regular points (see \autoref{thm:OP_bdry-reg}), and establish the optimal $C^{1+s}$-regularity of solutions (see \autoref{thm:OP_opt-reg} in case $\beta > 1 + 2s$).

\begin{proof}[Proof of \autoref{thm:OP_bdry-reg} and of \autoref{thm:OP_opt-reg} in case $\beta > 1 + 2s$]

Dividing the solution $u$ by a constant, and up to a translation, we can assume $x_0 = 0$ and $\Vert \phi \Vert_{C^{\beta}(B_1)} + \Vert u \Vert_{L^{\infty}(\R^n)} = 1$. Moreover, replacing $u$ by $w = (u-\phi)\eta$ for some cutoff-function $\eta \in C_c^{\infty}(B_1)$ with $0 \le \eta \le 1$ and $\eta \equiv 1$ in $B_{7/8}$ as in the proof of \autoref{thm:OP_opt-reg} in case $\beta < 1 + s$,  we have that $w$ satisfies the following properties:
\begin{align*}
w \ge 0 ~~ \text{ in } B_1, \quad \Vert \nabla w \Vert_{L^{\infty}(\R^n)} \le C, \quad w(0) = |\nabla w(0)| = 0,\\
\min\{L w - f , w \} = 0 ~~ \text{ in } B_{1/2}, ~~ \text{ with } |\nabla f| \le C.
\end{align*}
Now, we consider the rescalings for $\eta \in (0,1)$ and $k \in \N \cup \{ 0 \}$:
\begin{align*}
w_k(x) := \frac{\eta}{C}\frac{w(\eta^{-1} 2^{-k+1}x)}{( \eta^{-1} 2^{-k+1})^{1+s+\alpha}}, ~~ k \in \N,
\end{align*}
where $\alpha \in (0,s) \cap (0,1-s)$ (as in \autoref{lemma:Thm2.2}), to guarantee that $w_k$ solves 
\begin{align*}
\min\{ L w_k - f_k , w_k \} = 0 ~~ \text{ in } B_{1/\eta}, ~~\text{ with } |\nabla f_k| \le \eta.
\end{align*}
Moreover, for $k = 0$ it holds $\Vert \nabla w_0 \Vert_{L^{\infty}(\R^n)} \le 1$. This means that we have verified all the assumptions of \autoref{lemma:Thm2.2} for $w_k$, except for the growth condition, but this one holds true at least for $k = 0$. Hence, we are in the exact same position as in the proof of \cite[Corollary 2.17]{FRS23}. We will follow their proof, applying \autoref{lemma:Thm2.2} instead of \cite[Theorem 2.2]{FRS23}.

There are two cases: First, let us assume that $w_k$ satisfy
\begin{align}
\label{eq:growth-ass}
\Vert w_k \Vert_{L^{\infty}(B_R)} \le R^{s+\alpha} ~~ \forall R \ge 1, ~~ \forall k \ge 0.
\end{align}
Then, 
\begin{align*}
\Vert \nabla u \Vert_{L^{\infty}(B_{\eta^{-1} 2^{-k+1}})} = C \eta^{-1} (\eta^{-1} 2^{-k+1})^{s+\alpha} \Vert \nabla w_k \Vert_{L^{\infty}(B_1)} \le c (2^{-k})^{s+\alpha},
\end{align*}
which implies
\begin{align*}
|\nabla u(x)| \le c |x|^{s+\alpha}, \qquad |u(x)| \le c |x|^{1+s+\alpha} ~~ \forall x \in B_1,
\end{align*}
and the proof is complete in this case.\\
On the other hand, if \eqref{eq:growth-ass} does not hold, then there is a maximal number $k_0 \in \N$ such that
\begin{align}
\label{eq:growth-ass-2}
\Vert \nabla w_k \Vert_{L^{\infty}(B_R)} \le R^{s+\alpha} ~~ \forall R \ge 1, ~~ k \le k_0.
\end{align}
Then, choosing $\eta$ small enough, we can apply \autoref{lemma:Thm2.2} and obtain that
\begin{align*}
\Vert w_{k_0} - \kappa (x \cdot e)_+^{1+s} \Vert_{C^{0,1}(B_1)} \le \eps := \min \{ \eps_0 , 1/6 \}.
\end{align*}
for some $\kappa \ge 0$. Note that in case $\kappa \le \frac{1}{6}$ we have
\begin{align*}
\Vert \nabla w_{k_0} \Vert_{L^{\infty}(B_1)} \le \Vert \nabla \kappa (x \cdot e)_+^{1+s} \Vert_{L^{\infty}(B_1)} + \eps \le \frac{1}{2} \le 2^{-s-\alpha}.
\end{align*}
Since $\nabla w_{k_0 + 1}(x) = 2^{s+\alpha} \nabla w_{k_0}(\frac{x}{2})$ this implies
\begin{align*}
\Vert \nabla w_{k_0 + 1} \Vert_{L^{\infty}(B_2)} \le 1, \qquad \Vert \nabla w_{k_0 + 1} \Vert_{L^{\infty}(B_R)} = 2^{s+\alpha} \Vert \nabla w_{k_0} \Vert_{L^{\infty}(B_{R/2})} \le R^{s+\alpha} ~~ \forall R \ge 2.
\end{align*}
Hence, $w_{k_0+1}$ still satisfies the growth condition, contradicting the maximality of $k_0$. Thus, we can assume $\kappa \ge \frac{1}{6}$. Then, by \autoref{lemma:Thm2.2}, we have that the free boundary $\partial \{ w_{k_0} > 0 \}$ is a $C^{1,\gamma}$-graph in $B_{1/2}$ for some $\gamma > 0$, and
\begin{align*}
|\nabla w_{k_0} - \nabla \kappa (x \cdot e)_+^{1+s}| \le c |x|^{s+\gamma} ~~ \text{ in } B_1.
\end{align*}
By the boundary Harnack inequality (see \cite{CRS17}) we can replace $\gamma$ by $\alpha$ in the aforementioned estimate, and hence
\begin{align*}
|w_{k_0}(x) - \kappa (1+s)^{-1} (x \cdot e)_+^{1+s}| \le C_1 |x|^{1+s+\alpha} ~~ \text{ in } B_1. 
\end{align*}
This implies for $u$
\begin{align*}
|u(x) - c_0 (x \cdot e)_+^{1+s}| \le c |x|^{1+s+\alpha} ~~ \forall x \in B_{\eta^{-1} 2^{-k_0+1}}
\end{align*}
for some constant $c_0 = c_1 (\eta^{-1}2^{-k_0+1})^{\alpha} > 0$ with $c_1$ depending only $n,\eta,s,\kappa,C_1$, and that the free boundary $\partial \{ u > 0 \}$ is $C^{1,\alpha}$ in a ball of radius $\eta^{-1} 2^{-k_0+1}$. Hence, using \eqref{eq:growth-ass-2} we obtain for $x \in B_1 \setminus B_{\eta^{-1} 2^{-k_0+1}}$:
\begin{align*}
|u(x) - c_0 (x \cdot e)_+^{1+s}| \le |u(x)| + |c_0 (x \cdot e)_+^{1+s}| \le c |x|^{1+s+\alpha} + c_1 (\eta^{-1} 2^{-k_0+1})^{\alpha} |x|^{1+s} \le c |x|^{1+s+\alpha},
\end{align*}
which yields the desired expansion also in this case. Moreover, note that in this case we also have an estimate $|\nabla u(x)| \le C |x|^{s}$ for $x \in B_1$, and therefore we also have the $C^{1+s}$-regularity of $u$. The proof is complete.
\end{proof}

\section{Extensions of the Bernstein technique}
\label{sec:extensions}

In this section, we present several possible extensions of the Bernstein technique. To be more precise, in Subsection \ref{subsec:parabolic}, we explain how to establish Bernstein key estimates for operators of the form $\partial_t + L$ and explain how to derive a priori estimates for solutions to the corresponding parabolic equation $\partial_t u + L u = f$. Moreover, in Subsection \ref{subsec:nonsymm}, we extend the key estimates \autoref{thm:keyest} and \autoref{thm:keyest_pos-part} to operators with nonsymmetric kernels that possess first order drift terms.

\begin{remark}[Bernstein technique for convolution operators]
We mention that the Bernstein technique also works for nonlocal operators with integrable kernels of the form
\begin{align*}
Lu(x) = \int_{\R^n} (u(x) - u(y))K(x-y) \d y = u(x) - \int_{\R^n} u(y) K(x-y) \d y.
\end{align*}
In fact, we have the following:
\begin{itemize}
\item[(i)] Assume that $K \in C^{0,1}(\R^n)$ satisfies  for some $\Lambda > 0$
\begin{align*}
|\nabla K(y)| \le \Lambda K(y), \qquad \int_{\R^n} K(y) \d y = 1.
\end{align*}
Let $L$ be as before. Let $\eta \in C_c^{\infty}(B_1)$ be such that $\eta \equiv 1$ in $B_{1/2}$ and $0 \le \eta \le 1$.  Then, there exists $\sigma_0 = \sigma_0(n, \Lambda) > 0$ such that for every $\sigma \ge \sigma_0$ and $u \in L^{\infty}(\R^n)$  the key estimates for difference quotients \eqref{eq:keyest_diff-quot} and \eqref{eq:keyest_diff-quot_pos-part} hold true. Moreover, if $u \in C^{1+\eps}(\R^n)$ the key estimates \eqref{eq:keyest} and \eqref{eq:keyest_pos-part} hold true.
\item[(ii)] The proof of (i) follows along the lines of Step 1 in the proof of the key estimate \autoref{thm:keyest} (respectively its modifications), replacing $K_2$ by $K$.
\item[(iii)] Using the key estimate for difference quotients, one can prove that solutions $u$ to $L u = f$ in $\Omega$ with $f \in C^{0,1}(\Omega)$ satisfy $u \in C^{0,1}_{loc}(\Omega)$. In this way, one can argue that the Bernstein technique also extends to the limiting case $s = 0$.
\end{itemize}
\end{remark}

\subsection{Parabolic nonlocal equations}
\label{subsec:parabolic}

So far, in this work we have restricted ourselves to the study of elliptic problems. In this section, we explain how the Bernstein technique can be used to study regularity estimates for parabolic equations governed by nonlocal operators. As an application, we establish semiconvexity estimates for solutions to the parabolic obstacle problem.

\begin{theorem}
\label{thm:keyest_parabolic}
Let $s\in(0,1)$, $L_{K^{(t)}} \in \mathcal{L}_s(\lambda,\Lambda;1)$ for any $t \in [0,\infty)$. Let $\eta_1^{2s} \in C^{0,1}([0,1])$ be such that $\eta_1 \ge 0$. Let $\eta_2 \in C_c^{1,1}(B_1)$ be such that $\eta_2 \ge 0$. We define $\eta := \eta_1 \eta_2$. Then, there exists $\sigma_0 = \sigma_0(n, s, \lambda, \Lambda, \Vert \eta_1^{2s} \Vert_{C^{0,1}([0,1])}, \Vert \eta_2 \Vert_{C^{1,1}(B_1)}) > 0$ such that for every $\sigma \ge \sigma_0$ and $u,v \in C^{\infty}((0,\infty) \times \R^n)$
\begin{align}
\label{eq:keyest_parabolic}
(\partial_t + L_{K^{(t)}})(\eta^2 (\partial_e u)^2 + \sigma u^2) &\le 2 \eta^2 (\partial_t + L_{K^{(t)}})(\partial_e u) \partial_e u + 2\sigma (\partial_t + L_{K^{(t)}})(u) u,\\
\label{eq:keyest_parabolic2}
(\partial_t + L_{K^{(t)}})(\eta^2 (\partial_e v)_+^2 + \sigma u^2) &\le 2 \eta^2 (\partial_t + L_{K^{(t)}})(\partial_e v) (\partial_e v)_+ + 2\sigma (\partial_t + L_{K^{(t)}})(v) v.
\end{align}
\end{theorem}

\begin{remark}

By adapting the proofs in Section \ref{sec:diffquot} to the parabolic setting, it is also possible to establish key estimates in terms of difference quotients for $(\partial_t + L_{K^{(t)}})$:
\begin{align}
\label{eq:keyest_parabolic_diffquot1}
(\partial_t + L_{K^{(t)}}) (\eta^2 (D_{h} u)^2 + \sigma u_{h}^2) &\le 2 \eta^2 (\partial_t + L_{K^{(t)}}) (D_{h}u) D_{h}u + 2\sigma (\partial_t + L_{K^{(t)}}) (u_{h}) u_{h},\\
\label{eq:keyest_parabolic_diffquot2}
(\partial_t + L_{K^{(t)}})(\eta^2 (D_{h} u)^2 + \sigma [u^2]_{h}) &\le 2 \eta^2 (\partial_t + L_{K^{(t)}})(D_{h}u) D_{h}u + 2\sigma [(\partial_t + L_{K^{(t)}})(u) u]_{h}.
\end{align}
\end{remark}

\begin{proof}[Proof of \autoref{thm:keyest_parabolic}]
(i) First, we prove \eqref{eq:keyest_parabolic}. Let us compute
    \begin{align*}
        (\partial_t + L_{K^{(t)}}) &(\eta^2 (\partial_e u)^2 + \sigma u^2) = 2(\partial_t \eta_1)\eta_1 \eta_2^2 (\partial_e u)^2\\
        &+ 2 \eta^2 (\partial_t + L_{K^{(t)}}) (\partial_e u)(\partial_e u) + 2 \sigma (\partial_t + L_{K^{(t)}})(u) u\\
        &+ L_{K^{(t)}}(\eta^2)(\partial_e u)^2 - B_{K^{(t)}}(\eta^2 , (\partial_e u)^2) - \eta^2 B_{K^{(t)}}(\partial_e u, \partial_e u) - \sigma B_{K^{(t)}}(u,u).
    \end{align*}
    Therefore, it remains to prove
    \begin{align*}
        2(\partial_t \eta_1)\eta_1 \eta_2^2 (\partial_e u)^2 + L_{K^{(t)}}(\eta^2)(\partial_e u)^2 - B_{K^{(t)}}(\eta^2 , (\partial_e u)^2) \le \eta^2 B_{K^{(t)}}(\partial_e u, \partial_e u) + \sigma B_{K^{(t)}}(u,u).
    \end{align*}
    As in the proof of \autoref{thm:keyest}, let us choose $\eps = \gamma \eta(x)$, where $\gamma > 0$ is as before. Therefore, by Step 1 of the proof of \autoref{thm:keyest}, it remains to prove
    \begin{align*}
               2(\partial_t \eta_1)\eta_1 \eta_2^2 (\partial_e u)^2 + L_{K_1^{(t)}}(\eta^2)(\partial_e u)^2 - B_{K_1^{(t)}}(\eta^2 , (\partial_e u)^2) \le \eta^2 B_{K_1^{(t)}}(\partial_e u, \partial_e u) + \sigma B_{K_1^{(t)}}(u,u),
    \end{align*}
    where $K_1^{(t)}$ is defined as in \autoref{lemma:kernel-decomp}. Moreover, by carefully tracking Step 2 of the proof of \autoref{thm:keyest}, it becomes apparent that 
    \begin{align*}
         L_{K_1^{(t)}}(\eta^2)(\partial_e u)^2 - B_{K_1^{(t)}}(\eta^2 , (\partial_e u)^2) &= \eta_1^2 \left(L_{K_1^{(t)}}(\eta_2^2)(\partial_e u)^2 - B_{K_1^{(t)}}(\eta_2^2 , (\partial_e u)^2)\right)\\
         &\le \frac{3}{4}\eta^2 B_{K_1^{(t)}}(\partial_e u, \partial_e u) + \sigma_1 B_{K_1^{(t)}}(u,u)
    \end{align*}
    for some $\sigma_1 = \sigma_1(n,s,\lambda,\Lambda,\Vert \eta_2 \Vert_{C^{1,1}}(\R^n))$. Therefore, it remains to prove
    \begin{align}
    \label{eq:parabolic-help1}
        2(\partial_t \eta_1)\eta_1 \eta_2^2 (\partial_e u)^2 \le \frac{1}{4} \eta^2 B_{K_1^{(t)}}(\partial_e u, \partial_e u) + \sigma_2 B_{K_1^{(t)}}(u,u)
    \end{align}
    for some $\sigma_2 = \sigma_2(n,s,\lambda,\Lambda,\Vert \eta_1 \Vert_{C^{0,1}([0,1])} , \Vert \eta_2 \Vert_{C^{1,1}(\R^n)})$.
Note that \eqref{eq:parabolic-help1} is trivial once $\partial_t \eta_1(t) \le 0$. Thus, we can assume without loss of generality that $\partial_t \eta_1(t) > 0$ and choose $\delta =  \left(\frac{ \eta_1(t) }{ 8\partial_t \eta_1(t)}\right)^{\frac{1}{2s}} \wedge \frac{\eps}{2}$ and apply the interpolation estimate \autoref{lemma:interpol}. Then,
    \begin{align*}
        2(\partial_t \eta_1)\eta_1 \eta_2^2 (\partial_e u)^2 &\le 2 (\partial_t \eta_1) \eta_1 \eta_2^2 \left(\delta^{2s} B_{K_1^{(t)}}(\partial_e u,\partial_e u) + c \delta^{2s-2} B_{K_1^{(t)}}(u,u)\right) \\
        &\le \frac{1}{4}\eta^2 B_{K_1^{(t)}}(\partial_e u,\partial_e u) + c (\partial_t \eta_1)^{\frac{1}{s}} \eta_1^{2 - \frac{1}{s}} \eta_2^2 B_{K_1^{(t)}}(u,u)\\
        &\le \frac{1}{4}\eta^2 B_{K_1^{(t)}}(\partial_e u,\partial_e u) + c \eta_2^2 B_{K_1^{(t)}}(u,u).
    \end{align*}
Note that in the last step we used that $(\partial_t \eta_1)^{\frac{1}{s}} \eta_1^{2 - \frac{1}{s}} = c(s) (\partial_t(\eta_1^{2s}))^{\frac{1}{s}}$ and $\eta_2^2$ are bounded. This establishes \eqref{eq:parabolic-help1}, as desired.\\

(ii): Now, let us show \eqref{eq:keyest_parabolic2}. The proof follows by a combination of the arguments in the proof of \eqref{eq:keyest_parabolic} and \autoref{thm:keyest_pos-part}. In fact, it remains to prove
    \begin{align}
    \label{eq:parabolic-help2}
        2(\partial_t \eta_1)\eta_1 \eta_2^2 (\partial_e v)^2 \le \frac{1}{4} \eta^2 \left[ B_{K_1^{(t)}}((\partial_e v)_+, (\partial_e v)_+) - 2L_{K_1^{(t)}}((\partial_e v)_-)(\partial_e v)_+ \right] + \sigma_2 B_{K_1^{(t)}}(v,v)
    \end{align}
    for some $\sigma_2 = \sigma_2(n,s,\lambda,\Lambda,\Vert \eta_1 \Vert_{C^{0,1}([0,1])},  \Vert \eta_2 \Vert_{C^{1,1}(\R^n)})$.
    This will be achieved by application of \autoref{lemma:interpol_pos-part} with $\delta = \frac{1}{8} \left(\frac{ \eta_1(t) }{ \partial_t \eta_1(t)}\right)^{\frac{1}{2s}} \wedge \eps$, where $\eps > 0$ is as in the proof of \autoref{thm:keyest_pos-part}. The rest of the proof goes as in (i).
\end{proof}

We also have one-sided key estimates of second order:

\begin{corollary}
\label{cor:keyest_snd-order_parabolic}
Let $s\in(0,1)$, $L_{K^{(t)}} \in \mathcal{L}_s(\lambda,\Lambda;1)$ for any $t \in [0,\infty)$. Let $\bar{\eta}_1^{2s} \in C^{0,1}$ be such that $\bar{\eta}_1(0) = 0$, $\bar{\eta}_1 \equiv 1$ in $[1/2,1]$, and $0 \le \bar{\eta}_1 \le 1$. Let $\bar{\eta}_2 \in C_c^{1,1}(B_{1/4})$ be such that $\bar{\eta}_2 \equiv 1$ in $B_{1/8}$ and $0 \le \bar{\eta}_2 \le 1$. We define $\bar{\eta} := \bar{\eta}_1 \bar{\eta}_2$.  Then, there exists $\sigma_0 = \sigma_0(n, s, \lambda, \Lambda, \Vert \bar{\eta}_1^{2s} \Vert_{C^{0,1}([0,1])}, \Vert \bar{\eta}_2 \Vert_{C^{1,1}(B_1)}) > 0$ such that for every $\sigma \ge \sigma_0$ and $u \in C^{\infty}((0,1) \times \R^n)$ the following estimate holds true:
\begin{align}
\label{eq:keyest_snd-order_parabolic_one-sided}
\begin{split}
   (\partial_t + L_{K^{(t)}})&(\bar{\eta}^2 (\partial^2_{ee} u)_+^2 + \sigma (\partial_e u)^2) \\
   &\le 2\bar{\eta}^2 (\partial_t + L_{K^{(t)}})(\partial^2_{ee} u)(\partial^2_{ee} u)_+ + 2 \sigma (\partial_t + L_{K^{(t)}})(\partial_e u) \partial_e u.
\end{split}
\end{align}
Moreover, the following estimate holds true whenever all expressions are well-defined:
\begin{align}
\label{eq:keyest_snd-order_parabolic_one-sided-diffquot}
\begin{split}
&(\partial_t + L_{K^{(t)}})(\bar{\eta}^2 (D_{-h} D_{h} u)_+^2 + \sigma [(D_{h} u)^2]_{-h})\\
&~\le 2\bar{\eta}^2 (\partial_t + L_{K^{(t)}})(D_{-h} D_{h} u)(D_{-h} D_{h} u)_+ + 2 \sigma [ (\partial_t + L_{K^{(t)}})(D_{h} u) D_{h} u]_{-h}.
\end{split}
\end{align}
\end{corollary}

Let us now explain how to apply the parabolic Bernstein key estimate \autoref{thm:keyest_parabolic} in order to obtain semiconvexity estimates for solutions to the parabolic nonlocal obstacle problem.

By combination of the parabolic key estimate \autoref{thm:keyest_parabolic} and the parabolic maximum principle, we obtain semiconvexity estimates in space for smooth solutions to the parabolic nonlocal obstacle problem:

\begin{theorem}[Semiconvexity estimate]
Let $s\in(0,1)$, $L_{K^{(t)}} \in \mathcal{L}_s(\lambda,\Lambda;2)$ for any $t \in (0,1)$. Let $f \in L^{\infty}((0,1); C^{1,1}(B_1))$. Let $u\in C^{1+s}_x((0,1) \times B_1) \cap C^{0,1}_t((0,1) \times B_1)$ be a solution to the parabolic obstacle problem
\begin{align}
\label{eq:parOP}
\min\{\partial_t u + L_{K^{(t)}} u - f , u\} &= 0 ~~ \text{ in } (0,1) \times B_1.
\end{align}
Then, for any $e \in \mathbb S^{n-1}$, it holds
\begin{align}
\label{eq:parabolic_semiconvexity}
\inf_{(1/2,1) \times B_{1/2}} \partial_{ee}^2 u \ge -C \left( \Vert f \Vert_{L^{\infty}((0,1); C^{1,1}(B_1))} + \Vert u \Vert_{L^{\infty}((0,1) \times \R^n)}\right),
\end{align}
where $C = C(n,s,\lambda,\Lambda) > 0$ is a constant.
\end{theorem}

\begin{proof}
We split the proof into two parts.

\textbf{Step 1:} First order estimate: \\
Let $\eta$ be as in \autoref{thm:keyest_parabolic}. By the same arguments as in the proof of \eqref{eq:application_OP_first-order1}, we obtain in $(1/2,1) \times B_{1/2-|h|}$:
\begin{align*}
(\partial_t + L_{K^{(t)}})(D_{h} u) D_{h}u \le \Vert \nabla f \Vert_{L^{\infty}((0,1) \times B_1)}  |D_{h}u|.
\end{align*}
Moreover, since $u$ solves \eqref{eq:parOP}, we have
\begin{align*}
[(\partial_t + L_{K^{(t)}})(u) u]_{h} \le \Vert f \Vert_{L^{\infty}((0,1)\times B_1)} \Vert u \Vert_{L^{\infty}((0,1) \times B_1)}.
\end{align*}
Using \eqref{eq:keyest_parabolic_diffquot2}, we obtain
\begin{align*}
(\partial_t + L_{K^{(t)}})(\eta^2 (D_{h} u)^2 + \sigma [u^2]_{h}) &\le 2 \eta^2 (\partial_t + L_{K^{(t)}})(D_{h}u) D_{h}u + 2\sigma [(\partial_t + L_{K^{(t)}})(u) u]_{h}\\
&\le 2 \eta^2\Vert \nabla f \Vert_{L^{\infty}((0,1) \times B_1)}  |D_{h}u| + \Vert f \Vert_{L^{\infty}((0,1)\times B_1)} \Vert u \Vert_{L^{\infty}((0,1) \times B_1)}
\end{align*}
and deduce from the parabolic maximum principle:
\begin{align}
\label{eq:parOP-Lipschitz}
\Vert D_{h} u \Vert^2_{L^{\infty}((1/2,1) \times B_{1/2}-|h|)} \le C \left( \Vert f \Vert^2_{L^{\infty}((0,1);C^{0,1}(B_1))} + \Vert u \Vert^2_{L^{\infty}((0,1) \times \R^n)} \right).
\end{align}

\textbf{Step 2:} Second order estimate: \\
Let $\eta, \U$ be as in the proof of \autoref{thm:intro-OP-sc}, and $\bar{\eta}$ as in \autoref{cor:keyest_snd-order_parabolic}. By the same arguments as in the proof of \eqref{eq:application_OP_first-order1} and \eqref{eq:application_OP_snd-order1}, we obtain in $(1/2,1) \times B_{1/4-|h|}$:
\begin{align*}
(\partial_t + L_{K^{(t)}})(D_{h} u) D_{h}u &\le \Vert \nabla f \Vert_{L^{\infty}((0,1) \times B_1)}  |D_{h}u|,\\
(\partial_t + L_{K^{(t)}})(D_{-h}D_{h}u) (D_{-h}D_{h}u)_+ &\le \Vert D^2 f \Vert_{L^{\infty}((0,1) \times B_1)} (D_{-h}D_{h}u)_+.
\end{align*}
Moreover, we have (see \eqref{eq:application_OP_snd-order5})
\begin{align*}
|(\partial_t + L_{K^{(t)}})(D_{h}[(1-\eta)u])| &\le C \Vert u \Vert_{L^{\infty}((0,1) \times \R^n)},\\
|(\partial_t + L_{K^{(t)}})(D_{-h}D_{h}[(1-\eta)u])| &\le C \Vert u \Vert_{L^{\infty}((0,1) \times \R^n)}.
\end{align*}
Therefore, we obtain from \eqref{eq:keyest_snd-order_parabolic_one-sided-diffquot}
\begin{align*}
(\partial_t &+ L_{K^{(t)}}) (\bar{\eta}^2 (D_{-h} D_{h} \U)_+^2 + \sigma [(D_{h} \U)^2]_{-h})\\
&\le 2\bar{\eta}^2 (\partial_t + L_{K^{(t)}})(D_{-h} D_{h} \U)(D_{-h} D_{h} \U)_+ + 2 \sigma [(\partial_t + L_{K^{(t)}})(D_{h} \U) D_{h} \U]_{-h}\\
&\le C \left(\Vert f \Vert_{L^{\infty}((0,1);C^{1,1}(B_1))} + \Vert u \Vert_{L^{\infty}((0,1) \times \R^n)} \right)\left(\bar{\eta}^2(D_{-h}D_{h}u)_+ + [|D_{h} \U|]_{-h} \right),
\end{align*}
and deduce from the parabolic maximum principle:
\begin{align*}
\sup_{[1/2,1] \times B_{1/8}} &(D_{-h} D_{h} u)_+^2 \le \sup_{(1/2,1) \times B_{1/4-|h|}} [\bar{\eta}^2 (D_{-h} D_{h} \U)_+^2 + [(D_{h} \U)^2]_{-h}]\\
&\le C\left(\Vert f \Vert^2_{L^{\infty}((0,1);C^{1,1}(B_1))} + \Vert u \Vert^2_{L^{\infty}((0,1) \times \R^n)} + \Vert D_{h} u \Vert^2_{L^{\infty}([1/2,1) \times B_{1/4 - |h|})} \right).
\end{align*}
A combination of the previous estimate with the Lipschitz estimates \eqref{eq:parOP-Lipschitz} from Step 1 yields the desired semiconvexity estimate \eqref{eq:parabolic_semiconvexity} upon taking the limit $h \to 0$. This proves the claim.
\end{proof}

\subsection{Nonsymmetric operators and drift terms}
\label{subsec:nonsymm}

The goal of this section is to explain how the Bernstein technique developed in this article can be extended to nonsymmetric nonlocal operators. Our study covers nonlocal operators with nonsymmetric jumping kernels that might possess a first order drift term. \\
Recall that so far, we have considered operators of the form
\begin{align*}
Lu(x) = \text{p.v.} \int_{\R^n} (u(x) - u(y))K(x-y) \d y,
\end{align*}
where $K$ was a symmetric jumping kernel, i.e., $K(y) = K(-y)$.  From now on, we will drop this assumption. In order for the nonlocal operator to remain well-defined for smooth functions, if $s \ge 1/2$, we need to slightly adjust the definition of $L$, as follows:
\begin{align}
\label{eq:nonsymm-op}
Lu(x) = 
\begin{cases}
\int_{\R^n} (u(x) - u(y))K(x-y) \d y, ~~ &\text{ if } s \in (0,1/2),\\
\text{p.v.} \int_{\R^n} (u(x) - u(y)) K(x-y) \d y~~ &\text{ if } s = 1/2,\\
\int_{\R^n} (u(x) - u(y) - \nabla u(x) \cdot (x-y))K(x-y) \d y, ~~ &\text{ if } s \in (1/2,1).
\end{cases}
\end{align}
For the $\text{p.v.}$-integral in case $s = 1/2$ to be well-defined, we need to add the following cancellation condition (see \cite{DRSV22}, \cite{FRS23}):
\begin{align}
\label{eq:cancellation}
\int_{B_{R} \setminus B_r} y K(y) \d y = 0 \qquad \forall R > r > 0.
\end{align}
Moreover, we work under the assumption that $K$ satisfies \eqref{eq:coercive} and \eqref{eq:L1} for some constants $0 < \lambda \le \Lambda$. \\
Let us mention the works \cite{ChDa12}, \cite{KiKi12}, where $C^{1,\eps}$-regularity for solutions to fully nonlinear nonlocal equations governed by operators with nonsymmetric jumping kernels has been studied. Moreover, we refer to \cite{DRSV22}, where the interior and boundary regularity theory is developed for nonsymmetric stable operators.

A second way to introduce nonsymmetry to the picture is by considering operators with a drift, i.e., $L + b \cdot \nabla$, where $b \in \R^n$. The regularity theory for nonlocal operators with drifts has recently gained some attention, inspired by the works \cite{CaVa10}, \cite{KiNa10}, \cite{KNV07} on the critical dissipative SQG equation (see also \cite{Sil10}). The case $s \le 1/2$ is of particular interest, since the drift term becomes (super)critical with respect to the scaling of the equation. We refer the interested reader to \cite{Sil11}, \cite{Sil12} for a detailed study of higher regularity properties under the presence of a critical or supercritical drift.

Under these conditions, we are able to establish Bernstein key estimates reminiscent of \eqref{eq:keyest} and \eqref{eq:keyest_pos-part}

\begin{theorem}
\label{thm:keyest_nonsymm}
Let $s\in(0,1)$, $L$ be as in \eqref{eq:nonsymm-op}, $K$ be as before, and $b \in \R^n$ with $|b| \le \Lambda$. Let $\eta \in C^{1,1}(\R^n)$ be such that $\eta^{2s} \in C^{0,1}(\R^n)$ and $\eta \ge 0$.  Then, there exists $\sigma_0 = \sigma_0(n, s, \lambda, \Lambda, \Vert \eta \Vert_{C^{1,1}(\R^n)}, \Vert \eta^{2s} \Vert_{C^{0,1}(\R^n)}) > 0$ such that for every $\sigma \ge \sigma_0$ and $u,v \in C_{loc}^{1+2s+\eps}(\R^n) \cap L^{\infty}(\R^n)$:
\begin{align}
\label{eq:keyest-drift}
(L + b \cdot \nabla)(\eta^2 (\partial_e u)^2 + \sigma u^2) &\le 2 \eta^2 (L + b \cdot \nabla)(\partial_e u) \partial_e u + 2\sigma (L + b \cdot \nabla)(u) u,\\
\label{eq:keyest-drift_one-sided}
(L + b \cdot \nabla)(\eta^2 (\partial_e v)_+^2 + \sigma v^2) &\le 2 \eta^2 (L + b \cdot \nabla)(\partial_e v) (\partial_e v)_+ + 2\sigma (L + b \cdot \nabla)(v) v.
\end{align}
\end{theorem}

\begin{remark}
\begin{itemize}
\item[(i)] By following the same arguments as in Section \ref{sec:diffquot}, it is also possible to derive Bernstein key estimates in terms of difference quotients for $L + b \cdot \nabla$.
\item[(ii)] Recently, in \cite{FRS23}, the regularity theory for the nonlocal obstacle problem has been extended to integro-differential operators of the form \eqref{eq:nonsymm-op} with nonsymmetric kernels satisfying \eqref{eq:coercive}. Moreover, we refer to \cite{FeRo18}, \cite{FRS23} ($s = 1/2$), \cite{PePo15}, \cite{GPP17}, \cite{Kuk21} ($s > 1/2)$, for regularity results on the nonlocal obstacle problem under the presence of a drift term. We expect it to be possible to obtain analogues to \autoref{thm:OP_opt-reg} and \autoref{thm:OP_bdry-reg} also for such generalized problems by combining the ideas in the aforementioned papers with those in Section \ref{sec:obstacle}.
\end{itemize}

\end{remark}

Before we provide a proof of \autoref{thm:keyest_nonsymm}, let us point out that the product rule \autoref{lemma:productrule} and the interpolation \autoref{lemma:interpol} remain valid for nonsymmetric kernels and that the bilinear form $B$ keeps the same shape as in the symmetric case. Moreover, we have a nonsymmetric counterpart of \autoref{lemma:cutoff-est}, which will be applied exactly as before. However the proof of the first estimate changes slightly in case $s \le 1/2$. For completeness, we provide the results and a short proof below:

\begin{lemma}
\label{lemma:cutoff-est_nonsymm}
Let $s\in(0,1)$. Assume 
\begin{align*}
K(y) \le \Lambda |y|^{-n-2s}, ~~ \supp(K) \subset B_{\eps}
\end{align*}
for some $\Lambda > 0$ and $\eps \in (0,1)$. Let $\eta \in C^{1,1}(B_1)$. Then,
\begin{align*}
L(\eta^2)(x) &\le
\begin{cases}
 c_1 \Vert \eta^2 \Vert_{C^{1,1}(B_{\eps}(x))} (\eps + \eta(x)) \eps^{1-2s}, ~~ &\text{ if } s \le 1/2,\\
 c_1 \Vert \eta^2 \Vert_{C^{1,1}(B_{\eps}(x))} \eps^{2-2s}, ~~ &\text{ if } s > 1/2,\\
\end{cases}\\
B(\eta,\eta)(x) &\le c_2 \Vert \nabla \eta \Vert_{L^{\infty}(B_{\eps}(x))}^2 \eps^{2-2s},
\end{align*}
where $c_1, c_2 > 0$ are constants depending only on $n, s, \Lambda$.
\end{lemma}

\begin{proof}
The proof of the second estimate remains exactly the same as in the symmetric case. Also the proof of the first estimate does not change if $s > 1/2$. To prove the first estimate in case $s \le 1/2$, observe that
\begin{align*}
L(\eta^2)(x) &= \int_{B_{\eps}(x)} (\eta^2(x) - \eta^2(y) - \nabla(\eta^2)(x) \cdot (x-y)) K(x-y) \d y\\\
&+ \text{p.v.} \int_{B_{\eps}(x)} \nabla(\eta^2)(x) \cdot (x-y) K(x-y) \d y\\
&\le c_1 \Vert \eta^2 \Vert_{C^{1,1}(B_{\eps}(x))} \eps^{2-2s} + |\nabla(\eta^2)(x)| \left| \text{p.v.} \int_{B_{\eps}} y K(y) \d y \right|.
\end{align*}
In case $s < 1/2$, we can estimate the second term using the pointwise upper bound for $K$ by 
\begin{align*}
|\nabla(\eta^2)(x)| \left| \text{p.v.} \int_{B_{\eps}} y K(y) \d y \right| \le c |\nabla(\eta^2)(x)| \eps^{1-2s} \le c \Vert \eta \Vert_{C^{0,1}(B_{\eps}(x))} \eta(x) \eps^{1-2s}.
\end{align*}
In case $s = 1/2$, the second term can be estimated by 
\begin{align*}
|\nabla(\eta^2)(x)| \left| \text{p.v.} \int_{B_{\eps}} y K(y) \d y \right| \le c |\nabla(\eta^2)(x)| \le c \Vert \eta \Vert_{C^{0,1}(B_{\eps}(x))} \eta(x)
\end{align*}
due to the cancellation condition \eqref{eq:cancellation}. This proves the desired result.
\end{proof}

We are now in a position to prove \autoref{thm:keyest_nonsymm}:

\begin{proof}[Proof of \autoref{thm:keyest_nonsymm}]
Let us first explain how to prove \eqref{eq:keyest-drift}. 
As in the proof of \autoref{thm:keyest}, for any $x$, we split the kernel $K = K_1 + K_2$ according to \autoref{lemma:kernel-decomp} with $\eps = \gamma \eta(x)$, where $\gamma \in (0,1)$ is chosen as in the proof of \autoref{thm:keyest}.
By carefully tracing the arguments in Step 1 and Step 2 of the proof of \autoref{thm:keyest}, it becomes apparent that \autoref{lemma:cutoff-est_nonsymm} can be applied exactly in the same way as \autoref{lemma:cutoff-est} and apart from that no further changes are necessary due to the lack of symmetry, so that we obtain
\begin{align*}
\int_{\R^n} (\eta^2(\cdot) - \eta^2(y)) (\partial_e u(y))^2 K_2(\cdot -y) \d y \le \eta^2 \int_{\R^n} (\partial_e u(\cdot) - \partial_e u(y))^2 K_2(\cdot -y) \d y + \sigma_1  B_{K}(u,u)
\end{align*}
and
\begin{align*}
L_{K_1}(\eta^2) (\partial_e u)^2  - B_{K_1}(\eta^2, (\partial_e u)^2) \le \frac{3}{4}\eta^2 B_{K_1}(\partial_e u , \partial_e u) + \sigma_2 B_K(u,u).
\end{align*}
Let us explain how these two results imply the following estimate 
\begin{align}
\label{eq:keyest-nonsymm_equiv}
L(\eta^2)(\partial_e u)^2 - B(\eta^2 , (\partial_e u)^2) \le \frac{7}{8}\eta^2 B(\partial_e u, \partial_e u) + \sigma B(u,u),
\end{align}
which in turn implies
the key estimate \eqref{eq:keyest-drift} with $b \equiv 0$, namely,
\begin{align}
\label{eq:keyest-nonsymm}
L(\eta^2 (\partial_e u)^2 + \sigma u^2) &\le 2 \eta^2 L(\partial_e u) \partial_e u + 2\sigma L(u) u.
\end{align}

In case $s \le 1/2$, we are already done. In case $s > 1/2$, it remains to estimate 
\begin{align}
\label{eq:help-nonsymm}
\begin{split}
(\partial_e u(x))^2 \int_{\R^n} \nabla \eta^2(x) &\cdot (x-y)K_2(x-y) \d y \\
&\le \frac{1}{8}\eta^2(x) B_{K_1}(\partial_e u , \partial_e u)(x) + \sigma_3 B_K(u,u)(x),
\end{split}
\end{align}
in order to obtain \eqref{eq:keyest-nonsymm_equiv}. The estimate \eqref{eq:help-nonsymm} can be proved observing that by property $(iv)$ of $K_2$ (see \autoref{lemma:kernel-decomp}) and \eqref{eq:coercive}:  
\begin{align*}
\int_{B_1(x)} \nabla \eta^2(x)\cdot (x-y) K_2(x-y) \d y \le 2\Vert \eta \Vert_{C^{0,1}(\R^n)} \eta(x) \int_{B_1 \setminus B_{\eps}} |y|K(y) \d y \le C \eta(x)^{2-2s}.
\end{align*}
Next, we apply the interpolation estimate \autoref{lemma:interpol} with $\delta = \left(\frac{1}{8C}\right)^{\frac{1}{2s}} \eta(x) \wedge \frac{\eps}{2}$ and read off \eqref{eq:help-nonsymm}, as desired.
We have therefore established \eqref{eq:keyest-nonsymm} for all $s \in (0,1)$.

Let us finish the proof of \eqref{eq:keyest-drift} by explaining how to treat the drift term. To this end, let us compute
    \begin{align*}
        (L + b \nabla)(\eta^2 (\partial_e u)^2 + \sigma u^2) &= 2(b\nabla \eta)\eta (\partial_e u)^2\\
        &\quad + 2 \eta^2 (L + b \nabla) (\partial_e u)(\partial_e u) + 2 \sigma (L + b \nabla)(u) u\\
        &\quad + L(\eta^2)(\partial_e u)^2 - B(\eta^2 , (\partial_e u)^2) - \eta^2 B(\partial_e u, \partial_e u) - \sigma B(u,u).
    \end{align*}
Therefore, we need to prove
    \begin{align*}
        2(b\nabla \eta)\eta (\partial_e u)^2 + L(\eta^2)(\partial_e u)^2 - B(\eta^2 , (\partial_e u)^2) \le \eta^2 B(\partial_e u, \partial_e u) + \sigma B(u,u).
    \end{align*}
Since we already know \eqref{eq:keyest-nonsymm_equiv} it remains to establish
    \begin{align}
    \label{eq:drift-help1}
        2(b\nabla \eta)\eta (\partial_e u)^2 \le \frac{1}{8} \eta^2 B_{K_1}(\partial_e u, \partial_e u) + \sigma_4 B_{K_1}(u,u).
    \end{align}
To do so, assume without loss of generality that $b\nabla\eta \neq 0$, otherwise the inequality is trivial. Let us apply the interpolation estimate \autoref{lemma:interpol} with $\delta = \left(\frac{ \eta(x) }{ 16 b \nabla \eta(x)}\right)^{\frac{1}{2s}} \wedge \frac{\eps}{2}$ and proceed as in the proof of \autoref{thm:keyest_parabolic}:
    \begin{align*}
    2(b\nabla \eta)\eta(\partial_e u)^2 &\le 2 (b \nabla \eta) \eta \left(\delta^{2s} B_{K_1}(\partial_e u,\partial_e u) + c \delta^{2s-2} B_{K_1}(u,u)\right) \\
        &\le \frac{1}{8}\eta^2 B_{K_1}(\partial_e u,\partial_e u) + c (b \nabla \eta)^{\frac{1}{s}} \eta^{2 - \frac{1}{s}} B_{K_1}(u,u)\\
        &\le \frac{1}{8}\eta^2 B_{K_1}(\partial_e u,\partial_e u) + c \eta_2^2 B_{K_1}(u,u).
    \end{align*}
Since $(b \nabla\eta)^{\frac{1}{s}} \eta^{2 - \frac{1}{s}} = c(s) (b \nabla(\eta^{2s}))^{\frac{1}{s}}$ is bounded, this yields \eqref{eq:drift-help1}, as desired. This concludes the proof of the Bernstein key estimate \eqref{eq:keyest-drift}.

Finally, in order to prove \eqref{eq:keyest-drift_one-sided}, one makes the same adaptation to the aforementioned proof as in the derivation of \autoref{thm:keyest_pos-part}. Note that the existing of a drift term does not require any changes. However, in order to deal with the lack of symmetry of $K$, it is worth pointing out that in order to adapt \autoref{lemma:interpol_pos-part} to the nonsymmetric case for $s > 1/2$, one has to use that by the regularity of $v$ it holds $\nabla (\partial_e v(x))_- = 0$, whenever $\partial_e v(x) > 0$, in order to obtain
\begin{align*}
\int_{B_{\delta}(x)} & (\partial_e v(y))_- K(x-y) \d y \\
&= -\int_{B_{\delta}(x)} \big((\partial_e v(x))_- - (\partial_e v(y))_- - \nabla (\partial_e v(x))_- \cdot (x-y)\big) K(x-y) \d y\\
&\le - L((\partial_e v)_-)(x),
\end{align*}
which complements the computation in \eqref{eq:help-nonsymm-one-sided}.
\end{proof}

\subsection{General L\'evy kernels}
\label{subsec:Levy}

The goal of this section is to extend the Bernstein technique to nonlocal operators of the form
\begin{align*}
Lu(x) = \text{p.v.} \int_{\R^n} (u(x) - u(y))K(x-y) \d y,
\end{align*}
where $K$ is a kernel that satisfies a comparability condition with respect to a general radial L\'evy kernel but is not necessarily comparable to the fractional Laplacian. Instead of \eqref{eq:coercive}, we assume that 
\begin{align}
\label{eq:g-comp}
\lambda |y|^{-n} g(|y|) \le K(y) \le \Lambda |y|^{-n} g(|y|) \qquad \forall y \in \R^n
\end{align}
for some strictly decreasing function $g : (0,\infty) \to [0,\infty)$ satisfying 
\begin{align}
\label{eq:g-der}
2s_1 r^{-1}g(r) \le |g'(r)| \le 2s_2 r^{-1} g(r) \qquad \forall r > 0
\end{align}
for some $0 < s_1 \le s_2 < 1$. Note that \eqref{eq:g-der} implies the following doubling properties for $g$:
\begin{align*}
g(\lambda r) &\le \lambda^{-2 s_2} g(r) \qquad \forall \lambda < 1,\\
g(\lambda r) &\le \lambda^{-2 s_1} g(r) \qquad \forall \lambda > 1.
\end{align*}
In particular, we have
\begin{align}
\label{eq:Levy-ass-g}
\int_{\R^n} \min \{1 , |y|^2 \} |y|^{-n} g(y) \d y < \infty, \qquad  \lim_{r \searrow 0} g(r) = \infty,
\end{align}
which implies that $K$ is a L\'evy kernel. 

The study of nonlocal operators $L$ satisfying \eqref{eq:Levy-ass-g} whose jumping kernel $K$ is not comparable to the one of the fractional Laplacian are of general interest since they arise as generators of L\'evy processes and in particular of subordinate Brownian motion. In order to study properties of harmonic functions with respect to $L$, it is natural to impose some growth or scaling conditions such as \eqref{eq:g-der} on the kernel, or on the Fourier symbol of $L$.
Let us mention for instance the works \cite{KuRy16, GJZ22}, where gradient estimates for $L$-harmonic functions and \cite{ChKu08,CKS14,BGR14,KSV18,GrSz19,GKK20}, where estimates for fundamental solutions of the associated Cauchy problem are derived, using a probabilistic approach.

We claim that the key estimates \eqref{eq:keyest} and \eqref{eq:keyest_pos-part} remain true in this general setting:

\begin{theorem}
\label{thm:keyest-g}
Let $0 < s_1 \le s_2 < 1$ and $g$ be such that \eqref{eq:g-der} holds true and assume that $K$ satisfies \eqref{eq:g-comp} and \eqref{eq:L1}.
Let $\eta \in C^{1,1}(\R^n)$ be such that $\eta \ge 0$.  Then, there exists $\sigma_0 = \sigma_0(n, s_1, s_2, \Lambda/\lambda, \Vert \eta \Vert_{C^{1,1}(\R^n)}) > 0$ such that for every $\sigma \ge \sigma_0$ and every smooth enough $u,v \in L^{\infty}(\R^n)$
\begin{align}
\label{eq:keyest-g}
\qquad\qquad L\big(\eta^2 (\partial_e u)^2 + \sigma u^2\big) &\le 2 \eta^2 L(\partial_e u) \partial_e u + 2\sigma L(u) u \qquad \textrm{in}\quad \R^n,\\
\label{eq:keyest_pos-part-g}
\qquad\qquad L\big(\eta^2 (\partial_e v)_+^2 + \sigma v^2\big) &\le 2 \eta^2 L(\partial_e v) (\partial_e v)_+ + 2\sigma L(v) v \qquad \textrm{in}\quad \R^n.
\end{align}
\end{theorem}

In our setting, the statement of the interpolation estimate reads as follows:

\begin{lemma}
\label{lemma:interpol-g}
Let $0 < s_1 \le s_2 < 1$ and $\delta \in (0,1)$. Assume that $K$ satisfies for some $0 < \lambda \le \Lambda$:
\begin{align}
\label{eq:interpol-coercive-g}
\lambda |y|^{-n}g(|y|) &\le K(y) \le \Lambda |y|^{-n} g(|y|) ~~ \forall y \in B_{\delta},\\
\label{eq:interpol-C1-g}
|\nabla K(y)| &\le \Lambda |y|^{-1} K(y)~~ \forall y \in B_{\delta},
\end{align}
where $g$ satisfies \eqref{eq:g-der}.
Then, for every $x \in \R^n$ and $u \in C^{0,1}(B_{\delta}(x))$ it holds
\begin{align*}
\big(\partial_e u(x)\big)^2 &\le g(\delta)^{-1} B_K(\partial_e u, \partial_e u)(x) + c\delta^{-2} g(\delta)^{-1} B_K(u,u)(x),\\
(\partial_e v(x))_+^2 &\le g(\delta)^{-1} \left[ B((\partial_e v)_+, (\partial_e v)_+)(x) - L((\partial_e v)_-)(x)(\partial_e v)_+(x)\right] + c\delta^{-2}g(\delta)^{-1} B(v,v)(x),
\end{align*}
where $c = c(n,s_1,s_2,\lambda,\Lambda) > 0$ does not depend on $\delta$.
\end{lemma}

\begin{proof}
We only explain how to construct $K_\delta$. Then, the proof of the estimates goes in the same way as the proof of \autoref{lemma:interpol} and \autoref{lemma:interpol_pos-part}, replacing $\delta^{2s}$ by $g(\delta)^{-1}$.
We define
\begin{align*}
K_{\delta}(y) = \psi(|y|/\delta) K(y) |y|^{\frac{n}{2}+1}g(|y|)^{-\frac{1}{2}}.
\end{align*}
Note that the properties (1), (2), and (3) from the proof of \autoref{lemma:interpol} remain true for this choice of $K_{\delta}$. To see (2), we compute using \eqref{eq:g-der}, \eqref{eq:interpol-coercive-g}, and \eqref{eq:interpol-C1-g}:
\begin{align*}
&|\nabla K_{\delta}(y)|^2 \lesssim (|y|/\delta)^2 |\psi'(|y|/\delta)| K^2(y) |y|^{n} g(|y|)^{-1} \\
&~~\quad + \psi^2(|y|/\delta) \left[ |\nabla K(y)|^2 |y|^{n+2} g(|y|)^{-1} + K^2(y)|y|^{n} g(|y|)^{-1} + K^2(y) |y|^{n+2} |g'(|y|)|^2 g(|y|)^{-3} \right]\\
&~~\lesssim K(y).
\end{align*}
Let us make the following observation, which follows from \eqref{eq:g-der} (resp. its doubling properties):
\begin{align}
\label{eq:int-est-g}
\int_{B_{2r} \setminus B_r} K(y) \d y \asymp \int_r^{2r} t^{-1} g(t) \d t \asymp -\int_{r}^{2r} g'(t) \d t \asymp g(r) - g(2r) \asymp g(r), ~~ \forall r > 0,
\end{align}
Therefore, using again the doubling properties of $g$, as well as that $0 < s_1$, $s_2 < 1$, we have
\begin{align*}
\mu_{K_{\delta}}(B_{\delta}) &= \sum_{k = 0}^{\infty} \int_{B_{\delta 2^{-k}} \setminus B_{\delta 2^{-k-1}}} K_{\delta}(y) \d y \asymp \sum_{k = 0}^{\infty} (\delta 2^{-k})^{\frac{n}{2}+1} g(\delta 2^{-k})^{-1/2} \int_{B_{\delta 2^{-k}} \setminus B_{\delta 2^{-k-1}}} K(y) \d y\\
&\asymp \sum_{k = 0}^{\infty} (\delta 2^{-k})^{\frac{n}{2}+1} g(\delta 2^{-k})^{1/2} \asymp \delta^{\frac{n}{2}+1}g(\delta)^{1/2}.
\end{align*} 
The latter estimate can be seen as a counterpart of property (4) in the proof of \autoref{lemma:interpol}.
\end{proof}

Moreover, we have the following replacement of \autoref{lemma:cutoff-est}:

\begin{lemma}
\label{lemma:cutoff-est-g}
Let $0 < s_1 \le s_2 < 1$ and $K$ be symmetric, with
\begin{align}
\label{eq:cutoff-est-g-ass}
K(y) \le \Lambda |y|^{-n} g(y), ~~ \supp(K) \subset B_{\eps}
\end{align}
for some $\Lambda > 0$ and $\eps \in (0,1)$, where $g$ satisfies \eqref{eq:g-der}. Let $\eta \in C^{1,1}(B_1)$. Then, for any $x \in B_1$
\begin{align*}
L(\eta^2)(x) &\le c_1 \Vert D^2 \eta^2 \Vert_{L^{\infty}(B_{\eps}(x))} \eps^{2} g(\eps),\\
B(\eta,\eta)(x) &\le c_2 \Vert \nabla \eta \Vert_{L^{\infty}(B_{\eps}(x))}^2 \eps^{2} g(\eps),
\end{align*}
where $c_1, c_2 > 0$ are constants depending only on $n, s_1, s_2, \Lambda$.
\end{lemma}

\begin{proof}
The proof follows along the lines of the proof of \autoref{lemma:cutoff-est}, using \eqref{eq:cutoff-est-g-ass} and the doubling properties of $g$ to estimate:
\begin{align*}
\int_{B_{\eps}} |y|^2 K(y) \d y& \lesssim \sum_{k=0}^{\infty} \int_{B_{\eps 2^{-k}} \setminus B_{\eps 2^{-k-1}}} |y|^{2-n} g(|y|) \d y\\
&\lesssim \sum_{k=0}^{\infty} (\eps 2^{-k})^{2} g(\eps 2^{-k}) \lesssim \eps^2 g(\eps) \sum_{k=0}^{\infty} (2^{-k})^{2-2s_2} \lesssim \eps^2 g(\eps).
\end{align*}
\end{proof}

Finally, we can give the:

\begin{proof}[Proof of \autoref{thm:keyest-g}]
	The proof goes in the exact same way as the proofs of \autoref{thm:keyest} and \autoref{thm:keyest_pos-part}. We only need to replace the interpolation estimates \autoref{lemma:interpol} and \autoref{lemma:interpol_pos-part} by \autoref{lemma:interpol-g}, and the cut-off estimate \autoref{lemma:cutoff-est} by \autoref{lemma:cutoff-est-g}.
	Moreover, note that \autoref{lemma:kernel-decomp} remains true in this generalized setup. We apply the interpolation estimate and kernel decomposition with a suitable choice of $\delta$ and $\eps$,  which requires $g$ to be invertible. Note that it is possible to invert $g$ since it is strictly decreasing and by \eqref{eq:Levy-ass-g}, we have that $g(0) = \infty$ and $g(\infty) = 0$.
\end{proof}

\section{Appendix}
\label{sec:appendix}

The goal of this section is to give the proof of \autoref{lemma:CDS}.
Our proof relies on a modification of the ideas from \cite{AbRo20} and \cite{DSV19}. \\
The main auxiliary result in the proof of \autoref{lemma:CDS} is the following variant of Lemma 3.6 from \cite{AbRo20} for the nonlocal obstacle problem:

\begin{lemma}
\label{lemma:AbRo}
Let $s\in(0,1)$, $L \in \mathcal{L}_s(\lambda,\Lambda;1)$, and $\alpha \in (0,s)$. Assume that $u$ with $u \neq 0$ in $B_{1/2}$ is a solution to 
\begin{align*}
\min\{ L u - f , u\} = 0 ~~ \text{ in } B_1,
\end{align*}
where $f \in C^{\beta-2s}(B_1)$ for some $\beta \in (2s , 1+s)$.  Then, $\U = u\mathbbm{1}_{B_2}$ solves the obstacle problem
\begin{align*}
\min \{L \U -\F , \U \} = 0 ~~ \text{ in } B_1,
\end{align*}
for some $\F \in C^{\beta - 2s}(B_1)$ with
\begin{align}
\label{eq:polyRHS-Holder-est}
\Vert \F \Vert_{C^{\beta-2s}(B_1)} \le C \left([f]_{C^{\beta-2s}(B_1)} + \left\Vert \frac{u}{(1 + |\cdot|)^{1+s+\alpha}} \right\Vert_{L^{\infty}(\R^n)} \right) ,
\end{align}
where $C = C(n,s,\lambda,\Lambda) > 0$ is a constant.
\end{lemma}

\begin{proof}[Proof of \autoref{lemma:AbRo}]
We follow the proof of Lemma 3.6 in \cite{AbRo20}:\\
First, we observe that $\U$ satisfies by assumption:
\begin{align*}
\begin{cases}
L \U &= \F ~~ \text{ in } B_1 \cap \{ \U > 0 \},\\
L \U &\ge \F ~~ \text{ in } B_1 \cap \{ \U = 0\}.\\
\end{cases}
\end{align*}
where
\begin{align*}
\F = -L(u \1_{\R^n \setminus B_2}) + f.
\end{align*}
By following the same arguments as in the proof of Lemma 3.6 in \cite{AbRo20}, and observing that we can add and subtract constants to $\F$ without affecting the left hand side of the following estimate, we obtain using \eqref{eq:L1} that for any open set $U \subset B_1$ and $|h| < 2 \dist(U,\R^n \setminus B_1)$:
\begin{align}
\label{eq:F-est}
\Vert |h| D_{h} \F \Vert_{L^{\infty}(U)} \le C \left(\osc_U f + |h|\left\Vert \frac{u}{(1 + |\cdot|)^{1+s+\alpha}} \right\Vert_{L^{\infty}(\R^n)} \right).
\end{align}
Therefore, we can write
\begin{align*}
\F = g + p,
\end{align*}
where $p \in \R$ and $g \in L^{\infty}(B_1)$ satisfies
\begin{align*}
\Vert g \Vert_{L^{\infty}(B_1)} \le  C \left( \osc_{B_1} f + \left\Vert \frac{u}{(1 + |\cdot|)^{1+s+\alpha}} \right\Vert_{L^{\infty}(\R^n)} \right).
\end{align*}
Note that by \autoref{lemma:AbRo2}, $|p|$ satisfies the same upper estimate as $\Vert g \Vert_{L^{\infty}(B_1)}$ and therefore
\begin{align*}
\Vert \F \Vert_{L^{\infty}(B_1)} \le |p| + \Vert g \Vert_{L^{\infty}(B_1)} \le  C \left( \osc_{B_1} f + \left\Vert \frac{u}{(1 + |\cdot|)^{1+s+\alpha}} \right\Vert_{L^{\infty}(\R^n)} \right).
\end{align*}
Moreover, note that since $f \in C^{\beta -2s}(B_1)$, we can deduce from \eqref{eq:F-est} and since $\F = g + p$ for any $x \in B_1$ and $h$ with $|h| < 1 - |x|$: 
\begin{align*}
\frac{g(x) - g(x+h)}{|h|^{2s-\beta}} &\le C \left( |h|^{-2s+\beta} \osc_{B_{|h|}(x)} f + |h|^{1-2s+\beta}\left\Vert \frac{u}{(1 + |\cdot|)^{1+s+\alpha}} \right\Vert_{L^{\infty}(\R^n)} \right) \\
&\le C \left( [f]_{C^{2s-\beta}(B_{|h|}(x))} + \left\Vert \frac{u}{(1 + |\cdot|)^{1+s+\alpha}} \right\Vert_{L^{\infty}(\R^n)} \right)
\end{align*}
and hence
\begin{align*}
[ \F ]_{C^{\beta-2s}(B_1)} = [ g ]_{C^{\beta-2s}(B_1)} \le C \left([f]_{C^{\beta-2s}(B_1)} + \left\Vert \frac{u}{(1 + |\cdot|)^{1+s+\alpha}} \right\Vert_{L^{\infty}(\R^n)} \right) .
\end{align*} 
Altogether, we obtain the desired result.
\end{proof}

\begin{lemma}
\label{lemma:AbRo2}
Let $s\in(0,1)$, $L \in \mathcal{L}_s(\lambda,\Lambda)$. Assume that $u$ with $u \neq 0$ in $B_{1/2}$ is a solution to 
\begin{align*}
\min \{L u - g - p  , u \} &= 0 ~~ \text{ in } B_1,\\
u &\equiv 0 ~~ \text{ in } \R^n \setminus B_2,
\end{align*}
where $g \in C^{\beta-2s}(B_1)$ for some $\beta \in (2s,1+s)$ and $p \in \R$.
Then, it holds
\begin{align}
\label{eq:AbRo2-estimate}
|p| \le C \left(\Vert g \Vert_{L^{\infty}(B_1)} + \Vert u \Vert_{L^{\infty}(\R^n)} \right),
\end{align}
where $C = C(n,s,\lambda,\Lambda) > 0$ is a constant.
\end{lemma}

\begin{proof}
Assume that \eqref{eq:AbRo2-estimate} does not hold true. Then by contradiction there exist sequences $(L_k) \subset \mathcal{L}_s(\lambda,\Lambda)$, $(g_k) \subset L^{\infty}(B_1)$, $(p_k) \subset \R$, and $(u_k) \subset L^{\infty}(\R^n)$ with $u_k \neq 0$ in $B_{1/2}$, and $u_k \equiv 0$ in $\R^n \setminus B_2$, such that
\begin{align}
\label{eq:AbRo2help1}
\min \{L_k u_k - g_k - p_k  , u_k \}& = 0 ~~ \text{ in } B_1,\\
\nonumber
\Vert u_k \Vert_{L^{\infty}(\R^n)} &\to 0,\\
\nonumber
\Vert g_k \Vert_{L^{\infty}(B_1)} &\to 0,\\
\nonumber
|p_k| &= 1 ~~ \forall k \in \N.
\end{align}
The last property can be assumed by normalizing the sequence $(p_k)$. Then, the second and third property follow immediately from \eqref{eq:AbRo2help1} and since we assumed that the sequence violates  \eqref{eq:AbRo2-estimate}.
We can extract subsequences such that $L_{k_m} \to L$ weakly (using a similar argument as in \cite[Lemma 3.1]{RoSe16b}), $p_{k_m} \to p$ with $|p|=1$, $g_{k_m} \to g = 0$ in $L^{\infty}(B_1)$, $u_{k_m} \to u = 0$ in $L^{\infty}(\R^n)$. Moreover, since by \cite[Theorem 5.1]{CDS17} we get that $\Vert u_k \Vert_{C^{2s+\eps}(B_{1/2})} \le C$, it holds $u_k \to u$ in $C^{2s+\eps}(B_{1/2})$ (up to a subsequence) by Arzel\`a-Ascoli, and therefore $|L_k u_k| = |L_k u_k - L u| \to 0$ locally uniformly in $B_1$.\\
Consequently,
\begin{align*}
\min \{L u - g - p  , u \}& = 0 ~~ \text{ in } B_1.
\end{align*}
In particular, $\min \{-p, 0 \} = 0$ in $B_1$. This is a contradiction if $p = 1 > 0$. If $p = -1 < 0$, then, there must be $k \in \N$ such that $p_k < - \frac{3}{4}$, $\Vert g_k \Vert_{L^{\infty}(B_1)} < \frac{1}{8}$, and $\Vert L_k u_k \Vert_{L^{\infty}(B_{1/2})} < \frac{1}{8}$. However, since $u_k \neq 0$ in $B_{1/2}$ by assumption, this contradicts \eqref{eq:AbRo2help1}.
Therefore, \eqref{eq:AbRo2-estimate} holds true.
\end{proof}

We are now in a position to give the proof of \autoref{lemma:CDS}.

\begin{proof}[Proof of \autoref{lemma:CDS}]
Note that we can assume $u \neq 0$ in $B_{1/2}$ without loss of generality, since otherwise there is nothing to prove. We define $\U = u \1_{B_2}$ and deduce from \autoref{lemma:AbRo} that $\U$ solves
\begin{align*}
\min \{ L \U -\F , \U \} = 0 ~~ \text{ in } B_1.
\end{align*}
By application of Theorem 5.1 in \cite{CDS17} together with  \eqref{eq:polyRHS-Holder-est} we obtain
\begin{align*}
\Vert u \Vert_{C^{\max\{2s+\eps,1+\eps\}}(B_{1/2})} &= \Vert \U \Vert_{C^{\max\{2s+\eps,1+\eps\}}(B_{1/2})} \\
&\le C \left(\Vert \F\Vert_{C^{\beta-2s}(B_1)} + \left\Vert \frac{\U}{(1+|\cdot|)^{n + 2s}} \right\Vert_{L^{1}(\R^n)} \right)\\
&\le C \left([f]_{C^{\beta-2s}(B_1)} + \left\Vert \frac{u}{(1 + |\cdot|)^{1+s+\alpha}} \right\Vert_{L^{\infty}(\R^n)} \right).
\end{align*}
This proves the desired result.
\end{proof}

\enlargethispage{1.5cm}

\end{document}